\documentclass{article}
\usepackage{alphabeta,amsfonts,amsmath,amssymb,amsthm,authblk}
\usepackage[english]{babel}
\usepackage{booktabs}
\usepackage[makeroom]{cancel}
\usepackage{comment}
\usepackage{enumerate}
\usepackage[shortlabels]{enumitem}
\usepackage[margin = 1 in]{geometry}
\usepackage{graphicx}
\usepackage[usenames,dvipsnames]{xcolor}
\usepackage[linktocpage,colorlinks=true,linkcolor=Fuchsia!70!Magenta,urlcolor=blue!60!cyan,citecolor=Green!70!PineGreen,
	hypertexnames=false]{hyperref}
\usepackage[utf8]{inputenc}
\usepackage{leftidx}
\usepackage{mathrsfs,mathtools}
\usepackage{relsize}
\usepackage{scalerel,stackengine}
\usepackage{tikz,tikz-cd}
\usepackage[titles]{tocloft}
\usepackage{xspace}

\newcommand{\N}{\mathbb{N}}
\newcommand{\Z}{\mathbb{Z}}
\newcommand{\Q}{\mathbb{Q}}
\newcommand{\R}{\mathbb{R}}
\newcommand{\C}{\mathbb{C}}
\newcommand{\F}{\mathbb{F}}

\newcommand{\bA}{\boldsymbol{A}}
\newcommand{\cA}{\mathcal{A}}

\newcommand{\cB}{\mathcal{B}}
\newcommand{\sB}{\mathscr{B}}

\newcommand{\cC}{\mathcal{C}}
\newcommand{\cF}{\mathcal{F}}
\newcommand{\sF}{\mathscr{F}}
\newcommand{\sG}{\mathscr{G}}
\newcommand{\sH}{\mathscr{H}}
\newcommand{\cM}{\mathcal{M}}
\newcommand{\cS}{\mathcal{S}}

\newcommand{\e}{\varepsilon}
\newcommand{\Om}{\Omega} 
\newcommand{\om}{\omega}
\newcommand{\blambda}{\boldsymbol{\lambda}}

\newcommand{\MnC}{\operatorname{M}_n(\C)}

\DeclareMathOperator{\spn}{span}

\newcommand{\la}{\langle}
\newcommand{\ra}{\rangle}

\DeclareMathOperator{\sa}{sa}
\newcommand{\vertiii}[1]{{\left\vert\kern-0.25ex\left\vert\kern-0.25ex\left\vert #1 
    \right\vert\kern-0.25ex\right\vert\kern-0.25ex\right\vert}}
\DeclareMathOperator{\supp}{supp}

\DeclareMathOperator{\loc}{loc}

\newcommand{\iotimes}{\hat{\otimes}_i}

\DeclareMathOperator{\HS}{HS}
\newcommand{\wh}{\widehat}
\stackMath
\newcommand\wch[1]{%
\savestack{\tmpbox}{\stretchto{%
  \scaleto{%
    \scalerel*[\widthof{\ensuremath{#1}}]{\kern-.6pt\bigwedge\kern-.6pt}%
    {\rule[-\textheight/2]{1ex}{\textheight}}
  }{\textheight}%
}{0.5ex}}%
\stackon[1pt]{#1}{\scalebox{-1}{\tmpbox}}%
}

\newcommand\numberthis{\addtocounter{equation}{1}\tag{\theequation}}
\makeatletter
\newcommand{\oset}[3][0ex]{%
  \mathrel{\mathop{#3}\limits^{
    \vbox to#1{\kern-2\ex@
    \hbox{$\scriptstyle#2$}\vss}}}}
\makeatother

\theoremstyle{plain}
\newtheorem{prop}{Proposition}[subsection]
\theoremstyle{plain}

\theoremstyle{plain}
\newtheorem{lem}[prop]{Lemma}
\theoremstyle{plain}

\theoremstyle{plain}
\newtheorem{thm}[prop]{Theorem}
\theoremstyle{plain}

\theoremstyle{plain}
\newtheorem{cor}[prop]{Corollary}
\theoremstyle{plain}
\newtheorem{conj}[prop]{Conjecture}
\theoremstyle{plain}

\theoremstyle{plain}

\theoremstyle{plain}

\theoremstyle{definition}
\newtheorem{defi}[prop]{Definition}
\theoremstyle{definition}

\theoremstyle{definition}
\newtheorem{nota}[prop]{Notation}
\theoremstyle{definition}
\newtheorem*{notan}{Notation}
\theoremstyle{definition}

\theoremstyle{definition}

\theoremstyle{definition}
\newtheorem{ex}[prop]{Example}
\theoremstyle{definition}

\theoremstyle{definition}

\theoremstyle{definition}

\theoremstyle{definition}
\newtheorem{rem}[prop]{Remark}
\theoremstyle{definition}

\theoremstyle{definition}

\theoremstyle{definition}

\theoremstyle{definition}
\newtheorem*{ack}{Acknowledgements}

\makeatletter
\renewenvironment{proof}[1][\proofname]{%
  \par\pushQED{\qed}\normalfont%
  \topsep6\p@\@plus6\p@\relax
  \trivlist\item[\hskip\labelsep\bfseries#1\@addpunct{.}]%
  \ignorespaces
}{%
  \popQED\endtrivlist\@endpefalse
}
\makeatother
\begin{document}

\title{Noncommutative \texorpdfstring{$C^k$}{} functions and \\ Fr\'{e}chet derivatives of operator functions\vspace{-2.8mm}}
\author{Evangelos A. Nikitopoulos\thanks{Supported by NSF grant DGE 2038238 and partially supported by NSF grants DMS 1253402 and DMS 1800733}\vspace{-2.5mm}}
\affil{Department of Mathematics, University of California San Diego\protect\\
\noindent 9500 Gilman Drive, La Jolla, CA 92093-0112 (USA)\protect\\
Email: {\tt \href{mailto:enikitop@ucsd.edu}{enikitop@ucsd.edu}}}
\date{\vspace{-6.5ex}}

\maketitle

\begin{abstract}
Fix a unital $C^*$-algebra $\mathcal{A}$, and write $\mathcal{A}_{\mathrm{sa}}$ for the set of self-adjoint elements of $\mathcal{A}$.
Also, if $f \colon \mathbb{R} \to \mathbb{C}$ is a continuous function, then write $f_{\mathsmaller{\mathcal{A}}} \colon \mathcal{A}_{\mathrm{sa}} \to \mathcal{A}$ for the \textit{operator function} $a \mapsto f(a)$ defined via functional calculus.
In this paper, we introduce and study a space $NC^k(\mathbb{R})$ of $C^k$ functions $f \colon \mathbb{R} \to \mathbb{C}$ such that, no matter the choice of $\mathcal{A}$, the operator function $f_{\mathsmaller{\mathcal{A}}} \colon \mathcal{A}_{\mathrm{sa}} \to \mathcal{A}$ is $k$-times continuously Fr\'{e}chet differentiable.
In other words, if $f \in NC^k(\R)$, then $f$ ``lifts" to a $C^k$ map $f_{\mathsmaller{\mathcal{A}}} \colon \mathcal{A}_{\mathrm{sa}} \to \mathcal{A}$, for any (possibly noncommutative) unital $C^*$-algebra $\mathcal{A}$.
For this reason, we call $NC^k(\mathbb{R})$ the space of \textit{noncommutative $C^k$ functions}.
Our proof that $f_{\mathsmaller{\mathcal{A}}} \in C^k(\mathcal{A}_{\mathrm{sa}};\mathcal{A})$, which requires only knowledge of the Fr\'{e}chet derivatives of polynomials and operator norm estimates for ``multiple operator integrals" (MOIs), is more elementary than the standard approach;
nevertheless, $NC^k(\mathbb{R})$ contains all functions for which comparable results are known.
Specifically, we prove that $NC^k(\mathbb{R})$ contains the homogeneous Besov space $\dot{B}_1^{k,\infty}(\mathbb{R})$ and the H\"{o}lder space $C_{\mathrm{loc}}^{k,\varepsilon}(\mathbb{R})$.
We highlight, however, that the results in this paper are the first of their type to be proven for arbitrary unital $C^*$-algebras, and that the extension to such a general setting makes use of the author's recent resolution of certain ``separability issues" with the definition of MOIs.
Finally, we prove by exhibiting specific examples that $W_k(\mathbb{R})_{\mathrm{loc}} \subsetneq NC^k(\mathbb{R}) \subsetneq C^k(\mathbb{R})$, where $W_k(\mathbb{R})_{\mathrm{loc}}$ is the ``localized" $k^{\text{th}}$ Wiener space.\vspace{0.5mm}

\noindent \textbf{Keywords:}
functional calculus, $C^*$-algebra, Fr\'{e}chet derivative, multiple operator integral\vspace{0.5mm}

\noindent \textbf{MSC (2010):}
47A60, 47A56, 26E15 (primary);
46E35, 46L89 (secondary)
\end{abstract}
\vspace{-3.75mm}

\tableofcontents
\clearpage

\section{Introduction}\label{sec.intro}

\begin{notan}
Let $V$ and $W$ be normed vector spaces over $\F \in \{\R,\C\}$, $\cA$ be a unital $C^*$-algebra, $\Om$ be a topological space, $(\Xi,\sG)$ be a measurable space, $H$ be a complex Hilbert space, and $P \colon \sG \to B(H)$ be a projection-valued measure (Definition 1.1 in Chapter IX of \cite{conwayfunc}).
\begin{enumerate}[label=(\alph*),leftmargin=2\parindent]
    \item $B(V;W)$ is the space of bounded linear maps $V \to W$ with operator norm $\|\cdot\| = \|\cdot\|_{V \to W}$, and $B(V) = B(V;V)$.
    Also, $C(V;W) = C^0(V;W)$ is the space of continuous maps $V \to W$.
    Finally, if $k \in \N$, then $C^k(V;W)$ is the space of $k$-times continuously Fr\'{e}chet differentiable (Definition \ref{def.frechder}) maps $V \to W$.
    When $W=\C$, we shall suppress $W$ from this notation.
    \item $\cA_{\nu} = \{a \in \cA : a^*a=aa^*\}$ and $\cA_{\sa} = \{a \in \cA : a^*=a\}$.
    Also, if $f \in C(\R)$, then $f_{\mathsmaller{\cA}} \colon \cA_{\sa} \to \cA$ is the map defined via functional calculus by $\cA_{\sa} \ni a \mapsto f(a) \in \cA$.
    \item $\cB_{\Om}$ is the Borel $\sigma$-algebra on $\Om$.
    \item If $\varphi \colon \Xi \to \C$ is a function, then $\|\varphi\|_{\ell^{\infty}(\Xi)} = \sup_{\xi \in \Xi} |\varphi(\xi)| \in [0,\infty]$.
    Also, $\ell^{\infty}(\Xi,\sG)$ is the space of $(\sG,\cB_{\C})$-measurable functions $\varphi \colon \Xi \to \C$ such that $\|\varphi\|_{\ell^{\infty}(\Xi)} < \infty$.
    \item If $\varphi \in \ell^{\infty}(\Xi,\sG)$, then $P(\varphi) = \int_{\Xi} \varphi \, dP = \int_{\Xi} \varphi(\xi)\,P(d\xi) \in B(H)$ is the integral of $\varphi$ with respect to $P$.
    Also, if $a \in B(H)_{\nu}$, then $P^a \colon \cB_{\sigma(a)} \to B(H)$ is the projection-valued spectral measure of $a$ determined by $a = \int_{\sigma(a)}\lambda\,P^a(d\lambda)$.
    (Recall the existence and uniqueness of $P^a$ is precisely the Spectral Theorem.)
\end{enumerate}
Please see \cite{birmansolomyak5,conwayfunc,conwayop,sakai} for requisite material on operator theory --- e.g., $C^*$-algebras, von Neumann algebras, spectral theory, functional calculus, and projection-valued measure theory.
\end{notan}

The problem of computing higher derivatives of the operator function $f_{\mathsmaller{B(H)}} \colon B(H)_{\sa} \to B(H)$ has been studied extensively and is the reason multiple operator integrals (MOIs) were invented.
The current standard approach to this problem, which makes use of ``higher-order perturbation formulas," is tailored to differentiating the map $b \mapsto f(a+b)-f(a)$, where $a$ is an \textit{unbounded} self-adjoint operator on $H$.
Recently, however, there have been interesting applications --- e.g., to free stochastic calculus \cite{nikitopoulosIto} --- for which it is only necessary to consider bounded operators.
In this paper, we show that when one restricts to bounded operators, the classical approach of ``approximation by polynomials" can be adapted to achieve even better, more general results than the standard approach.

\subsection{Motivation}\label{sec.motiv}

Let $\cA$ be a unital $C^*$-algebra and $f \colon \R \to \C$ be a function.
We claim that if $f$ is continuous, then the \textbf{operator function} $f_{\mathsmaller{\cA}} \colon \cA_{\sa} \to \cA$ is continuous as well. Indeed, this is easy to see if $f$ is a polynomial.
For arbitrary $f \in C(\R)$, let $(q_n)_{n \in \N}$ be a sequence of polynomials such that $q_n \to f$ uniformly on compact sets as $n \to \infty$.
Such a sequence exists by the Weierstrass Approximation Theorem.
Then $(q_n)_{\mathsmaller{\cA}} \to f_{\mathsmaller{\cA}}$ uniformly on bounded sets as $n \to \infty$, from which it follows that $f_{\mathsmaller{\cA}}$ is continuous.
It is therefore natural to wonder whether $f \in C^k(\R)$ implies $f_{\mathsmaller{\cA}} \in C^k(\cA_{\sa};\cA)$ whenever $k \in \N$.
It turns out that this is not generally true.
For example, if $\dot{B}_q^{s,p}(\R)$ is the homogeneous $(s,p,q)$-Besov space (Definition \ref{def.Besov}), then results of V. V. Peller \cite{peller0} imply that if $f \in C^1(\R)$ is such that $f_{\mathsmaller{B(\ell^{\text{\scalebox{0.7}{$2$}}}(\N))}} \colon B(\ell^2(\N))_{\sa} \to B(\ell^2(\N))$ is Fr\'{e}chet differentiable, then $f$ must belong locally to the space $\dot{B}_1^{1,1}(\R)$.
(Please  see Sections 1.2 and 1.5 of \cite{aleksandrovOL} for more information about once differentiability of operator functions.)

To elucidate the difficulties with differentiating operator functions and to motivate our results, we begin by studying matrix functions.
In other words, we consider the case
\[
\cA = \MnC = \{n \times n \text{ complex matrices}\}, \,n \in \N,
\]
identified with $B(\C^n)$ in the usual way.
Under this identification, if $A \in \MnC_{\nu}$, $\lambda \in \sigma(A)$, and $g \colon \sigma(A) \to \C$ is a function, then $P_{\lambda}^A \coloneqq P^A(\{\lambda\})$ is (the standard representation of) the orthogonal projection onto the $\lambda$-eigenspace of $A$ and
\[
g(A) = \int_{\sigma(A)} g \,dP^A = \sum_{\lambda \in \sigma(A)}g(\lambda)\,P_{\lambda}^A. \numberthis\label{eq.fdfunccalc}
\]
Now, fix $A,B \in \MnC_{\sa}$.
If $f \in C^1(\R)$, then
\begin{align*}
    f(A)-f(B) & = \sum_{\lambda \in \sigma(A)}f(\lambda)\,P_{\lambda}^A - \sum_{\mu \in \sigma(B)} f(\mu)\,P_{\mu}^B = \sum_{\lambda \in \sigma(A)} \sum_{\mu \in \sigma(B)}(f(\lambda)-f(\mu))P_{\lambda}^AP_{\mu}^B \\
    & = \sum_{\lambda \in \sigma(A)}\sum_{\mu \in \sigma(B)}\frac{f(\lambda)-f(\mu)}{\lambda-\mu} P_{\lambda}^A(\lambda-\mu)P_{\mu}^B = \sum_{\lambda \in \sigma(A)}\sum_{\mu \in \sigma(B)} f^{[1]}(\lambda,\mu)\,P_{\lambda}^A(A-B)P_{\mu}^B,
\end{align*}
where \hspace{-0.2mm}$f^{[1]}(\lambda,\lambda) \hspace{-0.2mm}\coloneqq\hspace{-0.2mm} f'(\lambda)$\hspace{-0.2mm} so \hspace{-0.2mm}that\hspace{-0.2mm} $f^{[1]} \hspace{-0.2mm}\in\hspace{-0.2mm} C(\R^2)$.
\hspace{-0.2mm}In\hspace{-0.2mm} the \hspace{-0.2mm}first\hspace{-0.2mm} line, \hspace{-0.2mm}we\hspace{-0.2mm} used \hspace{-0.2mm}that\hspace{-0.2mm} $\sum_{\lambda \in \sigma(A)} \hspace{-0.2mm}P_{\lambda}^A \hspace{-0.2mm}=\hspace{-0.2mm} I_n \hspace{-0.2mm}=\hspace{-0.2mm} \sum_{\mu \in \sigma(B)} \hspace{-0.2mm}P_{\mu}^B$.
In the second line, we used that $P_{\lambda}^AP_{\mu}^A = \delta_{\lambda\mu}P_{\mu}^A$ and $A = \sum_{\lambda \in \sigma(A)}\lambda \,P_{\lambda}^A$ (and the same properties for $B$).
If we set the notation
\[
\big(I^{A,B}\varphi\big)[C] \coloneqq \sum_{\lambda \in \sigma(A)}\sum_{\mu \in \sigma(B)} \varphi(\lambda,\mu)\,P_{\lambda}^ACP_{\mu}^B
\]
for $\varphi \colon \sigma(A) \times \sigma(B) \to \C$ and $C \in \MnC$, then we have just shown that
\[
f(A)-f(B) = \big(I^{A,B}f^{[1]}\big)[A-B]. \numberthis\label{eq.pertform}
\]
This equation is called a \textbf{perturbation formula}.
Next, let us take for granted that if $\varphi \in C(\R^2)$ and $C \in \MnC$, then the map $\MnC_{\sa}^2 \ni (A,B) \mapsto (I^{A,B}\varphi)[C] \in \MnC$ is continuous, as can be seen by approximating $\varphi$ by two-variable polynomials.
Since $f^{[1]} \in C(\R^2)$, it follows that
\begin{align*}
    \partial_Bf_{\mathsmaller{\MnC}}(A) & = \lim_{t \to 0}\frac{1}{t}(f(A+tB)-f(A)) = \lim_{t \to 0}\big(I^{A+tB,A}f^{[1]}\big)[B] \\
    & = \big(I^{A,A}f^{[1]}\big)[B] = \sum_{\lambda,\mu \in \sigma(A)}f^{[1]}(\lambda,\mu)\,P_{\lambda}^ABP_{\mu}^A.
\end{align*}
A higher-order version of $f^{[1]}$ appears in a similar way in the higher derivatives of $f_{\mathsmaller{\MnC}}$.
We state the result presently.
For $k \in \N_0$ and $f \in C^k(\R)$, we define a function $f^{[k]} \in C(\R^{k+1})$, called the $\boldsymbol{k^{\textbf{th}}}$ \textbf{divided difference} of $f$, recursively by $f^{[0]} \coloneqq f$ and
\[
f^{[k]}(\blambda) \coloneqq \frac{f^{[k-1]}(\lambda_1,\ldots,\lambda_k)-f^{[k-1]}(\lambda_1,\ldots,\lambda_{k-1},\lambda_{k+1})}{\lambda_k-\lambda_{k+1}}
\]
for $k \in \N$ and $\blambda = (\lambda_1,\ldots,\lambda_{k+1}) \in \R^{k+1}$ with $\lambda_k \neq \lambda_{k+1}$.
(Please see Section \ref{sec.divdiff} for details.)

\begin{thm}[Derivatives of Matrix Functions]\label{thm.diffmatfunc}
Fix natural numbers $n,k \in \N$.
If $f \in C^k(\R)$, then $f_{\mathsmaller{\MnC}} \in C^k(\MnC_{\sa};\MnC)$ and
\[
\partial_{B_1}\cdots\partial_{B_k}f_{\mathsmaller{\MnC}}(A) = \sum_{\pi \in S_k}\sum_{\blambda \in \sigma(A)^{k+1}} f^{[k]}(\blambda) \, P_{\lambda_1}^AB_{\pi(1)}\cdots P_{\lambda_k}^AB_{\pi(k)}P_{\lambda_{k+1}}^A, \numberthis\label{eq.matfunccalcderform}
\]
for all $A,B_1,\ldots,B_k \in \MnC_{\sa}$, where $S_k$ is the symmetric group on $k$ letters.
\end{thm}

\begin{rem}
This result is due essentially to Yu. L. Daletskii and S. G. Krein \cite{daletskiikrein} (in 1956), though it was proven in approximately the above form by F. Hiai as Theorem 2.3.1 in \cite{hiai}.
\end{rem}

One way to prove Theorem \ref{thm.diffmatfunc} is to use induction and higher-order perturbation formulas --- i.e., versions of Eq. \eqref{eq.pertform} involving $f^{[k]}$ and $f^{[k+1]}$.
(For this kind of argument, please see the proof of Theorem 5.3.2 in \cite{skripka}.)
This is currently the standard approach to proving such results, since it can be adapted to differentiating operator functions at unbounded operators.
Please see, e.g., \cite{azamovetal,coine,depagtersukochev,lemerdymcdonald,lemerdyskripka,nikitopoulosOpDer,peller1,peller2}.
The classical approach is to use a polynomial approximation argument: first establish Eq. \eqref{eq.matfunccalcderform} when $f$ is a polynomial, and then approximate a general $f \in C^k(\R)$ by polynomials.
Since the details of this argument are important motivation for our main definitions and results, we go through them in Section \ref{sec.dermatfunc}.

Next, suppose $H$ is a complex Hilbert space (now possibly infinite-dimensional).
In view of Eqs. \eqref{eq.fdfunccalc} and \eqref{eq.matfunccalcderform}, it seems as though one should have
\[
\partial_{b_1}\cdots\partial_{b_k}f_{\mathsmaller{B(H)}}(a) = \sum_{\pi \in S_k}\underbrace{\int_{\sigma(a)}\cdots\int_{\sigma(a)}}_{k+1 \, \mathrm{times}}f^{[k]}(\blambda) \, P^a(d\lambda_1) \, b_{\pi(1)} \cdots P^a(d\lambda_k) \, b_{\pi(k)} \, P^a(d\lambda_{k+1}), \numberthis\label{eq.formalopderfunc}
\pagebreak
\]
for all $a,b_1,\ldots,b_k \in B(H)_{\sa}$.
However, in standard projection-valued measure theory, one only defines integrals of scalar-valued functions against projection-valued measures.
It is therefore unclear how even to interpret --- let alone prove --- Eq. \eqref{eq.formalopderfunc}.
In their seminal paper \cite{daletskiikrein}, Daletskii and Krein did so by using a Riemann--Stieltjes-type construction to define $\int_{\sigma(a)} F(\lambda) \, P^a(d\lambda)$ for certain operator-valued functions $F \colon [s,t] \to B(H)$ with $\sigma(a) \subseteq [s,t]$.
This approach, which requires rather stringent regularity conditions on $F$, allowed them to interpret the right hand side of Eq. \eqref{eq.formalopderfunc} as an iterated operator-valued integral --- in other words, a \textit{multiple operator integral} --- when $f \in C^{2k}(\R)$.
Moreover, they used the aforementioned polynomial approximation argument to prove Eq. \eqref{eq.formalopderfunc} (with $b_1=\cdots=b_k$) when $f \in C^{2k}(\R)$.
Hence the term ``classical approach."

In the end, it turns out that the assumption $f \in C^{2k}(\R)$ is far too strong and that the key to relaxing it is finding a different way to interpret the multiple operator integral (MOI) on the right hand side of Eq. \eqref{eq.formalopderfunc}.
For our purposes, the right way to do so is to use the ``separation of variables approach" developed originally for separable $H$ in \cite{azamovetal,peller1} and extended by the present author to not-necessarily-separable $H$ in \cite{nikitopoulosMOI}.
Using this approach, if $\cM \subseteq B(H)$ is a von Neumann algebra, then one is able to make sense of
\[
\int_{\sigma(a_{k+1})}\cdots\int_{\sigma(a_1)}\varphi(\blambda) \,P^{a_1}(d\lambda_1)\,b_1\cdots P^{a_k}(d\lambda_k) \,b_k\,P^{a_{k+1}}(d\lambda_{k+1}) \in \cM \numberthis\label{eq.formalMOI}
\]
for $(a_1,\ldots,a_{k+1}) \in \cM_{\nu}^{k+1}$, $(b_1,\ldots,b_k) \in \cM^k$, and functions $\varphi \colon \sigma(a_1) \times \cdots \times \sigma(a_{k+1}) \to \C$ belonging to the $\ell^{\infty}$\textit{-integral projective tensor product} $\ell^{\infty}(\sigma(a_1),\cB_{\sigma(a_1)}) \iotimes \cdots \iotimes \ell^{\infty}(\sigma(a_{k+1}),\cB_{\sigma(a_{k+1})})$ (Definition \ref{def.babyIPTP}).
We present a beginner-friendly review of how to do this in Section \ref{sec.MOIs}.
Henceforth, any MOI we write down is to be interpreted in accordance with Section \ref{sec.MOIs}.
For much more information about MOIs and their applications, please see A. Skripka and A. Tomskova's book \cite{skripka}.

We now state the weakest known regularity assumption on $f$ that guarantees $k$-times differentiability (in the operator norm) of $f_{\mathsmaller{B(H)}}$.
As was mentioned earlier, it is not common to differentiate only at bounded operators, so the best available result is about derivatives of operator functions at unbounded operators.

\begin{thm}[Peller \cite{peller1}, Nikitopoulos \cite{nikitopoulosOpDer}]\label{thm.PellerNik}
Fix a complex Hilbert space $H$ (not necessarily separable), a von Neumann algebra $\cM \subseteq B(H)$, and a (possibly unbounded) self-adjoint operator $a$ on $H$ affiliated with $\cM$.
If $f \in \dot{B}_1^{1,\infty}(\R) \cap \dot{B}_1^{k,\infty}(\R)$ is such that $f^{(k)}$ is bounded, then
\begin{enumerate}[label=(\roman*),font=\normalfont,leftmargin=2\parindent]
    \item $f_a(b) \coloneqq f(a+b)-f(a) \in \cM$ whenever $b \in \cM_{\sa}$,
    \item $f^{[k]}\big|_{\sigma(a_1) \times \cdots \times \sigma(a_{k+1})} \,\in\, \ell^{\infty}(\sigma(a_1),\cB_{\sigma(a_1)}) \,\iotimes\, \cdots \,\iotimes\, \ell^{\infty}(\sigma(a_{k+1}),\cB_{\sigma(a_{k+1})})$ whenever $a_1,\ldots,a_{k+1}$ are (possibly unbounded) self-adjoint operators on $H$, and\label{item.f[k]Peller}
    \item $f_a \in C^k(\cM_{\sa};\cM)$ with
    \[
    \partial_{b_1}\cdots\partial_{b_k}f_a(0) = \sum_{\pi \in S_k}\underbrace{\int_{\sigma(a)}\cdots\int_{\sigma(a)}}_{k+1 \, \mathrm{times}}f^{[k]}(\blambda) \, P^a(d\lambda_1) \, b_{\pi(1)} \cdots P^a(d\lambda_k) \, b_{\pi(k)} \, P^a(d\lambda_{k+1}),
    \]
    for all $b_1,\ldots,b_k \in \cM_{\sa}$.\label{item.PellerNikder}
\end{enumerate}
In particular, $f_{\mathsmaller{\cM}} \in C^k(\cM_{\sa};\cM)$, and Eq. \eqref{eq.formalopderfunc} holds for such $f$.
\end{thm}
\begin{rem}
In \cite{peller1}, Peller proved (something stronger than) item \ref{item.f[k]Peller} and item \ref{item.PellerNikder} (with $b_1=\cdots=b_k$) in the case $\cM = B(H)$ with $H$ separable.
In \cite{nikitopoulosOpDer}, the present author proved (something stronger than) the general case of item \ref{item.PellerNikder}.
Both papers use perturbation formulas to compute the higher derivatives of $f_a$.
\end{rem}

The goal of this paper may now be stated slightly more precisely as follows:
to use modern MOI theory to adapt the classical approach (polynomial approximation) to differentiating operator functions in such a way that recovers the result above in the case when $a$ is bounded.
As we see shortly, we do better.
Specifically, we are able to relax the regularity assumptions on $f$ slightly (when $a$ is bounded) and to generalize the setting from von Neumann algebras $\cM$ to arbitrary unital $C^*$-algebras $\cA$.
To do the latter generalization, we make use of MOIs on non-separable Hilbert spaces --- the subject of \cite{nikitopoulosMOI} --- in an interesting way.
Please see Theorems \ref{thm.main1}\ref{item.Besov} and \ref{thm.main2} and the discussion following Theorem \ref{thm.main2} for details.

\subsection{Main results}\label{sec.mainres}

For the duration of this section, fix $k \in \N$.
For readers who skipped the previous section, we recall now that
\begin{itemize}[leftmargin=2\parindent]
    \item If $k \in \N$ and $f \in C^k(\R)$, then $f^{[k]} \in C(\R^{k+1})$ denotes the $k^{\text{th}}$ divided difference of $f$ (Section \ref{sec.divdiff}).
    \item If $\Om_1,\ldots,\Om_{k+1}$ are compact subsets of $\R$, then $\ell^{\infty}(\Om_1,\cB_{\Om_1}) \iotimes \cdots \iotimes \ell^{\infty}(\Om_{k+1},\cB_{\Om_{k+1}})$ is the $\ell^{\infty}$-integral projective tensor product of $\ell^{\infty}(\Om_1,\cB_{\Om_1}),\ldots,\ell^{\infty}(\Om_{k+1},\cB_{\Om_{k+1}})$ (Section \ref{sec.IPTP}).
    \item If $H$ is a complex Hilbert space, $\cM \subseteq B(H)$ is a von Neumann algebra, $(a_1,\ldots,a_{k+1}) \in \cM_{\nu}^{k+1}$, $(b_1,\ldots,b_k) \in \cM^k$, and $\varphi \in \ell^{\infty}(\sigma(a_1),\cB_{\sigma(a_1)}) \iotimes \cdots \iotimes \ell^{\infty}(\sigma(a_{k+1}),\cB_{\sigma(a_{k+1})})$, then the multiple operator integral (MOI) in Eq. \eqref{eq.formalMOI} is to be interpreted in accordance with Section \ref{sec.MOIs}.
\end{itemize}
With these reminders out of the way, we introduce the space of functions we shall study.

\begin{defi}[Noncommutative $C^k$ Functions]
We say that $f \in C^k(\R)$ is \textbf{noncommutative} $\boldsymbol{C^k}$ and write $f \in NC^k(\R)$ if there is a sequence $(q_n)_{n \in \N}$ of polynomials such that for all $r > 0$ and $j \in \{0,\ldots,k\}$, we have $q_n^{[j]}\big|_{[-r,r]^{j+1}} \to f^{[j]}\big|_{[-r,r]^{j+1}}$ in $\ell^{\infty}([-r,r],\cB_{[-r,r]})^{\iotimes(j+1)}$ as $n \to \infty$.
\end{defi}

As we explain in Section \ref{sec.NCkdefandprop}, where we also develop some basic properties of $NC^k(\R)$, the name for the space $NC^k(\R)$ comes from parallel work by D. A. Jekel in \cite{jekel}.
Our first main result comes in the form of a list of examples (and a non-example) of noncommutative $C^k$ functions that paints the picture that a function only has to be ``slightly better than $C^k$" to be noncommutative $C^k$.
To state our result, we note that $W_k(\R)$ is the $k^{\text{th}}$ Wiener space (Definition \ref{def.Wk}), $\dot{B}_q^{s,p}(\R)$ is the homogeneous $(s,p,q)$-Besov space (Definition \ref{def.Besov}), and $C_{\loc}^{k,\e}(\R)$ is the space of $C^k$ functions such that $f^{(k)}$ is locally $\e$-H\"{o}lder continuous (Definition \ref{def.Holder}).
In addition, if $\cS \subseteq C^k(\R)$, then $\cS_{\loc}$ is defined to be the set of all $f \in C^k(\R)$ such that for all $r > 0$, there exists $g \in \cS$ such that $g|_{[-r,r]} = f|_{[-r,r]}$.

\begin{thm}[Examples of Noncommutative $C^k$ Functions]\label{thm.main1}
If $k \in \N$ and $\e > 0$, then
\begin{enumerate}[label=(\roman*),font=\normalfont,leftmargin=2\parindent]
    \item $\dot{B}_1^{k,\infty}(\R) \subseteq NC^k(\R)$, \label{item.Besov}
    \item $C_{\loc}^{k,\e}(\R) \subseteq NC^k(\R)$, and \label{item.Holder}
    \item $W_k(\R)_{\loc} \subsetneq NC^k(\R) \subsetneq C^k(\R)$. \label{item.stricts}
\end{enumerate}
\end{thm}
\begin{proof}
Item \ref{item.Besov} is part of Theorem \ref{thm.BesovNCk}, and item \ref{item.Holder} is Theorem \ref{thm.HolderNCk}.
The first containment in item \ref{item.stricts} (but not its strictness) follows from part of Theorem \ref{thm.WkNCk} and Proposition \ref{prop.cC[k]}\ref{item.Sloc};
an example demonstrating its strictness is given in Theorem \ref{thm.Wkloccex}.
An example demonstrating strictness of the second containment in item \ref{item.stricts} is given in Theorem \ref{thm.NCkcounterex}.
\end{proof}

Our second main result (proven in Section \ref{sec.diffopfunccalc}) is a statement about higher differentiability of operator functions associated to noncommutative $C^k$ functions, in particular functions in $\dot{B}_1^{k,\infty}(\R)$.

\begin{thm}[Derivatives of Operator Functions]\label{thm.main2}
Fix a unital $C^*$-algebra $\cA$ and a von Neumann algebra $\cM$ containing $\cA$ as a unital $C^*$-subalgebra.
If $k \in \N$ and $f \in NC^k(\R)$, then $f_{\mathsmaller{\cA}} \in C^k(\cA_{\sa};\cA)$ and
\[
\partial_{b_1}\cdots\partial_{b_k}f_{\mathsmaller{\cA}}(a) = \sum_{\pi \in S_k}\underbrace{\int_{\sigma(a)}\cdots\int_{\sigma(a)}}_{k+1 \, \mathrm{times}}f^{[k]}(\blambda) \, P^a(d\lambda_1) \, b_{\pi(1)} \cdots P^a(d\lambda_k) \, b_{\pi(k)} \, P^a(d\lambda_{k+1}), \numberthis\label{eq.main2}
\]
for all $a,b_1,\ldots,b_k \in \cA_{\sa}$.
\end{thm}

The purpose of $\cM$ is to allow us to make sense of the right hand side of Eq. \eqref{eq.main2}, since \textit{a priori} MOIs are only defined in von Neumann algebras.
(However, please see Lemma \ref{lem.CstarMOI}, which morally says that MOIs like the ones in Eq. \eqref{eq.main2} make sense in unital $C^*$-algebras.)
Of course, if $\cA$ happens to be a von Neumann algebra, then we may take $\cM = \cA$.
For an arbitrary (abstract) unital $C^*$-algebra $\cA$, a reasonable choice of $\cM$ is the double dual $\cA^{**}$ of $\cA$, which has a von Neumann algebra structure with respect to which the natural embedding $\cA \hookrightarrow \cA^{**}$ is a unital $\ast$-homomorphism (Theorem 1.17.2 in \cite{sakai}).
Therefore, if $a \in \cA_{\sa}$, then one may always interpret the projection-valued spectral measure $P^a$ in Eq. \eqref{eq.main2} as taking values in $\cA^{**}$, even when it does not make sense in $\cA$.
However, we highlight that the double dual $\cA^{**}$ of a $C^*$-algebra $\cA$ is frequently quite large;
specifically, it is frequently not representable on a separable Hilbert space.
This is why we must understand MOIs on non-separable Hilbert spaces (the subject of \cite{nikitopoulosMOI}).

\section{Preliminaries}\label{sec.prelim}

\subsection{Divided differences}\label{sec.divdiff}

In this section, we define divided differences and collect their relevant properties.

\begin{defi}[Divided Differences]\label{def.divdiff}
Let $f \colon \R \to \C$ be a function.
Define $f^{[0]} \coloneqq f$ and, for $k \in \N$ and distinct $\lambda_1,\ldots,\lambda_{k+1} \in \R$, recursively define
\[
f^{[k]}(\lambda_1,\ldots,\lambda_{k+1}) \coloneqq \frac{f^{[k-1]}(\lambda_1,\ldots,\lambda_k) - f^{[k-1]}(\lambda_1,\ldots,\lambda_{k-1},\lambda_{k+1})}{\lambda_k-\lambda_{k+1}}.
\]
We call $f^{[k]}$ the $\boldsymbol{k^{\textbf{th}}}$ \textbf{divided difference} of $f$.
\end{defi}

\begin{nota}
If $m \in \N$, then we define
\begin{align*}
     \Sigma_m & \coloneqq \Bigg\{\vec{s} = (s_1,\ldots,s_m) \in \R^m : |\vec{s}\,| = \sum_{j=1}^m |s_j| \leq 1 \text{ and } s_1,\ldots,s_m \geq 0\Bigg\} \; \text{ and} \\
     \Delta_m & \coloneqq \Bigg\{\boldsymbol{t} = (t_1,\ldots,t_{m+1}) \in \R^{m+1} : \sum_{j=1}^{m+1}t_j = 1 \text{ and } t_1,\ldots,t_{m+1} \geq 0\Bigg\}.
\end{align*}
Also, we denote by $\rho_m$ the pushforward of the $m$-dimensional Lebesgue measure on $\Sigma_m$ by the homeomorphism $\Sigma_m \ni \vec{s} \mapsto (\vec{s},1-|\vec{s}\,|) \in \Delta_m$.
Explicitly, $\rho_m$ is the finite Borel measure on $\Delta_m$ such that
\[
\int_{\Delta_m} \varphi(\boldsymbol{t}) \, \rho_m(d\boldsymbol{t}) = \int_{\Sigma_m} \varphi(\vec{s},1-|\vec{s}\,|) \, d\vec{s},
\]
for all $\varphi \in \ell^{\infty}(\Delta_m,\cB_{\Delta_m})$.
In particular, $\rho_m(\Delta_m) = \frac{1}{m!}$ as the reader may verify.
\end{nota}

The content of the following proposition is contained in Section 4.7 of \cite{devore}.
The first and third items are straightforward induction arguments.
We supply the proof of the second (and most important) item for the convenience of the reader.

\begin{prop}[Properties of $f^{[k]}$]\label{prop.divdiff}
Fix $k \in \N$, functions $f,g \colon \R \to \C$, and distinct $\lambda_1,\ldots,\lambda_{k+1} \in \R$.
Also, write $\blambda = (\lambda_1,\ldots,\lambda_{k+1})$.
\begin{enumerate}[label=(\roman*),font=\normalfont,leftmargin=2\parindent]
    \item We have
    \[
    f^{[k]}(\blambda) = \sum_{j=1}^{k+1}f(\lambda_j)\prod_{\ell \neq j}(\lambda_j-\lambda_{\ell})^{-1}.
    \]
    In particular, $f^{[k]}$ is symmetric.\label{item.symm}
    \item If $f \in C^k(\R)$, then
    \[
    f^{[k]}(\blambda) = \int_{\Delta_k}f^{(k)}(\boldsymbol{t} \cdot \blambda) \, \rho_k(d\boldsymbol{t}) = \int_{\Sigma_k}f^{(k)}\Bigg(\sum_{j=1}^ks_j\lambda_j+\Bigg(1-\sum_{j=1}^k s_j\Bigg)\lambda_{k+1}\Bigg) \,  ds_1 \cdots ds_k,
    \]
    where $\cdot$ is the Euclidean dot product.
    In particular, $f^{[k]}$ extends uniquely to a symmetric continuous function on all of $\R^{k+1}$.
    We shall use the same notation for this extension.\label{item.cont}
    \item The product rule
    \[
    (fg)^{[k]}(\blambda) = \sum_{j=0}^kf^{[j]}(\lambda_1,\ldots,\lambda_{j+1})\, g^{[k-j]}(\lambda_{j+1},\ldots,\lambda_{k+1})
    \]
    holds.
    If in addition $f,g \in C^k(\R)$, then this formula holds for all $\lambda_1,\ldots,\lambda_{k+1} \in \R$.\label{item.prod}
\end{enumerate}
\end{prop}
\begin{proof}[Proof of \ref{item.cont}]
If $f \in C^1(\R)$ and $\lambda_1,\lambda_2 \in \R$ are distinct, then
\[
f^{[1]}(\lambda_1,\lambda_2) = \frac{f(\lambda_1)-f(\lambda_2)}{\lambda_1-\lambda_2} = \int_0^1 f'(s_1\lambda_1+(1-s_1)\lambda_2) \, ds_1
\]
by the Fundamental Theorem of Calculus.
Now, assume the formula holds on $C^k(\R)$ with $k \geq 1$.
If $f \in C^{k+1}(\R)$ and $\lambda_1,\ldots,\lambda_{k+2} \in \R$ are distinct, then, writing $\vec{\lambda} \coloneqq (\lambda_1,\ldots,\lambda_k) \in \R^k$, we have
\begin{align*}
    f^{[k+1]}\big(\vec{\lambda},\lambda_{k+1},\lambda_{k+2}\big) & = \int_{\Sigma_k}\frac{f^{(k)}\big(\vec{\lambda}\cdot\vec{s} + (1-|\vec{s}\,|)\lambda_{k+1}\big) - f^{(k)}\big(\vec{\lambda}\cdot\vec{s} + (1-|\vec{s}\,|)\lambda_{k+2}\big)}{\lambda_{k+1}-\lambda_{k+2}} \, d\vec{s} \numberthis \label{eq.indhyp} \\
    & = \int_{\Sigma_k} (1-|\vec{s}\,|)\,\big(f^{(k)}\big)^{[1]}\big(\vec{\lambda} \cdot \vec{s} + (1-|\vec{s}\,|)\lambda_{k+1},\vec{\lambda} \cdot \vec{s} + (1-|\vec{s}\,|)\lambda_{k+2}\big)\,d\vec{s} \\
    & = \int_{\Sigma_k}(1-|\vec{s}\,|)\int_0^1f^{(k+1)}\big(\vec{\lambda}\cdot \vec{s} + (1-|\vec{s}\,|)(t\lambda_{k+1}+(1-t)\lambda_{k+2})\big) \, dt \, d\vec{s} \numberthis \label{eq.k=1} \\
    & = \int_{\Sigma_k}\int_0^{1-|\vec{s}\,|} f^{(k+1)}\big(\vec{\lambda}\cdot \vec{s} + s_{k+1}\lambda_{k+1} + (1-|\vec{s}\,|-s_{k+1})\lambda_{k+2}\big) \, ds_{k+1} \, d\vec{s} \numberthis \label{eq.deltavar} \\
    & = \int_{\Sigma_{k+1}}f^{(k+1)}\Bigg(\sum_{j=1}^{k+1}s_j\lambda_j+\Bigg(1-\sum_{j=1}^{k+1} s_j\Bigg)\lambda_{k+2}\Bigg) \,  ds_1 \cdots ds_{k+1},
\end{align*}
where Eq. \eqref{eq.indhyp} holds by definition of $f^{[k+1]}$ and the inductive hypothesis, Eq. \eqref{eq.k=1} holds by the $k=1$ case, and Eq. \eqref{eq.deltavar} holds by the change of variable $s_{k+1} \coloneqq (1-|\vec{s}\,|)t$.
\end{proof}

\begin{cor}\label{cor.easyboundfk}
Fix $k \in \N$.
If $f \in C^k(\R)$ and $\boldsymbol{\lambda} = (\lambda_1,\ldots,\lambda_{k+1}) \in \R^{k+1}$, then
\[
\big|f^{[k]}(\blambda)\big| \leq \frac{1}{k!}\sup_{|s| \leq |\blambda|_{\infty}}\big|f^{(k)}(s)\big| \; \text{ and } \; \big\|f^{[k]}\big\|_{\ell^{\infty}([-r,r]^{k+1})} = \frac{1}{k!}\big\|f^{(k)}\big\|_{\ell^{\infty}([-r,r])},
\]
for all $r > 0$, where $|\blambda|_{\infty} = \max\{|\lambda_j| : 1 \leq j \leq k+1\}$.
\end{cor}
\begin{proof}
The inequality follows from Proposition \ref{prop.divdiff}\ref{item.cont}, the observation that $\boldsymbol{t} \in \Delta_k \Rightarrow|\boldsymbol{t} \cdot \blambda| \leq |\blambda|_{\infty}$, and the fact that $\rho_k(\Delta_k) = \frac{1}{k!}$.
The identity follows from the inequality and the fact that if $\lambda \in \R$, then $f^{[k]}(\lambda,\ldots,\lambda) = f^{(k)}(\lambda) \,\rho_k(\Delta_k) = \frac{1}{k!}f^{(k)}(\lambda)$ also by the formula in Proposition \ref{prop.divdiff}\ref{item.cont}.
\end{proof}

We end this section with an important example calculation of a divided difference.

\begin{ex}[Polynomials]\label{ex.polydivdiff}
Let $n \in \N_0$ and $p_n(\lambda) \coloneqq \lambda^n \in \C[\lambda] \subseteq C^{\infty}(\R)$.
If $k \in \N$, then we claim that
\[
p_n^{[k]}(\blambda) = \sum_{|\gamma|=n-k} \blambda^{\gamma} = \sum_{\gamma \in \N_0^{k+1} : |\gamma| = n-k} \lambda_1^{\gamma_1}\cdots \lambda_{k+1}^{\gamma_{k+1}},\numberthis\label{eq.pmk}
\]
for all $\blambda = (\lambda_1,\ldots,\lambda_{k+1}) \in \R^{k+1}$, where $|\gamma| = \sum_{j=1}^{k+1}\gamma_j$ for $\gamma = (\gamma_1,\ldots,\gamma_{k+1}) \in \N_0^{k+1}$.
(Above and throughout, empty sums are defined to be zero.)
By continuity (Proposition \ref{prop.divdiff}\ref{item.cont}), it suffices to prove that Eq. \eqref{eq.pmk} holds for $\blambda = (\lambda_1,\ldots,\lambda_{k+1})$ such that $\lambda_1,\ldots,\lambda_{k+1}$ are distinct.
To do so, we proceed by induction on $k$.
Of course, the identity $p_n^{[1]}(\lambda_1,\lambda_2) = \frac{\lambda_1^n-\lambda_2^n}{\lambda_1-\lambda_2} = \sum_{j=0}^{n-1}\lambda_1^j\lambda_2^{n-1-j} = \sum_{\gamma_1+\gamma_2=n-1}\lambda_1^{\gamma_1}\lambda_2^{\gamma_2}$ is well-known.
Now, suppose Eq. \eqref{eq.pmk} holds.
If $\lambda_{k+2}$ is distinct from $\lambda_1,\ldots,\lambda_{k+1}$, then
\begin{align*}
    p_n^{[k+1]}(\lambda_1,\ldots,\lambda_{k+2}) & = \frac{p_n^{[k]}(\lambda_1,\ldots,\lambda_{k+1}) - p_n^{[k]}(\lambda_1,\ldots,\lambda_k,\lambda_{k+2})}{\lambda_{k+1}-\lambda_{k+2}} = \sum_{\gamma \in \N_0^{k+1} : |\gamma|=n-k} \lambda_1^{\gamma_1}\cdots\lambda_k^{\gamma_k}\frac{\lambda_{k+1}^{\gamma_{k+1}}-\lambda_{k+2}^{\gamma_{k+1}}}{\lambda_{k+1}-\lambda_{k+2}} \\
    & = \sum_{\gamma \in \N_0^{k+1} : |\gamma|=n-k}\Bigg(\sum_{\delta_1+\delta_2=\gamma_{k+1}-1}\lambda_1^{\gamma_1}\cdots\lambda_k^{\gamma_k}\lambda_{k+1}^{\delta_1}\lambda_{k+2}^{\delta_2}\Bigg) = \sum_{\tilde{\gamma} \in \N_0^{k+2} : |\tilde{\gamma}| = n-k-1} \lambda_1^{\tilde{\gamma}_1}\cdots\lambda_{k+2}^{\tilde{\gamma}_{k+2}},
\end{align*}
by the inductive hypothesis and the $k=1$ case.
This completes the proof of the claim.
In particular, if $p \in \C[\lambda]$, then $p^{[k]} \in \C[\lambda_1,\ldots,\lambda_{k+1}]$, for all $k \in \N$.
\end{ex}

\subsection{Integral projective tensor products}\label{sec.IPTP}

For the duration of this section, fix $m \in \N$ and Polish spaces (i.e., complete separable metric spaces) $\Om_1,\ldots,\Om_m$.
Also, write $\Om \coloneqq \Om_1 \times \cdots \times \Om_m$.
We first review the notion of the integral projective tensor product $\ell^{\infty}(\Om_1,\cB_{\Om_1}) \iotimes \cdots \iotimes \ell^{\infty}(\Om_m,\cB_{\Om_m})$, the idea for which is due to Peller \cite{peller1}.

\begin{lem}[Measurability]\label{lem.meas}
Let $\Xi$ be a Polish space and $(\Sigma,\sH,\rho)$ be a $\sigma$-finite measure space.
If $\varphi \colon \Xi \times \Sigma \to \C$ is product measurable, then the function
\[
\Sigma \ni \sigma \mapsto \|\varphi(\cdot,\sigma)\|_{\ell^{\infty}(\Xi)} \in [0,\infty]
\]
is $\big(\overline{\sH}^{\rho},\cB_{[0,\infty]}\big)$-measurable, where $\overline{\sH}^{\rho}$ is the $\rho$-completion of $\sH$.
\end{lem}
\begin{proof}
Since every $\sigma$-finite measure is equivalent to (i.e., has the same null sets as) a finite measure, we may assume $\rho$ is finite.
Let $C \colon \Sigma \to 2^{\Xi}$ be such that $\{(\sigma,\xi) : \sigma \in \Sigma, \,\xi \in C(\sigma)\} \in \sH \otimes \cB_{\Xi}$.
Then Corollary 2.13 in \cite{crauel} (which relies on the Measurable Projection Theorem, Theorem III.23 in \cite{castaing}) implies that
\[
\Sigma \ni \sigma \mapsto \sup_{\xi \in C(\sigma)}|\varphi(\xi,\sigma)| \in [0,\infty]
\]
is $\big(\overline{\sH}^{\rho},\cB_{[0,\infty]}\big)$-measurable.
Applying this to the constant function $C \equiv \Xi$ gives the desired result.
\end{proof}

This measurability lemma ensures that the integral in Eq. \eqref{eq.intfincond} below makes sense as a Lebesgue integral with respect to the completion of $\rho$.

\begin{defi}[IPTPs]\label{def.babyIPTP}
A $\boldsymbol{\ell^{\infty}}$\textbf{-integral projective decomposition} (IPD) of a function $\varphi \colon \Om \to \C$ is a choice $(\Sigma,\rho,\varphi_1,\ldots,\varphi_m)$ of a $\sigma$-finite measure space $(\Sigma,\sH,\rho)$ and, for each $j \in \{1,\ldots,m\}$, a product measurable function $\varphi_j \colon \Om_j \times \Sigma \to \C$ such that $\varphi_j(\cdot,\sigma) \in \ell^{\infty}(\Om_j,\cB_{\Om_j})$ whenever $\sigma\in \Sigma$,
\[
\int_{\Sigma} \|\varphi_1(\cdot,\sigma)\|_{\ell^{\infty}(\Om_1)}\cdots\|\varphi_m(\cdot,\sigma)\|_{\ell^{\infty}(\Om_m)} \, \rho(d\sigma) < \infty, \numberthis\label{eq.intfincond}
\]
and
\[
\varphi(\boldsymbol{\om}) = \int_{\Sigma} \varphi_1(\om_1,\sigma) \cdots \varphi_m(\om_m,\sigma) \, \rho(d\sigma), \; \text{ for all } \, \boldsymbol{\om} = (\om_1,\ldots,\om_m) \in \Om.
\]
Also, for any function $\varphi \colon \Om \to \C$, define
\[
\|\varphi\|_{\ell^{\infty}(\Om_1,\cB_{\Om_1}\hspace{-0.1mm}) \iotimes \cdots \iotimes \ell^{\infty}(\Om_m,\cB_{\Om_m}\hspace{-0.1mm})} \coloneqq \inf\Bigg\{\int_{\Sigma} \prod_{j=1}^m\|\varphi_j(\cdot,\sigma)\|_{\ell^{\infty}(\Om_j)}\,\rho(d\sigma) : (\Sigma,\rho,\varphi_1,\ldots,\varphi_m)  \text{ is a } \ell^{\infty}\text{-IPD of } \varphi\Bigg\},
\]
where $\inf \emptyset \coloneqq \infty$.
Finally, we define
\[
\ell^{\infty}(\Om_1,\cB_{\Om_1}) \iotimes \cdots \iotimes \ell^{\infty}(\Om_m,\cB_{\Om_m}) \coloneqq \big\{\varphi \in \ell^{\infty}(\Om,\cB_{\Om}) : \|\varphi\|_{\ell^{\infty}(\Om_1,\cB_{\Om_1}\hspace{-0.1mm}) \iotimes \cdots \iotimes \ell^{\infty}(\Om_m,\cB_{\Om_m}\hspace{-0.1mm})} < \infty\big\}
\]
to be the \textbf{integral projective tensor product of} $\boldsymbol{\ell^{\infty}(\Om_1,\cB_{\Om_1}),\ldots,\ell^{\infty}(\Om_m,\cB_{\Om_m})}$.
\end{defi}

In the next statement, a \textbf{$\boldsymbol{\ast}$-algebra} is a unital $\C$-algebra with a \textbf{$\boldsymbol{\ast}$-operation} (a unital, conjugate-linear, anti-multiplicative involution), and a \textbf{Banach $\boldsymbol{\ast}$-algebra} is a complex unital Banach algebra with an isometric $\ast$-operation.

\begin{prop}[Properties of IPTPs]\label{prop.babyIPTP}
The following hold.
\begin{enumerate}[label=(\roman*),font=\normalfont,leftmargin=2\parindent]
    \item If $\varphi \colon \Om \to \C$ is a function, then $\|\varphi\|_{\ell^{\infty}(\Om)} \leq \|\varphi\|_{\ell^{\infty}(\Om_1,\cB_{\Om_1}\hspace{-0.1mm}) \iotimes \cdots \iotimes \ell^{\infty}(\Om_m,\cB_{\Om_m}\hspace{-0.1mm})}$.\label{item.ellinfbd}
    \item $\ell^{\infty}(\Om_1,\cB_{\Om_1}) \iotimes \cdots \iotimes \ell^{\infty}(\Om_m,\cB_{\Om_m}) \subseteq \ell^{\infty}(\Om,\cB_{\Om})$ is a unital $\ast$-subalgebra, and
    \[
    \big(\ell^{\infty}(\Om_1,\cB_{\Om_1}) \iotimes \cdots \iotimes \ell^{\infty}(\Om_m,\cB_{\Om_m}),\|\cdot\|_{\ell^{\infty}(\Om_1,\cB_{\Om_1}\hspace{-0.1mm}) \iotimes \cdots \iotimes \ell^{\infty}(\Om_m,\cB_{\Om_m}\hspace{-0.1mm})}\big)
    \]
    is a Banach $\ast$-algebra under pointwise operations.\label{item.Banachstaralg}
    \item Suppose $1 \leq j,k \leq m$.
    If $\varphi \colon \Om_1 \times \cdots \times \Om_j \to \C$ and $\psi \colon \Om_k \times \cdots \times \Om_m \to \C$ are functions and we define $\chi(\om_1,\ldots,\om_m) \coloneqq \varphi(\om_1,\ldots,\om_j)\, \psi(\om_k,\ldots,\om_m)$, then
    \[
    \|\chi\|_{\ell^{\infty}(\Om_1,\cB_{\Om_1}\hspace{-0.1mm}) \iotimes \cdots \iotimes \ell^{\infty}(\Om_m,\cB_{\Om_m}\hspace{-0.1mm})} \leq \|\varphi\|_{\ell^{\infty}(\Om_1,\cB_{\Om_1}\hspace{-0.1mm}) \iotimes \cdots \iotimes \ell^{\infty}(\Om_j,\cB_{\Om_j}\hspace{-0.1mm})}\|\psi\|_{\ell^{\infty}(\Om_k,\cB_{\Om_k}\hspace{-0.1mm}) \iotimes \cdots \iotimes \ell^{\infty}(\Om_m,\cB_{\Om_m}\hspace{-0.1mm})},
    \]
    where as usual $0 \cdot \infty \coloneqq 0$.\label{item.tensorIPTP}
\end{enumerate}
\end{prop}
\begin{proof}
For ease of notation, write $\sB \coloneqq \ell^{\infty}(\Om_1,\cB_{\Om_1}) \iotimes \cdots \iotimes \ell^{\infty}(\Om_m,\cB_{\Om_m})$ and
\[
\|\cdot\|_{\sB} \coloneqq \|\cdot\|_{\ell^{\infty}(\Om_1,\cB_{\Om_1}\hspace{-0.1mm}) \iotimes \cdots \iotimes \ell^{\infty}(\Om_m,\cB_{\Om_m}\hspace{-0.1mm})}.
\]
We take each item in turn.

\ref{item.ellinfbd} Let $(\Sigma,\rho,\varphi_1,\ldots,\varphi_m)$ be a $\ell^{\infty}$-IPD of $\varphi$.
If $\boldsymbol{\om} = (\om_1,\ldots,\om_m) \in \Om$, then
\[
|\varphi(\boldsymbol{\om})| \leq \int_{\Sigma} |\varphi_1(\om_1,\sigma) \cdots \varphi_m(\om_m,\sigma)|\,\rho(d\sigma) \leq \int_{\Sigma} \|\varphi_1(\cdot,\sigma)\|_{\ell^{\infty}(\Om_1)}\cdots\|\varphi_m(\cdot,\sigma)\|_{\ell^{\infty}(\Om_m)} \,\rho(d\sigma).
\]
Taking the supremum over $\boldsymbol{\om} \in \Om$ and then the infimum over $\ell^{\infty}$-IPDs $(\Sigma,\rho,\varphi_1,\ldots,\varphi_m)$ gives the desired inequality.
Note that this inequality implies that $\|\varphi\|_{\sB} = 0 \Leftrightarrow \varphi \equiv 0$ on $\Om$.

\ref{item.Banachstaralg} We leave it to the reader to prove $\|c \, \varphi\|_{\sB} = |c|\,\|\varphi\|_{\sB} = |c|\,\|\overline{\varphi}\|_{\sB}$ whenever $c \in \C$ and $\varphi \in \sB$.
Next, suppose $(\varphi_n)_{n \in \N}$ is a sequence in $\sB$ such that $\sum_{n=1}^{\infty} \|\varphi_n\|_{\sB} < \infty$.
By the previous item, we have $\sum_{n=1}^{\infty}\|\varphi_n\|_{\ell^{\infty}(\Om)} \leq \sum_{n=1}^{\infty} \|\varphi_n\|_{\sB} < \infty$, so the series $\varphi \coloneqq \sum_{n=1}^{\infty}\varphi_n$ converges in $\ell^{\infty}(\Om,\cB_{\Om})$.
We claim that
\[
\|\varphi\|_{\sB} \leq \sum_{n=1}^{\infty} \|\varphi_n\|_{\sB}.
\]
To see this, fix $\e > 0$ and $n \in \N$.
Then, by definition of $\|\cdot\|_{\sB}$, there exists a $\ell^{\infty}$-IPD $(\Sigma_n,\rho_n,\varphi_{n,1},\ldots,\varphi_{n,m})$ of $\varphi_n$ such that
\[
\int_{\Sigma_n} \|\varphi_{n,1}(\cdot,\sigma_n)\|_{\ell^{\infty}(\Om_1)}\cdots\|\varphi_{n,m}(\cdot,\sigma_n)\|_{\ell^{\infty}(\Om_m)}\,\rho_n(d\sigma_n) < \|\varphi_n\|_{\sB} + \frac{\e}{2^n}.
\]
Define $(\Sigma,\sH,\rho)$ to be the disjoint union of the measure spaces $\{(\Sigma_n,\sH_n,\rho_n) : n \in \N\}$.
Also, for $1 \leq j \leq m$, define $\chi_j \colon \Om_j \times \Sigma \to \C$ to be the unique measurable function satisfying $\chi_j|_{\Om_j \times \Sigma_n} = \varphi_{n,j}$, for all $n \in \N$.
Then it is easy to see that $(\Sigma,\rho,\chi_1,\ldots,\chi_m)$ is a $\ell^{\infty}$-IPD of $\varphi$, so that
\begin{align*}
    \|\varphi\|_{\sB} & \leq \int_{\Sigma}\|\chi_1(\cdot,\sigma)\|_{\ell^{\infty}(\Om_1)}\cdots\|\chi_m(\cdot,\sigma)\|_{\ell^{\infty}(\Om_m)} \,\rho(d\sigma) \\
    & = \sum_{n = 1}^{\infty}\int_{\Sigma_n} \|\varphi_{n,1}(\cdot,\sigma_n)\|_{\ell^{\infty}(\Om_1)}\cdots\|\varphi_{n,m}(\cdot,\sigma_n)\|_{\ell^{\infty}(\Om_m)}\,\rho_n(d\sigma_n) \\
    & \leq \sum_{n = 1}^{\infty}\Big(\|\varphi_n\|_{\sB} + \frac{\e}{2^n}\Big) = \sum_{n = 1}^{\infty}\|\varphi_n\|_{\sB} + \e.
\end{align*}
Taking $\e \searrow 0$ results in the desired estimate.
Taking $\varphi_n \equiv 0$ for $n \geq 3$, we conclude that $\sB$ is closed under addition and that $\|\cdot\|_{\sB}$ satisfies the triangle inequality.
Applying the inequality we just proved to the sequence $(\varphi_{n+N})_{n \in \N}$ for fixed $N \in \N$ yields
\[
\Bigg\|\varphi-\sum_{n=1}^N\varphi_n\Bigg\|_{\sB} \leq \sum_{n=N+1}^{\infty}\|\varphi_n\|_{\sB} \to 0
\]
as $N \to \infty$.
Combining this with the observation from the end of the proof of the previous item, we conclude that $\sB$ is a Banach space.

Finally, we prove that if $\varphi,\psi \in \sB$, then $\|\varphi\,\psi\|_{\sB} \leq \|\varphi\|_{\sB}\|\psi\|_{\sB}$.
To do so, fix $\ell^{\infty}$-IPDs $(\Sigma_1,\rho_1,\varphi_1,\ldots,\varphi_m)$ and $(\Sigma_2,\rho_2,\psi_1,\ldots,\psi_m)$ of $\varphi$ and $\psi$, respectively.
Next, redefine $(\Sigma,\sH,\rho) \coloneqq (\Sigma_1 \times \Sigma_2,\sH_1 \otimes \sH_2,\rho_1 \otimes \rho_2)$ and $\chi_j(\om_j,\sigma) \coloneqq \varphi_j(\om_j,\sigma_1)\,\psi_j(\om_j,\sigma_2)$ whenever $1 \leq j \leq m$, $\om_j \in \Om_j$, and $\sigma = (\sigma_1,\sigma_2) \in \Sigma$.
We claim that $(\Sigma,\rho,\chi_1,\ldots,\chi_m)$ is a $\ell^{\infty}$-IPD of $\varphi \,\psi$.
Indeed, by Tonelli's Theorem,
\[
\int_{\Sigma} \prod_{j=1}^m\|\chi_j(\cdot,\sigma)\|_{\ell^{\infty}(\Om_j)}\,\rho(d\sigma) \leq \int_{\Sigma_1} \prod_{j=1}^m\|\varphi_j(\cdot,\sigma_1)\|_{\ell^{\infty}(\Om_j)}\,\rho_1(d\sigma_1)\,\int_{\Sigma_2} \prod_{j=1}^m\|\psi_j(\cdot,\sigma_2)\|_{\ell^{\infty}(\Om_j)}\,\rho_2(d\sigma_2) < \infty.\pagebreak
\]
Now, by Fubini's Theorem,
\[
\int_{\Sigma} \prod_{j=1}^m\chi_j(\om_j,\sigma)\,\rho(d\sigma) = \int_{\Sigma_1} \prod_{j=1}^m\varphi_j(\om_j,\sigma_1)\,\rho_1(d\sigma_1)\,\int_{\Sigma_2} \prod_{j=1}^m\psi_j(\om_j,\sigma_2)\,\rho_2(d\sigma_2) = \varphi(\boldsymbol{\om})\,\psi(\boldsymbol{\om}),
\]
for all $\boldsymbol{\om} = (\om_1,\ldots,\om_m) \in \Om$.
It follows that
\[
\|\varphi \, \psi\|_{\sB} \leq \int_{\Sigma_1} \prod_{j=1}^m\|\varphi_j(\cdot,\sigma_1)\|_{\ell^{\infty}(\Om_j)}\,\rho_1(d\sigma_1)\,\int_{\Sigma_2} \prod_{j=1}^m\|\psi_j(\cdot,\sigma_2)\|_{\ell^{\infty}(\Om_j)}\,\rho_2(d\sigma_2).
\]
Taking the infimum over $\ell^{\infty}$-IPDs of $\varphi$ and $\psi$ gives the desired result.

\ref{item.tensorIPTP} By the previous item, it suffices to consider the cases $\varphi \equiv 1$ and $\psi \equiv 1$.
We leave these cases to the reader, as they are easy consequences of the definitions.
\end{proof}

We shall work mainly with the case $\Om_1 = \cdots = \Om_m = [-r,r]$, for which we use the following notation.

\begin{nota}
If $\varphi \colon \R^m \to \C$ is a function, then we write
\[
\|\varphi\|_{r,m} \coloneqq \big\|\varphi|_{[-r,r]^m}\big\|_{\ell^{\infty}([-r,r],\cB_{[-r,r]})^{\iotimes m}} \in [0,\infty]
\]
for all $r > 0$.
\end{nota}

\begin{ex}[Multivariate Polynomials]\label{ex.mvarpoly}
Fix $m \in \N$, and suppose
\[
P(\blambda) = \sum_{|\alpha| \leq d}c_{\alpha} \blambda^{\alpha} = \sum_{\alpha \in \N_0^m : |\alpha| \leq d} c_{\alpha}\, \lambda_1^{\alpha_1}\cdots \lambda_m^{\alpha_m} \in \C[\lambda_1,\ldots,\lambda_m].
\]
For $r > 0$ and $\lambda_1,\ldots,\lambda_m \in [-r,r]$, let $\Sigma \coloneqq \{\alpha \in \N_0^m : |\alpha| \leq d\}$, $\rho$ be the counting measure on $\Sigma$, $\varphi_1(\lambda_1,\alpha) \coloneqq c_{\alpha}\,\lambda_1^{\alpha_1}$, and $\varphi_j(\lambda_j,\alpha) \coloneqq \lambda_j^{\alpha_j}$ for $j \in \{2,\ldots,m\}$ and $\alpha \in \Sigma$.
Then $(\Sigma,\rho,\varphi_1,\ldots,\varphi_m)$ is a $\ell^{\infty}$-IPD of $P|_{[-r,r]^m}$ that gives
\[
\|P\|_{r,m} \leq \int_{\Sigma} \prod_{j=1}^m\|\varphi_j(\cdot,\alpha)\|_{\ell^{\infty}([-r,r])}\,\rho(d\alpha) = \sum_{|\alpha| \leq d} |c_{\alpha}| \sup_{|\lambda_1| \leq r}\big|\lambda_1^{\alpha_1}\big|\cdots \sup_{|\lambda_m| \leq r}\big|\lambda_m^{\alpha_m}\big| = \sum_{|\alpha| \leq d} |c_{\alpha}| \, r^{|\alpha|} < \infty.
\]
Thus $P|_{[-r,r]^m} \in \ell^{\infty}([-r,r],\cB_{[-r,r]})^{\iotimes m}$.
In particular (by Example \ref{ex.polydivdiff}), if $p \in \C[\lambda]$, then
\[
p^{[k]}\big|_{[-r,r]^{k+1}} \in \ell^{\infty}\big([-r,r],\cB_{[-r,r]}\big)^{\iotimes (k+1)},
\]
for all $k \in \N$.
\end{ex}

\section{Noncommutative \texorpdfstring{$C^k$}{} functions}\label{sec.NCk}

\subsection{Definition and basic properties}\label{sec.NCkdefandprop}

In this section, we define our space of noncommutative $C^k$ functions and develop some of its basic properties.
First, however, we recall that if $k \in \N_0 \cup \{\infty\}$, then the space $C^k(\R)$ is a Fr\'{e}chet space with respect to the topology of locally uniform convergence of all derivatives of order at most $k$.
(Please see Proposition \ref{prop.Frechcomp}.)
This topology is called the \textbf{$\boldsymbol{C^k}$ topology}, and it is induced by the family
\[
\big\{f \mapsto \big\|f^{(j)}\big\|_{\ell^{\infty}([-r,r])} : 0 \leq j < k+1, \, r > 0\big\}
\]
of seminorms.
By Corollary \ref{cor.easyboundfk}, the $C^k$ topology is also induced by the family
\[
\big\{f \mapsto \big\|f^{[j]}\big\|_{\ell^{\infty}([-r,r]^{j+1})} : 0 \leq j < k+1, \, r > 0\big\}
\]
of seminorms.
For our space(s) of functions, we shall measure $f^{[j]}$ with a different family of seminorms.

\begin{nota}[The Space $\cC^{[k]}(\R)$]
For $k \in \N$, $f \in C^k(\R)$, and $r > 0$, define
\begin{align*}
    \|f\|_{\cC^{[k]},r} & \coloneqq \sum_{j=0}^k\big\|f^{[j]}\big\|_{r,j+1} \in [0,\infty] \; \text{ and} \\
    \cC^{[k]}(\R) & \coloneqq \big\{ f \in C^k(\R) : \|f\|_{\cC^{[k]},r}  < \infty, \text{ for all } r > 0\big\} \\
    & = \big\{f \in C^k(\R) : f^{[j]}\big|_{[-r,r]^{j+1}} \in \ell^{\infty}\big([-r,r],\cB_{[-r,r]}\big)^{\iotimes(j+1)} \text{ whenever } 1 \leq j \leq k \text{ and } r > 0\big\}.
\end{align*}
Also, let $\cC^{[\infty]}(\R) \coloneqq \bigcap_{k \in\N} \cC^{[k]}(\R)$.
\end{nota}

\begin{ex}[Polynomials]\label{ex.polycC[k]}
By Example \ref{ex.mvarpoly}, $\C[\lambda] \subseteq \cC^{[\infty]}(\R)$.
\end{ex}

Fix $k \in \N$.
By Proposition \ref{prop.babyIPTP}, $\cC^{[k]}(\R) \subseteq C^k(\R)$ is a linear subspace and $\{\|\cdot\|_{\cC^{[k]},r} : r > 0\}$ is a collection of seminorms on $\cC^{[k]}(\R)$.
Since these seminorms clearly separate points, they make $\cC^{[k]}(\R)$ into a Hausdorff locally convex topological vector space (LCTVS).
Similarly, $\cC^{[\infty]}(\R)$ is a Hausdorff LCTVS with topology induced by the family $\{\|\cdot\|_{\cC^{[k]},r} : k \in \N, \, r > 0\}$ of seminorms.
Here now are the basic properties of the spaces $\cC^{[k]}(\R)$, $k \in \N \cup \{\infty\}$.
In the result below, $\hookrightarrow$ indicates continuity of the relevant inclusion map, and a \textbf{Fr\'{e}chet $\boldsymbol{\ast}$-algebra} is a complex Fr\'{e}chet space with a $\ast$-algebra structure such that the $\ast$-operation and product are continuous.

\begin{prop}[Properties of $\cC^{[k]}(\R)$]\label{prop.cC[k]}
Fix $k \in \N \cup \{\infty\}$.
\begin{enumerate}[label=(\roman*),font=\normalfont,leftmargin=2\parindent]
    \item $\cC^{[k]}(\R) \hookrightarrow C^k(\R)$.\label{item.embed}
    \item For $\cS \subseteq C^k(\R)$, write
    \[
    \cS_{\loc} \coloneqq \big\{f \in C^k(\R) : \text{for all } r > 0, \text{ there exists } g \in \cS \text{ such that } g|_{[-r,r]} = f|_{[-r,r]}\big\}.
    \]
    Then $\cS \subseteq \cC^{[k]}(\R) \implies \cS_{\loc} \subseteq \overline{\cS} \subseteq \cC^{[k]}(\R)$.\label{item.Sloc}
    \item If $k < \infty$ and $f,g \in C^k(\R)$, then
    \[
    \big\|(fg)^{[k]}\big\|_{r,k+1} \leq \sum_{j=0}^k\big\|f^{[j]}\big\|_{r,j+1}\|g^{[k-j]}\big\|_{r,k-j+1} \; \text{ and } \; \|fg\|_{\cC^{[k]},r} \leq \|f\|_{\cC^{[k]},r}\|g\|_{\cC^{[k]},r},
    \]
    for all $r > 0$.\label{item.cC[k]prodrule}
    \item $\cC^{[k]}(\R)$ is a Fr\'{e}chet $\ast$-algebra under pointwise operations.\label{item.cC[k]Frechstaralg}
\end{enumerate}
\end{prop}
\begin{proof}
We take each item in turn.

\ref{item.embed} Fix $f \in C^k(\R)$.
If $0 \leq j < k+1$ and $r > 0$, then
\[
\big\|f^{[j]}\big\|_{\ell^{\infty}([-r,r]^{j+1})} \leq \big\|f^{[j]}\big\|_{r,j+1}
\]
by Proposition \ref{prop.babyIPTP}\ref{item.ellinfbd}.
Therefore, by the description of the $C^k$ topology given at the beginning of this section, the inclusion of $\cC^{[k]}(\R)$ into $C^k(\R)$ is continuous.

\ref{item.Sloc} Fix $\cS \subseteq \cC^{[k]}(\R)$.
If $f \in \cS_{\loc}$ and $n \in \N$, then --- by definition --- we can choose $g_n \in \cS \subseteq \cC^{[k]}(\R)$ such that $g_n|_{[-n,n]} = f|_{[-n,n]}$.
If $r > 0$, $n > r$, and $0 \leq j < k+1$, then $\|g_n-f\|_{\cC^{[j]},r} = 0$.
Thus $f \in \cC^{[k]}(\R)$ and $g_n \to f$ in $\cC^{[k]}(\R)$ as $n \to \infty$.
In other words, $\cS_{\loc} \subseteq \overline{\cS} \subseteq \cC^{[k]}(\R)$.

\ref{item.cC[k]prodrule} The claimed bound on $\big\|(fg)^{[k]}\big\|_{r,k+1}$ follows immediately from Propositions \ref{prop.divdiff}\ref{item.prod} and \ref{prop.babyIPTP}\ref{item.tensorIPTP}.
Consequently, we get $\|fg\|_{\cC^{[k]},r} \leq \|f\|_{\cC^{[k]},r}\|g\|_{\cC^{[k]},r}$.
Indeed,
\begin{align*}
    \|fg\|_{\cC^{[k]},r} & = \sum_{j=0}^k\big\|(fg)^{[j]}\big\|_{r,j+1} \leq \sum_{j=0}^k\sum_{\ell=0}^j\big\|f^{[\ell]}\big\|_{r,\ell+1} \,\big\|g^{[j-\ell]}\big\|_{r,j-\ell+1} \\
    & = \sum_{\ell=0}^k\big\|f^{[\ell]}\big\|_{r,\ell+1}\sum_{j=\ell}^k \big\|g^{[j-\ell]}\big\|_{r,j-\ell+1} \leq \|f\|_{\cC^{[k]},r}\|g\|_{\cC^{[k]},r},
\end{align*}
as desired.
\pagebreak

\ref{item.cC[k]Frechstaralg} For ease of notation, we assume $k < \infty$ and leave the case $k=\infty$ as an exercise.
First, note that the topology of $\cC^{[k]}(\R)$ is generated by the countable family $\{\|\cdot\|_{\cC^{[k]},N} : N \in \N\}$ of seminorms.
Thus $\cC^{[k]}(\R)$ is metrizable.
Next, we prove that $\cC^{[k]}(\R)$ is complete.
To this end, let $(f_n)_{n \in \N}$ be a Cauchy sequence in $\cC^{[k]}(\R)$.
By the first item, the sequence $(f_n)_{n \in \N}$ is also Cauchy in $C^k(\R)$.
By completeness of the latter space, there exists $f \in C^k(\R)$ such that $f_n \to f$ in the $C^k$ topology as $n \to \infty$.
In particular, if $1 \leq j \leq k$, then $f_n^{[j]} \to f^{[j]}$ uniformly on compact sets as $n \to \infty$.
Now, fix $j \in \{1,\ldots,k\}$ and $r > 0$.
Then
\[
\big(f_n^{[j]}\big|_{[-r,r]^{j+1}}\big)_{n \in \N}
\]
is Cauchy and therefore, by Proposition \ref{prop.babyIPTP}\ref{item.Banachstaralg}, convergent in $\ell^{\infty}([-r,r],\cB_{[-r,r]})^{\iotimes(j+1)}$.
Since we already know that $f_n^{[j]} \to f^{[j]}$ pointwise as $n \to \infty$, we conclude that $f^{[j]}|_{[-r,r]^{j+1}} \in \ell^{\infty}([-r,r],\cB_{[-r,r]})^{\iotimes(j+1)}$ and
\[
f_n^{[j]}|_{[-r,r]^{j+1}} \to f^{[j]}|_{[-r,r]^{j+1}} \text{ in } \ell^{\infty}([-r,r],\cB_{[-r,r]})^{\iotimes(j+1)}
\]
as $n \to \infty$ as well.
Thus $f \in \cC^{[k]}(\R)$ and $f_n \to f$ in $\cC^{[k]}(\R)$ as $n \to \infty$, so that $\cC^{[k]}(\R)$ is a Fr\'{e}chet space.

Finally, the previous item implies that $\cC^{[k]}(\R)$ is an algebra under pointwise multiplication and that pointwise multiplication is (jointly) continuous;
and Proposition \ref{prop.babyIPTP}\ref{item.Banachstaralg} implies that $\|\overline{f}\|_{\cC^{[k]},r} = \|f\|_{\cC^{[k]},r}$ whenever $f \in C^k(\R)$ and $r > 0$, so that complex conjugation is a continuous $\ast$-operation on $\cC^{[k]}(\R)$.
\end{proof}

This and Example \ref{ex.polycC[k]} bring us to our space of interest.
Inspired by the proof of Theorem \ref{thm.diffmatfunc} (specifically, Lemma \ref{lem.polydense}), we make the following definition of noncommutative $C^k$ functions.

\begin{defi}[Noncommutative $C^k$ Functions]\label{def.NCk}
For $k \in \N \cup \{\infty\}$, define $NC^k(\R) \coloneqq \overline{\C[\lambda]} \subseteq \cC^{[k]}(\R)$, where the closure takes place in $\cC^{[k]}(\R)$.
We call members of $NC^k(\R)$ \textbf{noncommutative} $\boldsymbol{C^k}$ \textbf{functions} and \textbf{noncommutative smooth functions} when $k=\infty$.
\end{defi}

\begin{prop}\label{prop.NCk}
If $k \in \N \cup \{\infty\}$, then $NC^k(\R)$ is a Fr\'{e}chet $\ast$-algebra under pointwise operations.
\end{prop}
\begin{proof}
Since $\C[\lambda] \subseteq \cC^{[k]}(\R)$ is a $\ast$-subalgebra, $NC^k(\R)$ is a closed $\ast$-subalgebra of the Fr\'{e}chet $\ast$-algebra $\cC^{[k]}(\R)$.
Thus $NC^k(\R)$ is a Fr\'{e}chet $\ast$-algebra in its own right.
\end{proof}

The idea for the name of $NC^k(\R)$ comes from parallel work by Jekel, who brought to the author's attention that he worked in Section 18 of \cite{jekel} with an abstract analog of $NC^k(\R)$ defined via completion and using classical projective tensor powers of $C([-r,r])$ in place of the integral projective tensor powers of $\ell^{\infty}([-r,r],\cB_{[-r,r]})$.
Jekel notates his space of noncommutative $C^k$ functions as $C_{\text{nc}}^k(\R)$.

\subsection{Elementary examples}\label{sec.NCkexs1}

So far it is only clear that $\C[\lambda] \subseteq NC^k(\R)$.
In this section, we use elementary harmonic analysis to exhibit more examples of noncommutative $C^k$ functions.
We begin with the \textit{Wiener spaces}.

\begin{nota}
Fix $m,k \in \N$.
\begin{enumerate}[label=(\alph*),leftmargin=2\parindent]
    \item $M(\R,\cB_{\R})$ is the space of Borel complex measures on $\R$.
    For $\mu \in M(\R,\cB_{\R})$, $|\mu|$ is the total variation measure of $\mu$, $\mu_{(0)} \coloneqq |\mu|(\R)$ is the total variation norm of $\mu$, and $\mu_{(k)} \coloneqq \int_{\R} |\xi|^k\,|\mu|(d\xi) \in [0,\infty]$ is the ``$k^{\text{th}}$ moment" of $|\mu|$.
    \item $BC(\R^m) = BC^0(\R^m) \coloneqq \{$bounded continuous functions $\R^m \to \C\}$ and
    \[
    BC^k(\R^m) \coloneqq \{f \in C^k(\R^m) : \partial^{\gamma}f \in BC(\R^m) \text{ whenever } \gamma \in \N_0^m \text{ and } |\gamma| \leq k\}
    \]
    with norm $\|f\|_{BC^k} \coloneqq \sum_{|\gamma| \leq k} \|\partial^{\gamma}f\|_{\ell^{\infty}(\R^m)}$.
    \item $\mathscr{S}(\R^m)$ is the Fr\'{e}chet space of Schwartz functions $\R^m \to \C$, and $\mathscr{S}'(\R^m) \coloneqq \mathscr{S}(\R^m)^*$ is the space of tempered distributions on $\R^m$.
    Also, the conventions we use for the Fourier transform and its inverse are, respectively,
    \[
    \wh{f}(\xi) = (\cF f)(\xi) = \int_{\R^m} e^{-ix \cdot \xi}f(x)\,dx \; \text{ and } \; \wch{f}(x) = \frac{1}{(2\pi)^m}\int_{\R^m}e^{ix\cdot\xi}f(\xi)\,d\xi
    \]
    for $f \in L^1(\R^m) = L^1(\R^m,\,\text{Lebesgue})$, with corresponding extensions to $\mathscr{S}'(\R^m)$.
\end{enumerate}
\end{nota}

\begin{defi}[Wiener Space]\label{def.Wk}
Fix $k \in \N_0$.
The $\boldsymbol{k^{\text{\textbf{th}}}}$ \textbf{Wiener space} $W_k(\R)$ is the set of functions $f \colon \R \to \C$ such that there exists (necessarily unique) $\mu \in M(\R,\cB_{\R})$ with $\mu_{(k)} < \infty$ and
\[
f(\lambda) = \int_{\R} e^{i\lambda \xi} \, \mu(d\xi),
\]
for all $\lambda \in \R$.
\end{defi}

Suppose that $\mu \in M(\R,\cB_{\R})$ and $\mu_{(k)} < \infty$ for some $k \in \N_0$.
If $0 \leq j \leq k$, then 
\begin{align*}
    \mu_{(j)} & = \int_{\R}|\xi|^j \, |\mu|(d\xi) \leq |\mu|([-1,1]) + \int_{\R} (1_{(-\infty,1)}(\xi) + 1_{(1,\infty)}(\xi)) |\xi|^j \, |\mu|(d\xi) \\
    & \leq |\mu|([-1,1]) + \int_{\R} (1_{(-\infty,1)}(\xi) + 1_{(1,\infty)}(\xi)) |\xi|^k \, |\mu|(d\xi) \leq |\mu|([-1,1]) + \mu_{(k)} < \infty.
\end{align*}
In particular, $W_k(\R) \subseteq W_j(\R)$.
Also, $W_k(\R) \subseteq BC^k(\R)$ because if $f = \int_{\R} e^{i\boldsymbol{\cdot}\xi}\,\mu(d\xi) \in W_k(\R)$, then
\[
f^{(j)}(\lambda) = \int_{\R}(i\xi)^je^{i\lambda \xi} \, \mu(d\xi), \numberthis\label{eq.Wkderform}
\]
for all $\lambda \in \R$;
thus $\|f\|_{BC^k} \leq \sum_{j=0}^k\mu_{(j)}$.
Here now is some important additional information about $W_k(\R)$.

\begin{lem}\label{lem.FourierFacts}
Fix $k \in \N$.
\begin{enumerate}[label=(\roman*),font=\normalfont,leftmargin=2\parindent]
    \item If $f \in BC(\R)$ and $\wh{f} \in L^1(\R)$, then\label{item.Wkcrit}
    \[
    \int_{\R} |\xi|^k\big|\wh{f}(\xi)\big|\,d\xi < \infty \iff f \in W_k(\R) \iff f \in C^k(\R) \text{ and } \wh{f^{(k)}} \in L^1(\R).
    \]
    \item If $f \in C^1(\R) \cap L^2(\R)$ and $f' \in L^2(\R)$, then $\wh{f} \in L^1(\R)$.\label{item.easyFourier}
    \item $C^{k+1}(\R) \subseteq W_k(\R)_{\loc}$.\label{item.Ck+1Wkloc}
\end{enumerate}
\end{lem}
\begin{proof}
We take each item in turn.

\ref{item.Wkcrit} Suppose $f \in BC(\R) \subseteq \mathscr{S}'(\R)$ and $\wh{f} \in L^1(\R)$.
By the Fourier Inversion Theorem for tempered distributions, the fact that $\wh{f} \in L^1(\R)$, and the continuity of $f$, we have
\[
f(\lambda) = \cF^{-1}\big(\wh{f}\,\big)(\lambda) = \frac{1}{2\pi}\int_{\R} e^{i\lambda\xi} \wh{f}(\xi)\,d\xi,
\]
for all $\lambda \in \R$.
Since $\mu(d\xi) \coloneqq \frac{1}{2\pi} \wh{f}(\xi) \,d\xi$ is a complex measure with $|\mu|(d\xi) = \frac{1}{2\pi} \big|\wh{f}(\xi)\big| \,d\xi$, we have
\[
\mu_{(k)} = \int_{\R} |\xi|^k\,|\mu|(d\xi) = \frac{1}{2\pi} \int_{\R} |\xi|^k\big|\wh{f}(\xi)\big|\,d\xi.
\]
The first equivalence immediately follows from this observation.
If $f \in C^k(\R)$, then $\big|\wh{f^{(k)}}(\xi)\big| = |\xi|^k\big|\wh{f}(\xi)\big|$, from which the second equivalence follows.

\ref{item.easyFourier} If $f \in C^1(\R) \cap L^2(\R)$ and $f' \in L^2(\R)$, then
\begin{align*}
    \big\|\wh{f}\big\|_{L^1} & = \int_{\R}\frac{1}{1+|\xi|} (1+|\xi|)\,\big|\wh{f}(\xi)\big| \, d\xi = \int_{\R}\frac{1}{1+|\xi|} \big(\big|\wh{f}(\xi)\big| + \big|\wh{f'}(\xi)\big|\big) \, d\xi \\
    & \leq \big\| (1+|\boldsymbol{\cdot}|)^{-1}\big\|_{L^2}\big(\big\|\wh{f}\big\|_{L^2} + \big\|\wh{f'}\big\|_{L^2}\big) = 2\sqrt{\pi}\big(\|f\|_{L^2} + \big\|f'\big\|_{L^2}\big) < \infty
\end{align*}
by the Cauchy--Schwarz Inequality and Plancherel's Theorem.

\ref{item.Ck+1Wkloc} Fix $f \in C^{k+1}(\R)$ and $r > 0$.
Let $\psi_r \in C_c^{\infty}(\R)$ be such that $\psi_r \equiv 1$ on $[-r,r]$.
We claim that $g \coloneqq \psi_r f \in W_k(\R)$.
Indeed, since $g \in C^{k+1}(\R)$ and $g$ has compact support, we have $g,g^{(k)} \in C^1(\R) \cap L^2(\R)$ and $g',g^{(k+1)} \in L^2(\R)$.
Thus $\cF(g),\cF(g^{(k)}) \in L^1(\R)$ by the previous item.
Since $g \in BC(\R)$ as well, the claim then follows from the first item.
But $g|_{[-r,r]} = f|_{[-r,r]}$ and $r > 0$ was arbitrary, so $f \in W_k(\R)_{\loc}$.
\end{proof}

We now get to work on our first classes of examples of noncommutative $C^k$ functions.

\begin{lem}\label{lem.WkcCk}
If $k \in \N$, $0 \leq j \leq k$, and $f = \int_{\R} e^{i\boldsymbol{\cdot}\xi}\,\mu(d\xi) \in W_k(\R)$, then
\[
\big\|f^{[j]}\big\|_{\ell^{\infty}(\R,\cB_{\R})^{\iotimes(j+1)}} \leq \frac{\mu_{(j)}}{j!} < \infty. \numberthis\label{eq.WkcCk}
\]
In particular, $W_k(\R) \subseteq \cC^{[k]}(\R)$ and $C^{k+1}(\R) \subseteq \overline{W_k(\R)} \subseteq \cC^{[k]}(\R)$.
\end{lem}
\begin{proof}
Fix $f = \int_{\R} e^{i\boldsymbol{\cdot}\xi}\,\mu(d\xi) \in W_k(\R)$ and $j \in \{1,\ldots,k\}$.
By Proposition \ref{prop.divdiff}\ref{item.cont} and Eq. \eqref{eq.Wkderform},
\[
f^{[j]}(\blambda) = \int_{\Delta_j} f^{(j)}(\boldsymbol{t} \cdot \blambda) \, \rho_j(d\boldsymbol{t}) =
\int_{\Delta_j}\int_{\R}e^{it_1\lambda_1\xi}\cdots e^{it_{j+1}\lambda_{j+1}\xi} (i\xi)^j\,\mu(d\xi)\, \rho_j(d \boldsymbol{t}) , \numberthis\label{eq.divdiffWk}
\]
for all $\blambda = (\lambda_1,\ldots,\lambda_{j+1}) \in \R^{j+1}$.
This expression, together with the Radon--Nikodym Theorem and Fubini's Theorem, says that if $(\Sigma,\sH,\rho) \coloneqq (\Delta_j \times \R, \cB_{\Delta_j \times \R}, \rho_j \otimes |\mu|)$, $\varphi_{\ell}(\lambda_{\ell},\sigma) \coloneqq e^{it_{\ell}\lambda_{\ell}\xi}$ for $1 \leq \ell \leq j$, and
\[
\varphi_{j+1}(\lambda_{j+1},\sigma) \coloneqq e^{i t_{j+1}\lambda_{j+1}\xi}(i\xi)^j\frac{d\mu}{d|\mu|}(\xi)
\]
for $\sigma = (t_1,\ldots,t_{j+1},\xi) \in \Delta_j \times \R$, then $(\Sigma,\rho,\varphi_1,\ldots,\varphi_{j+1})$ is a $\ell^{\infty}$-IPD of $f^{[j]}$ that gives
\[
\big\|f^{[j]}\big\|_{\ell^{\infty}(\R,\cB_{\R})^{\iotimes(j+1)}} \leq \int_{\Delta_j \times \R} |\xi|^j \,\rho(d\sigma) = \rho_j(\Delta_j) \, \mu_{(j)} = \frac{\mu_{(j)}}{j!},
\]
as claimed.
This proves the containment $W_k(\R) \subseteq \cC^{[k]}(\R)$.
The containment $C^{k+1}(\R) \subseteq \overline{W_k(\R)} \subseteq \cC^{[k]}(\R)$ then follows from Lemma \ref{lem.FourierFacts}\ref{item.Ck+1Wkloc} and Proposition \ref{prop.cC[k]}\ref{item.Sloc}.
\end{proof}
\begin{rem}
We can extract from the proofs above that if $f \in C^{k+1}(\R)$ has compact support, then
\[
\big\|f^{[j]}\big\|_{\ell^{\infty}(\R,\cB_{\R})^{\iotimes(j+1)}} \leq \frac{1}{2\pi j!}\int_{\R} |\xi|^j\big|\wh{f}(\xi)\big|\,d\xi \leq \frac{1}{j!\sqrt{\pi}}\big(\big\|f^{(j)}\big\|_{L^2}+\big\|f^{(j+1)}\big\|_{L^2}\big)
\]
whenever $1 \leq j \leq k$.
\end{rem}

\begin{thm}\label{thm.WkNCk}
If $k \in \N$, then $\overline{W_k(\R)} = \overline{C^{k+1}(\R)} = NC^k(\R)$.
(To be clear, the closures in the previous sentence take place in $\cC^{[k]}(\R)$.)
\end{thm}
\begin{proof}
We know from Lemma \ref{lem.WkcCk} that $\C[\lambda] \subseteq C^{k+1}(\R) \subseteq \overline{W_k(\R)}$.
Thus
\[
NC^k(\R) \subseteq \overline{C^{k+1}(\R)} \subseteq \overline{W_k(\R)}.
\]
It therefore suffices to prove $W_k(\R) \subseteq \overline{\C[\lambda]} = NC^k(\R)$.
To this end, fix a function $f = \int_{\R} e^{i\boldsymbol{\cdot}\xi}\,\mu(d\xi) \in W_k(\R)$.

We first reduce to the case when $\supp |\mu|$ is compact.
For $n \in \N$, define $\mu_n(d\xi) \coloneqq 1_{[-n,n]}(\xi) \, \mu(d\xi)$ and
\[
f_n(\lambda) \coloneqq \int_{\R} e^{i\lambda\xi}\,\mu_n(d\xi) = \int_{\R} e^{i\lambda \xi}1_{[-n,n]}(\xi) \, \mu(d\xi).
\]
Then $f_n \in W_k(\R)$ and $\supp |\mu_n| \subseteq [-n,n]$.
Also, if $0 \leq j \leq k$, then applying Eq. \eqref{eq.WkcCk} to $f-f_n$ yields
\[
\sup_{r > 0}\big\|(f-f_n)^{[j]}\big\|_{r,j+1} \leq \big\|(f-f_n)^{[j]}\big\|_{\ell^{\infty}(\R,\cB_{\R})^{\iotimes(j+1)}} \leq \frac{1}{j!}\int_{\R}|\xi|^j(1-1_{[-n,n]}(\xi)) \, |\mu|(d\xi) \to 0
\]
as $n \to \infty$ by the Dominated Convergence Theorem.
In particular, $f_n \to f$ in $\cC^{[k]}(\R)$ as $n \to \infty$. It now suffices to assume $\supp |\mu|$ is compact.

Suppose $R > 0$ and $\supp |\mu| \subseteq [-R,R]$.
Then $\int_{\R} |f|\,d|\mu| \leq \mu_{(0)}\,\|f\|_{\ell^{\infty}([-R,R])}$, for all Borel measurable functions $f \colon \R \to \C$.
In particular, $\mu_{(m)} \leq R^m\mu_{(0)} < \infty$, for all $m \in \N$.
Therefore, we may define
\[
q_n(\lambda) \coloneqq \int_{\R}\sum_{m=0}^n\frac{(i\lambda \xi)^m}{m!} \, \mu(d\xi) = \sum_{m=0}^n\frac{(i\lambda)^m}{m!}\int_{\R} \xi^m \, \mu(d\xi) \in \C[\lambda],
\]
for all $n \in \N$.
We claim that $q_n \to f$ in $\cC^{[k]}(\R)$ as $n \to \infty$. Indeed, note
\[
f(\lambda) - q_n(\lambda) = \int_{\R}\sum_{m=n+1}^{\infty}\frac{(i\lambda \xi)^m}{m!} \, \mu(d\xi),
\]
so that if $r > 0$, then
\[
\|f-q_n\|_{\ell^{\infty}([-r,r])} \leq \int_{\R}\sum_{m=n+1}^{\infty}\frac{(r|\xi|)^m}{m!} |\mu|(d\xi) \leq \mu_{(0)} \sum_{m=n+1}^{\infty}\frac{(rR)^m}{m!} \to 0
\]
as $n \to \infty$.
Also, if $j \in \N$ and $\blambda = (\lambda_1,\ldots,\lambda_{j+1}) \in \R^{j+1}$, then we know from Example \ref{ex.polydivdiff} that
\[
q_n^{[j]}(\blambda) = \sum_{m=0}^n\frac{i^m}{m!}\int_{\R}\xi^m \, \mu(d\xi)\sum_{|\gamma| = m-j} \blambda^{\gamma}.
\]
Similarly, after noting that $f(\lambda) = \sum_{m=0}^{\infty} \frac{(i\lambda)^m}{m!}\int_{\R}\xi^m \, \mu(d\xi)$, we have
\[
f^{[j]}(\blambda) = \sum_{m=0}^{\infty}\frac{i^m}{m!}\int_{\R}\xi^m \, \mu(d\xi)\sum_{|\gamma| = m-j} \blambda^{\gamma}
\]
as well.
Therefore, using the fact that $\{\gamma \in \N_0^{j+1} : |\gamma| = m-j\}$ has $\binom{m}{m-j} \leq 2^m$ elements, we get
\[
\big\|(f-q_n)^{[j]}\big\|_{r,j+1} \leq \sum_{m=n+1}^{\infty}\binom{m}{m-j}\frac{r^{m-j}}{m!}\mu_{(m)} \leq \frac{\mu_{(0)}}{r^j}\sum_{m=n+1}^{\infty}\frac{(2rR)^m}{m!} \to 0
\]
as $n \to \infty$, for all $r > 0$.
In particular, $q_n \to f$ in $\cC^{[k]}(\R)$ as $n \to \infty$.
This completes the proof.
\end{proof}

\begin{cor}\label{cor.NCinf}
$NC^{\infty}(\R) = \cC^{[\infty]}(\R) = C^{\infty}(\R)$. \qed
\end{cor}

\subsection{Advanced examples}\label{sec.NCkexs2}

As we just saw, only elementary methods are required to prove $C^{k+1}(\R) \subseteq NC^k(\R)$.
However, $NC^k(\R)$ is much closer to $C^k(\R)$ than that.
In this section, we use more advanced harmonic analysis done by Peller \cite{peller1} to exhibit two classes of examples of noncommutative $C^k$ functions that illustrate this point more precisely.
We begin by defining Besov spaces and stating their relevant properties.
For (much) more information on Besov spaces, please see \cite{leoni,peetre,sawano,triebel1,triebel2}.

\begin{defi}[Besov Spaces]\label{def.Besov}
Fix $m \,\in\, \N$ and $\eta \,\in\, C_c^{\infty}(\R^m)$ such that $0 \leq \eta \leq 1$ everywhere, $\supp \eta \subseteq \{\xi \in \R^m : |\xi|_2 \leq 2\}$, and $\eta \equiv 1$ on $\{\xi \in \R^m : |\xi|_2 \leq 1\}$.
(Here and throughout, $|\cdot|_2$ is the Euclidean norm.)
For $j \in \Z$, define
\[
\eta_j(\xi) \coloneqq \eta(2^{-j}\xi) - \eta(2^{-j+1}\xi), \; \xi \in \R^m.
\]
Now, for $(s,p,q) \in \R \times [1,\infty]^2$ and $f \in \mathscr{S}'(\R^m)$, define
\begin{align*}
    \|f\|_{\dot{B}_q^{s,p}} & \coloneqq \big\|\big(2^{js}\|\wch{\eta}_j \ast f\|_{L^p}\big)_{j \in \Z}  \big\|_{\ell^q(\Z)} \in [0,\infty] \; \text{ and} \\
    \|f\|_{B_q^{s,p}} & \coloneqq \|\wch{\eta} \ast f\|_{L^p} + \big\|\big(2^{js}\|\wch{\eta}_j \ast f\|_{L^p}\big)_{j \in \N}  \big\|_{\ell^q(\N)} \in [0,\infty].
\end{align*}
We call
\[
\dot{B}_q^{s,p}(\R^m) \coloneqq \big\{f \in \mathscr{S}'(\R^m) : \|f\|_{\dot{B}_q^{s,p}} < \infty\big\}
\]
the \textbf{homogeneous $\boldsymbol{(s,p,q)}$-Besov space} and
\[
B_q^{s,p}(\R^m) \coloneqq \big\{f \in \mathscr{S}'(\R^m) : \|f\|_{B_q^{s,p}} < \infty\big\}
\]
the \textbf{inhomogeneous $\boldsymbol{(s,p,q)}$-Besov space}.
\end{defi}
\begin{rem}
Note that $\wch{\eta} \ast f, \wch{\eta}_j \ast f$ have compactly supported Fourier transforms and so, by the Paley--Wiener Theorem, are smooth;
it therefore makes sense to apply the $L^p$-norm to them.
Also, beware that the positions of $p$ and $q$ in the notation for $B_q^{s,p}(\R^m)$ and $\dot{B}_q^{s,p}(\R^m)$ vary in the literature.
\end{rem}
\pagebreak

Here are the properties of Besov spaces that we shall use.
Below, the symbol $\hookrightarrow$ indicates (as usual) continuous inclusion, and $\sim$ indicates equivalence of (possibly infinite) norms.

\begin{thm}[Properties of Besov Spaces]\label{thm.Besov}
Fix $s,s_1,s_2 \in \R$ and $p,q,q_1,q_2 \in [1,\infty]$.
\begin{enumerate}[label=(\roman*),font=\normalfont,leftmargin=2\parindent]
    \item $(B_q^{s,p}(\R^m),\|\cdot\|_{B_q^{s,p}})$ is a Banach space that is independent of the choice of $\eta$.\label{item.Bspindep}
    \item $s_1 > s_2 \Rightarrow B_{q_1}^{s_1,p}(\R^m) \hookrightarrow B_{q_2}^{s_2,p}(\R^m)$ and $q_1 < q_2 \Rightarrow B_{q_1}^{s,p}(\R^m) \hookrightarrow B_{q_2}^{s,p}(\R^m)$.\label{item.easyemb}
    \item If $s \geq 0$, then $B_1^{s,\infty}(\R^m) \hookrightarrow BC^{\lfloor s \rfloor}(\R^m)$.\label{item.BCembed}
    \item For $f \in \mathscr{S}'(\R^m)$, define $\|f\|_{\mathrm{h},B_q^{s,p}} \coloneqq \|f\|_{L^p} + \|f\|_{\dot{B}_q^{s,p}} \in [0,\infty]$.
    (Of course, we declare $\|f\|_{L^p} \coloneqq \infty$ if $f$ is not induced by a locally integrable function.)
    If $s > 0$, then $\|\cdot\|_{B_q^{s,p}}  \sim \|\cdot\|_{\mathrm{h},B_q^{s,p}}$ on $\mathscr{S}'(\R^m)$.
    In particular, $B_q^{s,p}(\R^m) = L^p(\R^m) \cap \dot{B}_q^{s,p}(\R^m)$ when $s > 0$.\label{item.hIBesov}
    \item For vector spaces $V$ and $W$, a function $g \colon V \to W$, and vectors $x,h \in V$, recursively define
    \begin{align*}
        \Delta_h^1g(x) & = \Delta_h g(x) \coloneqq g(x+h)-g(x) \; \text{ and} \\
        \Delta_h^kg(x) & \coloneqq \Delta_h(\Delta_h^{k-1}g)(x) \; \text{ for } k \geq 2.
    \end{align*}
    Now, suppose $s > 0$.
    For $f \in L_{\loc}^1(\R^m)$, define
    \[
    \|f\|_{\mathrm{cl},B_q^{s,p}} \coloneqq \begin{cases}
    \|f\|_{L^p} + \big(\int_{\R^m} |h|_2^{-sq-m}\|\Delta_h^{\lfloor s \rfloor + 1} f\|_{L^p}^q \,dh\big)^{\frac{1}{q}}, & \text{if } q < \infty \\
    \|f\|_{L^p} + \sup_{h \in \R^m \setminus \{0\}} |h|_2^{-s}\|\Delta_h^{\lfloor s \rfloor + 1} f\|_{L^p}, & \text{if } q=\infty.
    \end{cases}
    \]
    Then $B_q^{s,p}(\R^m) = \{f \in L_{\loc}^1(\R^m) : \|f\|_{\mathrm{cl},B_q^{s,p}} < \infty\}$ and $\|\cdot\|_{B_q^{s,p}}  \sim \|\cdot\|_{\mathrm{cl},B_q^{s,p}}$ on $\mathscr{S}'(\R^m) \cap L_{\loc}^1(\R^m)$.\label{item.clIBesov}
\end{enumerate}
\end{thm}
\begin{proof}
Items \ref{item.Bspindep} and \ref{item.easyemb} are proven in Sections 2.3.2 and 2.3.3 of \cite{triebel1}, item \ref{item.BCembed} is proven in Section 2.1.2.4 of \cite{sawano}, item \ref{item.hIBesov} is proven in Section 2.3.3 of \cite{triebel2}, and item \ref{item.clIBesov} is proven in  Section 2.5.12 of \cite{triebel1}.
\end{proof}
\begin{rem}\label{rem.Sobolev}
It is also the case that $B_{\min\{p,2\}}^{s,p}(\R^m) \hookrightarrow W^{s,p}(\R^m) \hookrightarrow B_{\max\{p,2\}}^{s,p}(\R^m)$ for $s \in \R$ and $1 < p < \infty$, where $W^{s,p}(\R^m) = L_s^p(\R^m) = H_p^s(\R^m)$ is the fractional Sobolev (Bessel potential) space.
Please see Sections 2.2.2, 2.3.2, and 2.5.6 of \cite{triebel1}.
Also, in Chapter 17 of \cite{leoni}, $B_q^{s,p}(\R^m)$ (with $s > 0$) is defined and studied using $\|\cdot\|_{\mathrm{cl},B_q^{s,p}}$.
The equivalence $\|\cdot\|_{\mathrm{cl},B_q^{s,p}} \sim \|\cdot\|_{\mathrm{h},B_q^{s,p}}$ is proven in Section 17.7 thereof.
\end{rem}

As might be clear from the introduction, the most important indices for us are $(s,p,q) = (k \in \N,\infty,1)$.
It turns out in this case that $\dot{B}_1^{k,\infty}(\R) \subseteq C^k(\R)$.
(Please see Section A.2 of \cite{nikitopoulosOpDer}.)
With this in mind, we now state an important result of Peller that we shall use to prove $\dot{B}_1^{k,\infty}(\R) \subseteq NC^k(\R)$.

\begin{thm}[Peller \cite{peller1}]\label{thm.Peller}
If $k \in \N$, then there is a constant $a_k < \infty$ such that for all $f \in \dot{B}_1^{k,\infty}(\R)$ with $f^{(k)} \in BC(\R)$, we have
\[
\big\|f^{[k]}\big\|_{\ell^{\infty}(\R,\cB_{\R})^{\iotimes(k+1)}} \leq a_k\big(\big\|f^{(k)}\big\|_{\ell^{\infty}(\R)}+\|f\|_{\dot{B}_1^{k,\infty}}\big) < \infty.
\]
In particular, by Theorem \ref{thm.Besov}\ref{item.BCembed}--\ref{item.hIBesov}, there exists a constant $c_k < \infty$ such that
\[
\big\|f^{[k]}\big\|_{\ell^{\infty}(\R,\cB_{\R})^{\iotimes(k+1)}} \leq c_k\|f\|_{B_1^{k,\infty}},
\]
for all $f \in B_1^{k,\infty}(\R)$.
\end{thm}

A slightly stronger form of this result is Theorem 5.5 in \cite{peller1} (or Theorem 2.2.1 in \cite{peller2}, or Theorem 4.3.13 in \cite{nikitopoulosOpDer}).
Appendix A of \cite{nikitopoulosOpDer} contains a detailed and self-contained proof.
With this result in hand, we now begin the proof that $\dot{B}_1^{k,\infty}(\R) \subseteq NC^k(\R)$.

\begin{lem}[Inhomogeneous Littlewood--Paley Decomposition]\label{lem.LPapprox}
Fix $(s,p,q) \in \R \times [1,\infty] \times [1,\infty)$.
If $f \in \dot{B}_q^{s,p}(\R^m)$ and
\[
(f_n)_{n \in \N} \coloneqq \Bigg(\wch{\eta} \ast f + \sum_{j=1}^n \wch{\eta}_j \ast f\Bigg)_{n \in \N}
\]
is the \textbf{inhomogeneous Littlewood--Paley sequence} of $f$, then $f-f_n \in B_q^{s,p}(\R^m)$, for all $n \in \N$, and $\|f-f_n\|_{B_q^{s,p}} \to 0$ as $n \to \infty$.
\end{lem}
\begin{proof}
If $n \in \N$, then $\eta+\sum_{j=1}^n\eta_j = \eta(2^{-n}\cdot)$ by definition, so that
\[
f_n = \Bigg(\wch{\eta} + \sum_{j=1}^n\wch{\eta}_j \Bigg) \ast f = \wch{\eta(2^{-n}\cdot)} \ast f.
\]
Since $\eta(2^{-n}\cdot) \equiv 1$ on $\{\xi \in \R^m : |\xi|_2 \leq 2^n\}$, we have that if $1 \leq j \leq n-1$, then
\[
\wch{\eta}_j \ast f_n = \wch{\eta}_j \ast f\; \text{ and } \;\wch{\eta} \ast f_n = \wch{\eta} \ast f,
\]
as can be seen by taking Fourier transforms of both sides and using the fact that $\eta_j$ is supported in the annulus $\{\xi \in \R^m : 2^{j-1} \leq |\xi|_2 \leq 2^{j+1}\}$.
Next, note that
\[
\wch{\eta(2^{-n}\cdot)} = (2^n)^m \wch{\eta}(2^n\cdot) \implies
\Big\|\wch{\eta(2^{-n}\cdot)}\Big\|_{L^1} =\|\wch{\eta}\|_{L^1}.
\]
Therefore, by Young's Convolution Inequality, $\|\chi \ast f_n\|_{L^p} \leq \|\wch{\eta}\|_{L^1}\|\chi \ast f\|_{L^p}$ for all $\chi \in \mathscr{S}(\R^m)$.
Applying this to $\chi = \wch{\eta}_j$ and using the definition of $\|\cdot\|_{B_q^{s,p}}$, we have
\[
\|f-f_n\|_{B_q^{s,p}} = \Bigg(\sum_{j=n}^{\infty}2^{jsq}\|\wch{\eta}_j \ast (f-f_n)\|_{L^p}^q\Bigg)^{\frac{1}{q}} \leq (1+\|\wch{\eta}\|_{L^1})\Bigg(\sum_{j=n}^{\infty}2^{jsq}\|\wch{\eta}_j \ast f\|_{L^p}^q\Bigg)^{\frac{1}{q}} \to 0
\]
as $n \to \infty$ because $f \in \dot{B}_q^{s,p}(\R^m)$.
\end{proof}

\begin{thm}\label{thm.BesovNCk}
If $k \in \N$, then $\dot{B}_1^{k,\infty}(\R) \subseteq NC^k(\R)$.
Moreover, if $f \in \dot{B}_1^{k,\infty}(\R)$ and $(f_n)_{n \in \N}$ is the inhomogeneous Littlewood--Paley sequence of $f$, then $f_n \to f$ in $NC^k(\R)$ as $n \to \infty$.
\end{thm}
\begin{proof}
First, by Theorem \ref{thm.Peller}, $B_1^{k,\infty}(\R) \subseteq \cC^{[k]}(\R)$ and if $0 \leq j \leq k$, then
\[
\sup_{r > 0}\big\|f^{[j]}\big\|_{r,j+1} \leq \big\|f^{[j]}\big\|_{\ell^{\infty}(\R,\cB_{\R})^{\iotimes(j+1)}} \leq c_j\|f\|_{B_1^{j,\infty}} \leq c_j \|f\|_{B_1^{k,\infty}} < \infty, \numberthis\label{eq.Besovbound}
\]
for all $f \in B_1^{k,\infty}(\R)$.
(The case $j=0$ actually comes from Theorem \ref{thm.Besov}\ref{item.BCembed}.)
Next, fix $f \in \dot{B}_1^{k,\infty}(\R)$, and let $(f_n)_{n \in \N}$ be the inhomogeneous Littlewood--Paley sequence of $f$.
Note that if $n \in \N$, then $f_n$ has a compactly supported Fourier transform.
Therefore, by the Paley--Wiener Theorem, $f_n \in C^{\infty}(\R)$.
In particular, $f_n \in NC^k(\R)$ by Theorem \ref{thm.WkNCk}.
Now, by Lemma \ref{lem.LPapprox} and Eq. \eqref{eq.Besovbound}, 
\[
\sup_{r > 0}\|f-f_n\|_{\cC^{[k]},r} \lesssim_k \|f-f_n\|_{B_1^{k,\infty}} \to 0
\]
as $n \to \infty$.
Thus $f=f_n+(f-f_n) \in \cC^{[k]}(\R)$ and $f_n \to f$ in $\cC^{[k]}(\R)$ as $n \to \infty$.
Since we already noted $f_n \in NC^k(\R)$, for all $n \in \N$, this completes the proof.
\end{proof}

We observe somewhat parenthetically that the containment $\dot{B}_1^{k,\infty}(\R) \subseteq NC^k(\R)$ generalizes the containment $W_k(\R) \subseteq NC^k(\R)$.
(It should be noted, however, that our proof of the former used the latter in a crucial way.)
This is because $W_k(\R) \subseteq B_1^{k,\infty}(\R)$.
Indeed, if $k \in \N$, $f = \int_{\R} e^{i\boldsymbol{\cdot}\xi}\,\mu(d\xi) \in W_k(\R)$, and $\chi \in \mathscr{S}(\R)$, then
\[
(\chi \ast f)(\lambda) = \int_{\R}\int_{\R}e^{i(\lambda-y)\xi} \chi(y) \, \mu(d\xi) \, dy = \int_{\R}e^{i\lambda \xi} \int_{\R}e^{-iy\xi} \chi(y)\, dy \, \mu(d\xi) = \int_{\R} e^{i\lambda \xi} \wh{\chi}(\xi) \, \mu(d\xi),
\]
for all $\lambda \in \R$, by definition of convolution and Fubini's Theorem.
In particular, $(\wch{\eta} \ast f)(\lambda) = \int_{\R} e^{i\lambda \xi} \eta(\xi) \, \mu(d\xi)$ and $(\wch{\eta}_j \ast f)(\lambda) = \int_{\R} e^{i\lambda \xi} \eta_j(\xi) \, \mu(d\xi)$, for all $j \in \N$.
It follows that
\begin{align*}
    \|f\|_{B_1^{k,\infty}} & = \|\wch{\eta} \ast f\|_{L^{\infty}} + \sum_{j=1}^{\infty}2^{jk}\|\wch{\eta}_j \ast f\|_{L^{\infty}} \leq \int_{\R}|\eta(\xi)| \, |\mu|(d\xi) + \sum_{j=1}^{\infty}2^{jk}\int_{\R}|\eta_j(\xi)| \, |\mu|(d\xi)\\
    & = \int_{\{\xi \in \R : |\xi| \leq 2\}} |\eta(\xi)| \, |\mu|(d\xi) + \sum_{j=1}^{\infty}\int_{\{\xi \in \R : 2^{j-1} \leq |\xi| \leq 2^{j+1}\}}2^{jk} |\eta_j(\xi)| \, |\mu|(d\xi) \\
    & \leq |\mu|([-2,2]) + 2^k\int_{\R}\sum_{j=1}^{\infty}1_{\{\xi \in \R : 2^{j-1} \leq |\xi| \leq 2^{j+1}\}}|\xi|^k \, |\mu|(d\xi) \leq \mu_{(0)} + 3\cdot 2^k\mu_{(k)} < \infty,
\end{align*}
as claimed.
\pagebreak

We end this section by defining the H\"{o}lder spaces, describing their relationship to the Besov spaces, and proving that $C_{\loc}^{k,\e}(\R) \subseteq NC^k(\R)$.
For more information about H\"{o}lder spaces, please see \cite{fiorenza}.

\begin{defi}[H\"{o}lder Spaces]\label{def.Holder}
Let $(X,d_X)$ and $(Y,d_Y)$ be metric spaces and $g \colon X \to Y$ be a function.
For $\e > 0$, define
\[
[g]_{C^{0,\e}(X;Y)} \coloneqq \sup\Bigg\{\frac{d_Y(g(x),g(y))}{d_X(x,y)^{\e}} : x,y \in X, \,x \neq y\Bigg\}.
\]
If $[g]_{C^{0,\e}(X;Y)} < \infty$, then we say $g$ is $\boldsymbol{\e}$\textbf{-H\"{o}lder continuous} and write $g \in C^{0,\e}(X;Y)$.
As usual, we shall leave $Y$ out of the notation when $Y=\C$.
Next, if $k \in \N_0$ and $g \in C^k(\R^m)$, then we define
\[
[g]_{C^{k,\e}} \coloneqq \sqrt{\sum_{|\alpha| = k}\big[\partial^{\alpha}g \big]_{C^{0,\e}(\R^m)}^2}
\]
and $C^{k,\e}(\R^m) \coloneqq \{g \in C^k(\R^m) : [g]_{C^{k,\e}} < \infty\}$.
Also, we write $\|g\|_{BC^{k,\e}} \coloneqq \|g\|_{BC^k} + [g]_{C^{k,\e}}$ for $g \in C^k(\R^m)$ and define $BC^{k,\e}(\R^m) \coloneqq \{g \in BC^k(\R^m) : [g]_{C^{k,\e}} < \infty\}$.
Finally, we define $C_{\loc}^{k,\e}(\R^m)$ to be the set of $g \in C^k(\R^m)$ such that $\partial^{\gamma}g|_{[-r,r]^m} \in C^{0,\e}([-r,r]^m)$, for all $r > 0$ and $\gamma \in \N_0^m$ with $|\gamma| = k$.
\end{defi}

If $\e > 1$ and $g \in C_{\loc}^{0,\e}(\R^m)$, then $g$ is constant.
In particular, if $\e > 1$, $k \in \N$, and $g \in C_{\loc}^{k,\e}(\R^m)$, then $g \in \C[\lambda_1,\ldots,\lambda_m]$.
Also, the use of the $\ell^2$-norm (as opposed to the $\ell^1$-norm) in the definition of $[\cdot]_{C^{k,\e}}$ is atypical.
We made this choice so that the proof of the following proposition is more pleasant --- specifically, so that Eq. \eqref{eq.keyholderestim} below holds.

\begin{prop}\label{prop.HolderBesov}
If $k \in \N_0$ and $0 < \e \leq 1$, then $BC^{k,\e}(\R^m) \hookrightarrow B_{\infty}^{k+\e,\infty}(\R^m)$.
\end{prop}
\begin{proof}
Fix $k \in \N_0$ and $0 < \e \leq 1$.
For $g = (g_1,\ldots,g_n) \in C^k(\R^m;\C^n)$, define $[g]_{C^{k,\e}} \coloneqq \big(\sum_{j=1}^n[g_j]_{C^{k,\e}}^2\big)^{\frac{1}{2}}$.
We claim that if $g \in C^k(\R^m;\C^n)$, then
\[
\sup_{x \in \R^m}\big|\Delta_h^{k+1}g(x)\big|_2 \leq [g]_{C^{k,\e}}|h|_2^{k+\e}, \numberthis\label{eq.keyholderestim}
\]
for all $h \in \R^m \setminus \{0\}$.
First, observe that if $g \in C^k(\R^m;\C^n)$ and $\nabla g = (\partial_{\ell}g_j)_{1 \leq j \leq n,1 \leq \ell \leq m} \in C^{k-1}(\R^m;\C^{n \times m})$, then $[\nabla g]_{C^{k-1,\e}} = [g]_{C^{k,\e}}$.
Now, we prove Eq. \eqref{eq.keyholderestim} by induction.
If $k=0$, then it is immediate from the definition.
Now, assume the desired result holds when $k_0 \geq 0$, and let $k \coloneqq k_0+1$.
If $g \in C^k(\R^m;\C^n)$, $x,h \in \R^m$, and $h \neq 0$, then by Taylor's Theorem (or the Fundamental Theorem of Calculus),
\[
\Delta_h^{k+1}g(x) = \Delta_h^kg(x+h) - \Delta_h^k g(x) = \int_0^1 \nabla (\Delta_h^k g)(x+th) h \, dt = \int_0^1 \Delta_h^k (\nabla g)(x+th) h \, dt,
\]
where the juxtapositions $\nabla (\Delta_h^k g)(x+th) h, \, \Delta_h^k (\nabla g)(x+th) h$ above are matrix multiplications.
It then follows from the inductive hypothesis, the Cauchy--Schwarz Inequality, and our initial observation that
\[
\big|\Delta_h^{k+1}g(x)\big|_2 \leq \int_0^1\big|\Delta_h^k (\nabla g)(x+th)\big|_2 \,|h|_2 \, dt \leq [\nabla g]_{C^{k-1,\e}} |h|_2^{k-1+\e}|h|_2 = [g]_{C^{k,\e}}|h|_2^{k+\e},
\]
as desired.

Next, suppose $0 < \e < 1$.
Then $\lfloor k+\e \rfloor = k$, so Eq. \eqref{eq.keyholderestim} gives
\begin{align*}
    \|f\|_{\mathrm{cl},B_{\infty}^{k+\e,\infty}} & = \|f\|_{L^{\infty}} + \sup_{h \in \R^m \setminus \{0\}} |h|_2^{-k-\e}\|\Delta_h^{k+1}f\|_{L^{\infty}} \leq \|f\|_{L^{\infty}} + [f]_{C^{k,\e}} < \infty,
\end{align*}
for all $f \in BC^{k,\e}(\R^m)$.
Now, if $\e = 1$, then $\lfloor k+\e \rfloor = k+1$.
Combining Eq. \eqref{eq.keyholderestim} with the obvious fact that $\|\Delta_h^{k+2}f\|_{L^{\infty}} \leq 2\|\Delta_h^{k+1}f\|_{L^{\infty}}$ then gives
\begin{align*}
    \|f\|_{\mathrm{cl},B_{\infty}^{k+1,\infty}} & \leq \|f\|_{L^{\infty}} + 2\sup_{h \in \R^m \setminus \{0\}} |h|_2^{-k-1}\|\Delta_h^{k+1}f\|_{L^{\infty}} \leq \|f\|_{L^{\infty}} + 2 [f]_{C^{k,1}} < \infty,
\end{align*}
for all $f \in BC^{k,1}(\R^m)$.
An appeal to Theorem \ref{thm.Besov}\ref{item.clIBesov} completes the proof.
\end{proof}
\begin{rem}
In fact, $BC^{k,\e}(\R^m) = B_{\infty}^{k+\e,\infty}(\R^m)$ when $0 < \e < 1$ and $k \in \N_0$.
In general, the space $B_{\infty}^{s,\infty}(\R^m)$ is the \textit{H\"{o}lder--Zygmund space} $\mathcal{C}^s(\R^m)$ when $s > 0$.
For more information, please see Sections 2.2.2, 2.3.5, 2.5.7, and 2.5.12 of \cite{triebel1};
Sections 1.2.2, 1.5.1, and 2.6.5 of \cite{triebel2};
or Section 2.2.2 of \cite{sawano}.
\end{rem}
\pagebreak

As a consequence, we obtain the inclusion $C_{\loc}^{k,\e}(\R) \subseteq NC^k(\R)$.

\begin{thm}\label{thm.HolderNCk}
If $k \in \N$ and $\e > 0$, then $C_{\loc}^{k,\e}(\R) \subseteq NC^k(\R)$.
\end{thm}
\begin{proof}
Fix $\e,r > 0$, $k \in \N$, and $f \in C_{\loc}^{k,\e}(\R)$.
If $\e > 1$, then $f \in \C[\lambda] \subseteq NC^k(\R)$, so we assume $0 < \e \leq 1$.
Now, let $\psi_r \in C_c^{\infty}(\R)$ be such that $\psi_r \equiv 1$ on $[-r,r]$.
Then
\[
\psi_r f \in BC^{k,\e}(\R) \subseteq B_{\infty}^{k+\e,\infty}(\R) \subseteq B_1^{k,\infty}(\R) \subseteq \dot{B}_1^{k,\infty}(\R)
\]
by Proposition \ref{prop.HolderBesov} and Theorem \ref{thm.Besov}\ref{item.easyemb},\ref{item.hIBesov}.
Since $(\psi_r f)|_{[-r,r]} = f|_{[-r,r]}$ and $r > 0$ was arbitrary, we get
\[
f \in B_1^{k,\infty}(\R)_{\loc} \subseteq \dot{B}_1^{k,\infty}(\R)_{\loc} \subseteq \overline{\dot{B}_1^{k,\infty}(\R)} \subseteq NC^k(\R)
\]
by Theorem \ref{thm.BesovNCk} and Proposition \ref{prop.cC[k]}\ref{item.Sloc}.
\end{proof}

\subsection{Demonstration that \texorpdfstring{$W_k(\R)_{\loc} \subsetneq NC^k(\R)$}{}}\label{sec.Wkloc}

The formula \eqref{eq.divdiffWk} for the divided difference(s) of a function in $W_k(\R)$ is quite easy to work with, so it is reasonable to ask whether all examples of interest can be dealt with by ``localizing" $W_k(\R)$.
We already saw (Theorem \ref{thm.WkNCk}) that $W_k(\R)$ is dense in $NC^k(\R)$, but what we are really asking is whether the stronger statement $W_k(\R)_{\loc} = NC^k(\R)$ holds as well.
The goal of this section is to prove that this is not the case.

\begin{thm}\label{thm.Wkloccex}
If $k \in \N$, then $W_k(\R)_{\loc} \subsetneq NC^k(\R)$.
Specifically, we have the following counterexample.
Fix $\psi \in C_c^{\infty}(\R)$ such that $\psi \equiv 1$ on $[-1,1]$ and $\supp \psi \subseteq [-2,2]$, and define
\[
\kappa(x) \coloneqq 1_{(0,\infty)}(x) \,\psi(x) \sqrt{x}\,e^{-\frac{i}{x}}, \; x \in \R.
\]
If $f \in C^k(\R)$ and $f^{(k)} = \kappa$, then $f \in C^{k,\frac{1}{4}}(\R) \setminus W_k(\R)_{\loc} \subseteq NC^k(\R) \setminus W_k(\R)_{\loc}$.
\end{thm}

We break the proof into a few lemmas.
Fix $k \in \N$.

\begin{lem}\label{lem.Wkloc}
$W_k(\R)_{\loc} = \{f \in C^k(\R) : \eta f \in W_k(\R)$, for all $\eta \in C_c^{\infty}(\R)\}$.
\end{lem}
\begin{proof}
We first observe that if $f = \int_{\R} e^{i\boldsymbol{\cdot}\xi}\,\mu(d\xi) \in W_k(\R)$ and $\eta \in C_c^{\infty}(\R)$, then $\eta f \in W_k(\R)$.
Indeed, note
\[
\cF(\eta f)(\xi) = \int_{\R} e^{-ix\xi} \eta(x)f(x)\,dx = \int_{\R} \int_{\R} e^{-ix(\xi-y)} \eta(x)\,\mu(dy)\,dx = \int_{\R} \wh{\eta}(\xi-y)\,\mu(dy),
\]
for all $\xi \in \R$, by Fubini's Theorem.
Thus $\cF(\eta f) \in L^1(\R)$ with $\|\cF(\eta f)\|_{L^1} \leq \mu_{(0)}\|\widehat{\eta}\|_{L^1}$.
In addition,
\begin{align*}
    \int_{\R} |\xi|^k|\cF(\eta f)(\xi)|\,d\xi & \leq \int_{\R}\int_{\R}|\xi|^k|\wh{\eta}(\xi-y)|\,|\mu|(dy)\,d\xi = \int_{\R} \int_{\R} |\xi|^k|\wh{\eta}(\xi-y)|\,d\xi \,|\mu|(dy) \\
    & = \int_{\R} \int_{\R} |\zeta+y|^k|\wh{\eta}(\zeta)|\,d\zeta \,|\mu|(dy) \leq 2^k\int_{\R}\int_{\R} (|\zeta|^k+|y|^k)\,|\wh{\eta}(\zeta)| \,d\zeta \,|\mu|(dy) \\
    & = 2^k\big(\mu_{(0)}\,\big\|\cF\big(\eta^{(k)}\big)\big\|_{L^1} + \mu_{(k)}\,\|\wh{\eta}\|_{L^1}\big) < \infty
\end{align*}
by Tonelli's Theorem.
It follows from Lemma \ref{lem.FourierFacts}\ref{item.Wkcrit} that $\eta f \in W_k(\R)$.

Next, fix $f \in W_k(\R)_{\loc}$ and $\eta \in C_c^{\infty}(\R)$.
Suppose $\supp \eta \subseteq [-r,r]$ with $r > 0$.
By definition of $W_k(\R)_{\loc}$, there exists $g \in W_k(\R)$ such that $g|_{[-r,r]} = f|_{[-r,r]}$. But then $\eta f = \eta g \in W_k(\R)$ by the previous paragraph.
This proves $W_k(\R)_{\loc} \subseteq \{f \in C^k(\R) : \eta f \in W_k(\R)$, for all $\eta \in C_c^{\infty}(\R)\}$.

Finally, suppose that $f \in C^k(\R)$ is such that $\eta f \in W_k(\R)$, for all $\eta \in C_c^{\infty}(\R)$.
If $r > 0$, then let $\eta \in C_c^{\infty}(\R)$ be such that $\eta \equiv 1$ on $[-r,r]$.
Taking $g \coloneqq \eta f \in W_k(\R)$, we have $g|_{[-r,r]} = (\eta f)|_{[-r,r]} = f|_{[-r,r]}$.
We conclude that $f \in W_k(\R)_{\loc}$, which completes the proof.
\end{proof}

\begin{lem}\label{lem.Fourierint}
Fix $g \in C_c(\R)$ and $h \in C^k(\R)$ such that $h^{(k)} = g$.
Then
\[
h \in W_k(\R)_{\loc} \iff \wh{g} \in L^1(\R).
\]
\end{lem}
\begin{proof}
Let $\eta \in C_c^{\infty}(\R)$ be arbitrary.
Then $f_{\eta} \coloneqq \eta h \in C_c^k(\R)$.
In particular, by Lemma \ref{lem.FourierFacts}\ref{item.easyFourier}, $\wh{f_{\eta}} \in L^1(\R)$.
Now, by the product rule,
\[
f_{\eta}^{(k)} = \sum_{j=0}^k \binom{k}{j} \eta^{(j)}h^{(k-j)} = \eta\,h^{(k)} + \underbrace{\sum_{j=1}^k\binom{k}{j} \eta^{(j)}h^{(k-j)}}_{\coloneqq \chi} = \eta \, g+\chi.
\]
Because no more than $k-1$ derivatives fall on $h$ in the definition of $\chi$, we have $\chi \in C^1(\R)$.
Since $\chi$ has compact support, Lemma \ref{lem.FourierFacts}\ref{item.easyFourier} gives $\wh{\chi} \in L^1(\R)$.
It then follows from Lemma \ref{lem.FourierFacts}\ref{item.Wkcrit} that
\[
f_{\eta} = \eta h \in W_k(\R) \iff \cF\big(f_{\eta}^{(k)}\big) \in L^1(\R) \iff \cF(\eta g) \in L^1(\R).
\]
We combine this observation with the characterization of $W_k(\R)_{\loc}$ in Lemma \ref{lem.Wkloc} to finish the proof.
Suppose $h \in W_k(\R)_{\loc}$, and choose $\eta \in C_c^{\infty}(\R)$ such that $\eta \equiv 1$ on $\supp g$.
Then $\wh{g} = \cF(\eta g) \in L^1(\R)$.
Now, suppose $\wh{g} \in L^1(\R)$, and let $\eta \in C_c^{\infty}(\R)$ be arbitrary.
Then $\cF(\eta g) = \frac{1}{2\pi}\wh{\eta} \ast \wh{g} \in L^1(\R)$ because $\wh{\eta} \in L^1(\R)$ and $\wh{g} \in L^1(\R)$.
Thus $\eta h \in W_k(\R)$.
Since $\eta \in C_c^{\infty}(\R)$ was arbitrary, we conclude that $h \in W_k(\R)_{\loc}$.
\end{proof}

\begin{lem}\label{lem.chirp}
If $\kappa \in C_c(\R)$ is as in Theorem \ref{thm.Wkloccex}, then $\wh{\kappa} \not\in L^1(\R)$.
\end{lem}
\begin{proof}
Let $\xi > 0$.
Then
\[
\wh{\kappa}(\xi) = \int_0^{\infty} e^{-i(x\xi+x^{-1})} \psi(x)\sqrt{x}\,dx = \xi^{-\frac{3}{4}}\int_0^{\infty} e^{-i\sqrt{\xi}(y+y^{-1})}\psi\big(\xi^{-\frac{1}{2}}y\big)\sqrt{y}\,dy \numberthis\label{eq.rescaleFT}
\]
by the change of variable $y \coloneqq \sqrt{\xi}\,x$.
We use the method of stationary phase to analyze the oscillatory integral on the right hand side.
First, we notice that the phase $\phi(y) \coloneqq y+y^{-1}$ ($y > 0$) has a unique critical (``stationary") point at $y=1$, and this critical point is non-degenerate because $\phi''(1) = 2 \neq 0$.
Next, let $\chi \in C_c^{\infty}(\R)$ be such that $\chi \equiv 1$ on $\big[\frac{3}{4},\frac{3}{2}\big]$ and $\supp \chi \subseteq \big[\frac{1}{2},2\big]$.
Then
\[
I_1(\zeta) \coloneqq \int_0^{\infty} e^{-i\zeta \phi(y)}\chi(y)\,\psi\big(\zeta^{-1}y\big)\sqrt{y}\,dy = \int_0^{\infty} e^{-i\zeta \phi(y)}\chi(y)\sqrt{y}\,dy
\]
whenever $\zeta \geq 2$ because $\psi \equiv 1$ on $[0,1]$. Therefore, by Theorem 7.7.5 (and Equation (3.4.6)) in \cite{hormander},
\[
I_1(\zeta) = \chi(1) \sqrt{1}\,e^{-i\zeta \phi(1)-i\,\mathrm{sgn}(\phi''(1))\frac{\pi}{4}} \sqrt{\frac{2\pi}{\zeta |\phi''(1)|}} + O(\zeta^{-1}) = \sqrt{\pi}\,e^{-i(2\zeta +\frac{\pi}{4})}\zeta^{-\frac{1}{2}} + O(\zeta^{-1}) \, \text{ as } \zeta \to \infty.\numberthis\label{eq.I1assymp}
\]
Second, note that $\phi'(y) \neq 0$ for $0 < y \in \supp (1-\chi)$.
One can therefore apply the ``method of nonstationary phase" (integration by parts) to prove that
\[
I_2(\zeta) \coloneqq \int_0^{\infty} e^{-i\zeta \phi(y)}(1-\chi(y))\,\psi\big(\zeta^{-1}y\big)\sqrt{y}\,dy = O\big(\zeta^{-1}\big) \, \text{ as } \zeta \to \infty. \numberthis\label{eq.I2assymp}
\]
Due to the singularities of $\phi$ and the square root function at zero, standard theorems do not apply directly, so we need to prove this by hand.
The calculations necessary to do so are elementary but rather tedious, so we relegate them to Section \ref{sec.IBP}.
In the end, combining Eqs. \eqref{eq.rescaleFT}--\eqref{eq.I2assymp} gives
\[
\wh{\kappa}(\xi) = \xi^{-\frac{3}{4}}\big(I_1\big(\xi^{\frac{1}{2}}\big)+I_2\big(\xi^{\frac{1}{2}}\big)\big) = \sqrt{\pi}\,e^{-i(2\sqrt{\xi} +\frac{\pi}{4})}\xi^{-1} + O\big(\xi^{-\frac{5}{4}}\big) \, \text{ as } \xi \to \infty.
\]
It follows that $\wh{\kappa} \not\in L^1(\R)$, as claimed.
\end{proof}

\begin{proof}[Proof of Theorem \ref{thm.Wkloccex}]
It is an elementary exercise to show that $\kappa \in C^{0,\frac{1}{4}}(\R)$.
(For instance, one can adapt the argument from Example 1.1.8 in \cite{fiorenza}.)
In particular, if $f \in C^k(\R)$ and $f^{(k)} = \kappa$, then $f \in C^{k,\frac{1}{4}}(\R)$.
Thus $f \in NC^k(\R)$ by Theorem \ref{thm.HolderNCk}.
But $f \not\in W_k(\R)_{\loc}$ by Lemmas \ref{lem.Fourierint} and \ref{lem.chirp}.
\end{proof}

The above development provides a recipe for constructing functions in $\mathcal{N}_k \coloneqq NC^k(\R) \setminus W_k(\R)_{\loc}$.
Indeed, any compactly supported $g \in C^{0,\e}(\R)$ with $\wh{g} \not\in L^1(\R)$ can be used to produce a function in $\mathcal{N}_k$ via Lemma \ref{lem.Fourierint};
J. Sterbenz suggested $g=\kappa$ as an example.
(In general, for such a $g$ to exist, one must have $\e \leq \frac{1}{2}$.
This can be proven using Remark \ref{rem.Sobolev} and an argument like the one in the proof of Lemma \ref{lem.FourierFacts}\ref{item.easyFourier}.)

\section{Differentiating operator functions}

\subsection{Fr\'{e}chet derivatives}\label{sec.frechder}

Here, we briefly review some definitions and facts about Fr\'{e}chet derivatives.
For the duration of this section, fix $k \in \N$ and normed vector spaces $V_1,\ldots,V_k,V,W$ over $\F \in \{\R,\C\}$.

\begin{nota}[Bounded Multilinear Maps]\label{nota.bddmultilin}
If $T \colon V_1 \times \cdots \times V_k \to W$ is a $k$-linear map, then we write
\[
\|T\|_{B_k(V_1 \times \cdots \times V_k;W)} \coloneqq \sup\{\|T(v_1,\ldots,v_k)\|_W : v_j \in V_j, \; \|v_j\|_{V_j}\leq 1, \; 1 \leq j \leq k\} \in [0,\infty]
\]
for the operator norm of $T$ and $B_k(V_1 \times \cdots \times V_k;W)$ for the space of $k$-linear maps $V_1 \times \cdots \times V_k \to W$ with finite operator norm.
As usual, $(B(V_1;W),\|\cdot\|_{V_1 \to W}) \coloneqq (B_1(V_1;W),\|\cdot\|_{B_1(V_1;W)})$ and $B(W) \coloneqq B(W;W)$.
\end{nota}

Note that $B_k(V_1 \times \cdots \times V_k ; W) \cong B(V_1;B_{k-1}(V_2 \times \cdots \times V_k;W))$ isometrically via the map
\[
T \mapsto (v_1 \mapsto ((v_2,\ldots,v_k) \mapsto T(v_1,\ldots,v_k))).
\]
(In particular, by induction, if $W$ is a Banach space, then so is $B_k(V_1 \times \cdots \times V_k;W)$.)
We shall use this identification in the definition below.

\begin{defi}[Fr\'{e}chet Derivatives]\label{def.frechder}
Fix an open set $U \subseteq V$ and a map $F \colon U \to W$.
For $p \in U$, we say that $F$ is \textbf{Fr\'{e}chet differentiable at} $\boldsymbol{p}$ if there exists (necessarily unique) $DF(p) \in B(V;W)$ such that
\[
\frac{\|F(p+h)-F(p)-DF(p)h\|_W}{\|h\|_V} \to 0
\]
as $h \to 0$ in $V$.
If $F$ is Fr\'{e}chet differentiable at all $p \in U$, then we say $F$ is \textbf{Fr\'{e}chet differentiable in} $\boldsymbol{U}$ and write $D^1F = DF \colon U \to B(V;W)$ for its \textbf{Fr\'{e}chet derivative} map $U \ni p \mapsto DF(p) \in B(V;W)$.
For $k\geq 2$, we say that $F$ is $\boldsymbol{k}$\textbf{-times Fr\'{e}chet differentiable at} $\boldsymbol{p}$ if it is $(k-1)$-times Fr\'{e}chet differentiable in a neighborhood --- for simplicity, say $U$ --- of $p$ and $D^{k-1}F \colon U \to B_{k-1}(V^{k-1};W)$ is Fr\'{e}chet differentiable at $p$.
In this case, we write
\[
D^kF(p) \coloneqq D(D^{k-1}F)(p) \in  B(V;B_{k-1}(V^{k-1};W)) \cong B_k(V^k;W).
\]
If $F$ is $k$-times Fr\'{e}chet differentiable at all $p \in U$, then we say $F$ is $\boldsymbol{k}$\textbf{-times Fr\'{e}chet differentiable in} $\boldsymbol{U}$ and write $D^kF \colon U \to B_k(V^k;W)$ for its $\boldsymbol{k^{\text{\textbf{th}}}}$ \textbf{Fr\'{e}chet derivative} map $U \ni p \mapsto D^kF(p) \in B_k(V^k;W)$.
Finally, if in addition $D^kF$ is continuous, then we say $F$ is $\boldsymbol{k}$\textbf{-times continuously differentiable in} $\boldsymbol{U}$ and write $F \in C^k(U;W)$.
As usual, we also write $C^{\infty}(U;W) \coloneqq \bigcap_{k \in \N} C^k(U;W)$.
\end{defi}

Concretely, if $F \colon U \to W$ is $k$-times Fr\'{e}chet differentiable, then it is easy to show by induction that
\[
D^kF(p)[h_1,\ldots,h_k] = \partial_{h_1}\cdots \partial_{h_k}F(p) = \frac{d}{ds_1}\Big|_{s_1=0}\cdots \frac{d}{ds_k}\Big|_{s_k=0} F(p+s_1h_1+\cdots+s_kh_k).
\]
The only nontrivial fact about Fr\'{e}chet derivatives that we shall need is a certain completeness property, which we formulate and prove presently.

\begin{defi}
We define $BC_{\loc}(V;W) = BC_{\loc}^0(V;W)$ to be the space of continuous maps $V \to W$ that are bounded on bounded subsets of $V$.
For $k \in \N \cup \{\infty\}$, we define
\[
BC_{\loc}^k(V;W) \coloneqq \{F \in C^k(V;W) : D^jF \in BC_{\loc}(V;B_j(V^j;W)), \; 0 \leq j < k+1\}.
\]
Finally, if $k \in \N_0 \cup \{\infty\}$, then we define the \textbf{$\boldsymbol{BC_{\loc}^k}$ topology} on $BC_{\loc}^k(V;W)$ to be the locally convex topology induced by the family
\[
\Big\{F \mapsto \sup_{\|p\|_V \leq r}\|D^jF(p)\|_{B_j(V^j;W)} : 0 \leq j < k+1, \; r > 0\Big\}
\]
of seminorms, where $B_0(V^0;W) \coloneqq W$ and $D^0F \coloneqq F$.
\end{defi}

\begin{prop}[Completeness of $BC_{\loc}^k(V;W)$]\label{prop.Frechcomp}
If $W$ is a Banach space and $k \in \N_0 \cup \{\infty\}$, then $BC_{\loc}^k(V;W)$ is a Fr\'{e}chet space under the $BC_{\loc}^k$ topology.
\end{prop}
\begin{proof}
First, note that the $BC_{\loc}^k$ topology is generated by the countable family
\[
\Big\{F \mapsto \sup_{\|p\|_V \leq N}\|D^jF(p)\|_{B_j(V^j;W)} : 0 \leq j < k+1, \; N \in \N\Big\}
\]
of seminorms.
Since this family clearly separates points, we know that $BC_{\loc}^k(V;W)$ is a metrizable LCTVS.
It remains to prove that $BC_{\loc}^k(V;W)$ is complete.

We begin by proving that $BC_{\loc}(V;W)$ is complete.
To this end, let $(F_n)_{n \in \N}$ be a Cauchy sequence in $BC_{\loc}(V;W)$.
Note that if $p \in V$, then the linear map $BC_{\loc}(V;W) \ni F \mapsto F(p) \in W$ is continuous.
Therefore, $(F_n(p))_{n \in \N}$ is Cauchy in $W$.
Since $W$ is complete, $(F_n(p))_{n \in \N}$ converges to some $F(p) \in W$.
We claim that $F \in BC_{\loc}(V;W)$ and that $F_n \to F$ in the $BC_{\loc}^0$ topology as $n \to \infty$.
In fact, the former follows from the latter because the latter is precisely the statement that for all $r > 0$, the sequence $(F_n|_{\{q \in V : \|q\|_V \leq r\}})_{n \in \N}$ of bounded continuous functions converges uniformly to $F|_{\{q \in V : \|q\|_V \leq r\}}$;
and uniform limits of sequences of bounded continuous functions are themselves bounded and continuous.
To see that $F_n \to F$ in the $BC_{\loc}^0$ topology as $n \to \infty$, let $r > 0$ and $p \in V$ be such that $\|p\|_V \leq r$.
Then 
\[
\|F_n(p) - F(p)\|_W = \lim_{m \to \infty}\|F_n(p)-F_m(p)\|_W \leq \sup_{m \geq n}\sup_{\|q\|_V \leq r}\|F_n(q)-F_m(q)\|_W.
\]
Thus
\[
\limsup_{n \to \infty}\sup_{\|p\|_V \leq r}\|F_n(p)-F(p)\|_W \leq \lim_{n,m \to \infty}\sup_{\|p\|_V \leq r}\|F_n(p)-F_m(p)\|_W = 0
\]
because $(F_n)_{n \in \N}$ is Cauchy in $BC_{\loc}(V;W)$.
This completes the proof that $BC_{\loc}(V;W)$ is a Fr\'{e}chet space.

Next, suppose that $k \in \N \cup \{\infty\}$ and $(F_n)_{n \in \N}$ is a Cauchy sequence in $BC_{\loc}^k(V;W)$.
If $0 \leq j < k+1$, then the sequence $(D^jF_n)_{n \in \N}$ is Cauchy in $BC_{\loc}(V;B_j(V^j;W))$.
By the previous paragraph, there exists $G_j \in BC_{\loc}(V;B_j(V^j;W))$ such that $D^jF_n \to G_j$ in $BC_{\loc}(V;B_j(V^j;W))$ as $n \to \infty$.
By Theorem 85 in Section 1.3 of \cite{hajek}, it follows that $G_0 \in C^k(V;W)$ and $D^jG_0 = G_j$ whenever $0 \leq j < k+1$.
Unraveling the definitions, we conclude that $F_n \to G_0$ in the $BC_{\loc}^k$ topology as $n \to \infty$.
This completes the proof.
\end{proof}

\subsection{Multiple operator integrals (MOIs)}\label{sec.MOIs}

For the duration of this section, fix a complex Hilbert space $H$, a von Neumann algebra $\cM \subseteq B(H)$, and a natural number $k \in \N$.
(Readers less familiar with von Neumann algebras may take $\cM = B(H)$ throughout.)
We now describe the ``separation of variables" approach, from \cite{azamovetal,nikitopoulosMOI,peller1,peller2}, to defining the multiple operator integral (MOI)
\[
\big(I^{\boldsymbol{a}}\varphi\big)[b] = \int_{\sigma(a_{k+1})}\cdots\int_{\sigma(a_1)} \varphi(\blambda)\,P^{a_1}(d\lambda_1)\,b_1\cdots P^{a_k}(d\lambda_k)\,b_k\,P^{a_{k+1}}(d\lambda_{k+1}) \in \cM
\]
for $\boldsymbol{a} = (a_1,\ldots,a_{k+1}) \in \cM_{\nu}^{k+1}$, $b = (b_1,\ldots,b_k) \in \cM^k$, and
\[
\varphi \in \ell^{\infty}\big(\sigma(a_1),\cB_{\sigma(a_1)}\big) \iotimes \cdots \iotimes \ell^{\infty}\big(\sigma(a_{k+1}),\cB_{\sigma(a_{k+1})}\big).
\]
Heuristically, if $(\Sigma,\rho,\varphi_1,\ldots,\varphi_{k+1})$ is a $\ell^{\infty}$-IPD of $\varphi$, then we should have
\begin{align*}
    \big(I^{\boldsymbol{a}}\varphi\big)[b] & = \int_{\sigma(a_{k+1})}\cdots\int_{\sigma(a_1)} \int_{\Sigma} \varphi_1(\lambda_1,\sigma)\cdots\varphi_{k+1}(\lambda_{k+1},\sigma)\,\rho(d\sigma)\,P^{a_1}(d\lambda_1)\,b_1\cdots P^{a_k}(d\lambda_k)\,b_k\,P^{a_{k+1}}(d\lambda_{k+1}) \\
    & = \int_{\Sigma}\int_{\sigma(a_{k+1})}\cdots\int_{\sigma(a_1)}  \varphi_1(\lambda_1,\sigma)\cdots\varphi_{k+1}(\lambda_{k+1},\sigma)\,P^{a_1}(d\lambda_1)\,b_1\cdots P^{a_k}(d\lambda_k)\,b_k\,P^{a_{k+1}}(d\lambda_{k+1})\,\rho(d\sigma) \\
    & = \int_{\Sigma}\Bigg(\int_{\sigma(a_1)}  \varphi_1(\cdot,\sigma)\,dP^{a_1}\Bigg)\,b_1\cdots \Bigg(\int_{\sigma(a_k)}\varphi_k(\cdot,\sigma)\,dP^{a_k}\Bigg)\,b_k\Bigg(\int_{\sigma(a_{k+1})}\varphi_{k+1}(\cdot,\sigma)\,dP^{a_{k+1}}\Bigg)\,\rho(d\sigma) \\
    & = \int_{\Sigma} \varphi_1(a_1,\sigma)\,b_1\cdots \varphi_k(a_k,\sigma)\,b_k\,\varphi_{k+1}(a_{k+1},\sigma)\,\rho(d\sigma). \numberthis\label{eq.integ}
\end{align*}
Accordingly, we shall use Eq. \eqref{eq.integ} as a definition.
To do so, we must address exactly what kind of integral $\int_{\Sigma} \boldsymbol{\cdot} \,d\rho$ is being used above and whether this integral depends on the chosen $\ell^{\infty}$-IPD of $\varphi$.
\pagebreak

We begin with the former.
Let $(\Sigma,\sH,\rho)$ be a measure space.
A map $F \colon \Sigma \to B(H)$ is called \textbf{pointwise weakly measurable} if $\la F(\cdot)h_1,h_2 \ra \colon \Sigma \to \C$ is $(\sH,\cB_{\C})$-measurable, for all $h_1,h_2 \in H$.
Also, for $\cS \subseteq B(H)$, write $W^*(\cS) \subseteq B(H)$ for the smallest von Neumann algebra containing $\cS$.

\begin{lem}[Pointwise Pettis Integrability]\label{lem.integ}
Let $(\Sigma,\sH,\rho)$ be a measure space.
If $F \colon \Sigma \to B(H)$ is pointwise weakly measurable and $\int_{\Sigma} |\la F(\sigma)h_1,h_2 \ra|\,\rho(d\sigma) < \infty$, for all $h_1,h_2 \in H$, then there exists unique $T \in B(H)$ such that
\[
\la Th_1,h_2 \ra = \int_{\Sigma} \la F(\sigma)h_1,h_2 \ra\,\rho(d\sigma), \text{ for all } h_1,h_2 \in H.
\]
Moreover, $T \in W^*(F(\sigma) : \sigma \in \Sigma)$.
In this case, we say that $F$ is \textbf{pointwise Pettis ($\boldsymbol{\rho}$-)integrable}, we call the operator $T$ the \textbf{pointwise Pettis ($\boldsymbol{\rho}$-)integral} of $F$, and we write $\int_{\Sigma} F\,d\rho = \int_{\Sigma} F(\sigma)\,\rho(d\sigma) \coloneqq T$.
\end{lem}
\begin{proof}
Write $K$ for the complex conjugate space of $H$, i.e., the Hilbert space with the same underlying set and addition but with the scalar multiplication $c\cdot h\coloneqq \overline{c} \,h$ and the inner product $\la h_1,h_2 \ra_K \coloneqq \la h_2,h_1 \ra$.
By a standard Riesz Representation Theorem argument, the correspondence
\[
B(H) \ni T \mapsto ((h_1,h_2) \mapsto \la Th_1,h_2 \ra) \in B_2(H \times K;\C)
\]
is an isometric isomorphism.
The proof of the first claim will therefore be complete if we can show that the bilinear map $H \times K \ni (h_1,h_2) \mapsto \int_{\Sigma} \la F(\sigma)h_1,h_2 \ra \,\rho(d\sigma) \in \C$ is bounded.

To this end, define a bilinear map $I \colon H \times K \to L^1(\Sigma,\rho)$ by $(h_1,h_2) \mapsto \la F(\cdot)h_1,h_2 \ra$.
By a standard Principle of Uniform Boundedness argument, if we can show that for every fixed $(h_1,h_2) \in H \times K$, we have $I[h_1,\cdot] \in B(K;L^1(\Sigma,\rho))$ and $I[\cdot,h_2] \in B(H;L^1(\Sigma,\rho))$, then we get $I \in B_2(H \times K;L^1(\Sigma,\rho))$.
Accordingly, fix $h_1 \in H$, and suppose that $(k_n)_{n \in \N}$ is a sequence in $K$ converging to $k \in K$.
Now, if $\sigma \in \Sigma$, then $\la F(\sigma)h_1,k_n \ra \to \la F(\sigma)h_1,k \ra$ as $n \to \infty$.
In other words, $I[h_1,k_n] \to I[h_1,k]$ pointwise as $n \to \infty$.
In particular, if we assume that the sequence $(I[h_1,k_n])_{n \in \N}$ converges in $L^1(\Sigma,\rho)$, then it must converge to $I[h_1,k]$.
We have just proven that $I[h_1,\cdot]$ is closed.
By the Closed Graph Theorem, $I[h_1,\cdot] \in B(K;L^1(\Sigma,\rho))$.
The proof that $I[\cdot,h_2] \in B(H;L^1(\Sigma,\rho))$, for all $h_2 \in K$, is identical.
Finally, note that the integral map $L^1(\Sigma,\rho) \ni f \mapsto \int_{\Sigma} f \,d\rho \in \C$ is bounded.
Composing this integral map with the bounded bilinear map $I$ and unraveling the definitions yield the desired conclusion.

For the second claim, write $\mathcal{N} \coloneqq W^*(F(\sigma) : \sigma \in \Sigma) \subseteq B(H)$.
To see that $\int_{\Sigma} F\,d\rho \in \mathcal{N}$, fix $a \in \mathcal{N}'$, where $\mathcal{N}'$ is the commutant of $\mathcal{N}$.
Then $a\int_{\Sigma} F\,d\rho = \int_{\Sigma} a\,F(\sigma)\,\rho(d\sigma) = \int_{\Sigma} F(\sigma)\,a \,\rho(d\sigma) = \big(\int_{\Sigma} F\,d\rho\big)a$, as the reader may readily verify.
Thus $\int_{\Sigma} F\,d\rho \in \mathcal{N}'' = \mathcal{N}$ by von Neumann's Bicommutant Theorem.
\end{proof}
\begin{rem}
For an explanation of the terms \textit{pointwise weakly measurable} and \textit{pointwise Pettis integrable}, please see Section 3.3 of \cite{nikitopoulosMOI}.
\end{rem}

\begin{lem}\label{lem.wmeas}
Suppose $(\Sigma,\sH)$ and $(\Om,\sF)$ are measurable spaces, $P \colon \sF \to B(H)$ is a projection-valued measure, and $\varphi \colon \Om \times \Sigma \to \C$ is product measurable.
If $\varphi(\cdot,\sigma) \in \ell^{\infty}(\Om,\sF)$, for all $\sigma \in \Sigma$, and $F \colon \Sigma \to B(H)$ is pointwise weakly measurable, then
\[
F(\cdot)\int_{\Om} \varphi(\om,\cdot)\,P(d\om) \colon \Sigma \to B(H) \; \text{ and } \; \int_{\Om} \varphi(\om,\cdot)\,P(d\om)\,F(\cdot) \colon \Sigma \to B(H)
\]
are pointwise weakly measurable.
\end{lem}
\begin{proof}
For $\sigma \in \Sigma$ and $h_1,h_2 \in H$, define $\nu_{\sigma} \coloneqq \la F(\sigma)P(\cdot)h_1,h_2 \ra = \la P(\cdot)h_1,F(\sigma)^*h_2\ra = P_{h_1,F(\sigma)^*h_2}$.
Then $\nu_{\sigma}$ is a complex measure, and $\Sigma \ni \sigma \mapsto \nu_{\sigma}(G) \in \C$ is measurable, for all $G \in \sF$, because $F$ is pointwise weakly measurable.
Also, by definition of integration with respect to $P$,
\begin{align*}
    \Bigg\la F(\sigma)\int_{\Om} \varphi(\om,\sigma)\,P(d\om)h_1,h_2\Bigg\ra & = \Bigg\la \int_{\Om} \varphi(\om,\sigma)\,P(d\om)h_1,F(\sigma)^*h_2\Bigg\ra \\
    & = \int_{\Om} \varphi(\om,\sigma) \,P_{h_1,F(\sigma)^*h_2}(d\om) = \int_{\Om} \varphi(\om,\sigma)\,\nu_{\sigma}(d\om).
\end{align*}
By Lemma 4.2.2 in \cite{nikitopoulosMOI}, the function $\Sigma \ni \sigma \mapsto \int_{\Om} \varphi(\om,\sigma)\,\nu_{\sigma}(d\om) \in \C$ is measurable.
This proves that $F(\cdot)\int_{\Om}\varphi(\om,\cdot)\,P(d\om)$ is pointwise weakly measurable.
The argument for $\int_{\Om} \varphi(\om,\cdot)\,P(d\om)\,F(\cdot)$ is identical.
\end{proof}

\begin{thm}[Definition of MOIs]\label{thm.babyMOI}
Fix $\boldsymbol{a} = (a_1,\ldots,a_{k+1}) \in \cM_{\nu}^{k+1}$,
\[
\varphi \in \ell^{\infty}\big(\sigma(a_1),\cB_{\sigma(a_1)}\big) \iotimes \cdots \iotimes \ell^{\infty}\big(\sigma(a_{k+1}),\cB_{\sigma(a_{k+1})}\big),
\]
and $b = (b_1,\ldots,b_k) \in \cM^k$.
\begin{enumerate}[label=(\roman*),font=\normalfont,leftmargin=2\parindent]
    \item If $(\Sigma,\rho,\varphi_1,\ldots,\varphi_{k+1})$ is a $\ell^{\infty}$-IPD of $\varphi$, then the map
    \[
    \Sigma \ni \sigma \mapsto F(\sigma) \coloneqq \varphi_1(a_1,\sigma)\,b_1\cdots \varphi_k(a_k,\sigma)\,b_k\,\varphi_{k+1}(a_{k+1},\sigma) \in \cM
    \]
    is pointwise Pettis integrable.\label{item.integ}
    \item If $F$ is as in the previous item, then the pointwise Pettis integral
    \[
    \big(I^{\boldsymbol{a}}\varphi\big)[b] = \int_{\sigma(a_{k+1})}\cdots\int_{\sigma(a_1)} \varphi(\blambda)\,P^{a_1}(d\lambda_1)\,b_1\cdots P^{a_k}(d\lambda_k)\,b_k\,P^{a_{k+1}}(d\lambda_{k+1}) \coloneqq \int_{\Sigma} F\,d\rho \in \cM
    \]
    is independent of the chosen $\ell^{\infty}$-IPD of $\varphi$.\label{item.welldef}
    \item The assignment $\cM^k \ni b \mapsto (I^{\boldsymbol{a}}\varphi)[b] \in \cM$ is $k$-linear and bounded.
    Moreover, the assignment
    \[
   \ell^{\infty}\big(\sigma(a_1),\cB_{\sigma(a_1)}\big) \iotimes \cdots \iotimes \ell^{\infty}\big(\sigma(a_{k+1}),\cB_{\sigma(a_{k+1})}\big) \ni \varphi \mapsto I^{\boldsymbol{a}}\varphi \in B_k(\cM^k;\cM)
    \]
    is linear and has operator norm at most one.\label{item.bddlin}
\end{enumerate}
\end{thm}
\begin{proof}[Sketch of proof]
We take each item in turn.
We write $\|\cdot\| = \|\cdot\|_{H \to H}$ throughout.

\ref{item.integ} By Lemma \ref{lem.wmeas} and induction, $F \colon \Sigma \to B(H)$ is pointwise weakly measurable.
To prove integrability, fix $h_1,h_2 \in H$ and $\sigma \in \Sigma$.
Then
\[
|\la F(\sigma)h_1,h_2 \ra| \leq \|h_1\|\,\|h_2\|\Bigg(\prod_{\ell=1}^k\|b_{\ell}\|\Bigg)\prod_{j=1}^{k+1}\|\varphi_j(a_j,\sigma)\| \leq \|h_1\|\,\|h_2\|\Bigg(\prod_{\ell=1}^k\|b_{\ell}\|\Bigg)\prod_{j=1}^{k+1}\|\varphi_j(\cdot,\sigma)\|_{\ell^{\infty}(\sigma(a_j))}.
\]
Therefore,
\[
\int_{\Sigma}|\la F(\sigma)h_1,h_2 \ra| \,\rho(d\sigma) \leq \|h_1\|\,\|h_2\|\Bigg(\prod_{\ell=1}^k\|b_{\ell}\|\Bigg)\int_{\Sigma}\prod_{j=1}^{k+1}\|\varphi_j(\cdot,\sigma)\|_{\ell^{\infty}(\sigma(a_j))}\,\rho(d\sigma) < \infty. \numberthis\label{eq.MOIbd}
\]
Thus $F$ is pointwise Pettis integrable.
Also, $\int_{\Sigma} F\,d\rho \in W^*(F(\sigma) : \sigma \in \Sigma) \subseteq \cM$ by Lemma \ref{lem.integ}.

\ref{item.welldef} For this item, it suffices to assume $\cM = B(H)$.
First, suppose $h_1,\tilde{h}_1,\ldots,h_k,\tilde{h}_k \in H$ and $b_j = \la \cdot, h_j \ra \tilde{h}_j$ whenever $1 \leq j \leq k$.
Also, for $h_0,\tilde{h}_{k+1} \in H$, define
\[
\nu \coloneqq P_{\tilde{h}_1,h_0}^{a_1} \otimes \cdots \otimes P_{\tilde{h}_{k+1},h_k}^{a_{k+1}}.
\]
(This is a product of complex measures.)
Then one can show without much difficulty that
\[
\Bigg\la \Bigg(\int_{\Sigma} F\,d\rho\Bigg)\tilde{h}_{k+1},h_0 \Bigg \ra = \int_{\sigma(a_1) \times \cdots \times \sigma(a_{k+1})} \varphi \, d\nu. \numberthis\label{eq.indep}
\]
For this calculation or similar ones, please see the proof of Theorem 4.2.12 in \cite{nikitopoulosMOI}, the proof of Lemma 4.3 in \cite{azamovetal}, or the proof of Theorem 2.1.1 in \cite{peller2}.
From Eq. \eqref{eq.indep} and $k$-linearity, we conclude that $\int_{\Sigma} F\,d\rho$ is independent of the chosen $\ell^{\infty}$-IPD of $\varphi$ when $b_1,\ldots,b_k$ are finite-rank operators.
Now, if $H$ is separable and $c \in B(H)$, then $c$ is the strong operator limit of a sequence of finite-rank operators.
This allows one to use an operator-valued Dominated Convergence Theorem (e.g., from Section 2.1 of \cite{nikitopoulosOpDer}) to extend the claimed independence to arbitrary $b_1,\ldots,b_k \in B(H)$.
This is what is done in \cite{azamovetal,peller2}.
The case when $H$ is not separable, which is treated in \cite{nikitopoulosMOI}, requires much more care.
The claim is again extended from finite-rank to arbitrary bounded operators by density but in a different topology:
the ultraweak topology.
Indeed, one proves that, for fixed $b_1,\ldots,b_k \in B(H)$ and $j \in \{1,\ldots,k\}$, the assignment
\[
B(H) \ni c \mapsto \int_{\Sigma}\varphi_1(a_1,\sigma) \Bigg(\prod_{\ell=1}^{j-1} b_{\ell}\,\varphi_{\ell+1}(a_{\ell+1},\sigma)\Bigg)c\Bigg(\prod_{\ell=j+1}^k\varphi_{\ell}(a_{\ell},\sigma)\,b_{\ell}\Bigg)\varphi_{k+1}(a_{k+1},\sigma)\,\rho(d\sigma) \in B(H)
\]
is ultraweakly continuous.
(Above, empty products are declared to be $1$.)
Proving this is quite technical.
We refer the interested reader to Section 4.2, specifically Corollary 4.2.11 and its leadup, in \cite{nikitopoulosMOI} for the details.
\pagebreak

In the present setting, which is less general than that of \cite{nikitopoulosMOI}, we can employ a different argument to deduce the non-separable case from the separable case treated in \cite{azamovetal,peller2}.
First, suppose that $\cA \subseteq B(H)$ is a unital subalgebra.
We claim that if $\cA$ is SOT-separable and $h_1,\ldots,h_n \in H$, then there is a closed, separable linear subspace $K \subseteq H$ such that $h_1,\ldots,h_n\in K$ and $\cA K \subseteq K$, i.e., $K$ is $\cA$-invariant.
Indeed, define
\[
K \coloneqq \overline{\spn}\big(\cA h_1 \cup \cdots \cup \cA h_n\big) \subseteq H.
\]
Certainly, $K$ is a closed linear subspace of $H$ containing $h_1,\ldots,h_n$.
Also, $K$ is separable because if $\cA_0 \subseteq \cA$ is a countable SOT-dense subset, then the $\Q[i]$-span of $\cA_0 h_1 \cup \cdots \cup \cA_0h_n$ is dense in $K$.
Finally, $K$ is $\cA$-invariant because $\cA$ is a subalgebra and closed linear spans of $\cA$-invariant subsets are $\cA$-invariant.
Next, fix $h_1,h_2 \in H$, and apply the previous result to $\cA \coloneqq W^*(a_1,\ldots,a_{k+1},b_1,\ldots,b_k) \subseteq B(H)$ to obtain a closed, separable, $\cA$-invariant linear subspace $K \subseteq H$ that contains $h_1$ and $h_2$.
(Note that $\cA$ is SOT-separable because the $\Q[i]$-span of noncommutative monomials in $a_1,a_1^*,\ldots,a_{k+1},a_{k+1}^*,b_1,b_1^*,\ldots,b_k,b_k^*$ is SOT-dense in $\cA$.)
If we write $\pi_K \colon H \to K$ for the orthogonal projection onto $K$, $\iota_K \colon K \to H$ for the inclusion of $K$ into $H$, $\tilde{a}_j \coloneqq \pi_Ka_j\iota_K \in B(K)_{\nu}$ for $j \in \{1,\ldots,k+1\}$, and $\tilde{b}_{\ell} \coloneqq \pi_Kb_{\ell}\iota_K \in B(K)$ for $\ell \in \{1,\ldots,k\}$, then
\[
F(\sigma)h = \varphi_1\big(\tilde{a}_1,\sigma\big) \,\tilde{b}_1\cdots \varphi_k\big(\tilde{a}_k,\sigma\big)\,\tilde{b}_k\,\varphi_{k+1}\big(\tilde{a}_{k+1},\sigma\big)h, \; \text{ for all } \sigma \in \Sigma \text{ and } h \in K,
\]
as we encourage the reader to verify using the $\cA$-invariance of $K$.
Therefore,
\begin{align*}
    \Bigg\la \Bigg(\int_{\Sigma} F\,d\rho\Bigg)h_1,h_2 \Bigg\ra_H & = \int_{\Sigma} \la F(\sigma)h_1,h_2 \ra_H \, \rho(d\sigma) \\
    & = \int_{\Sigma}\big\la \varphi_1\big(\tilde{a}_1,\sigma\big) \,\tilde{b}_1\cdots \varphi_k\big(\tilde{a}_k,\sigma\big)\,\tilde{b}_k\,\varphi_{k+1}\big(\tilde{a}_{k+1},\sigma\big)h_1,h_2 \big\ra_K \, \rho(d\sigma) \\
    & = \Bigg\la \Bigg(\int_{\Sigma} \varphi_1\big(\tilde{a}_1,\sigma\big) \,\tilde{b}_1\cdots \varphi_k\big(\tilde{a}_k,\sigma\big)\,\tilde{b}_k\,\varphi_{k+1}\big(\tilde{a}_{k+1},\sigma\big)\,\rho(d\sigma)\Bigg)h_1,h_2 \Bigg\ra_K.
\end{align*}
By the separable case, the last quantity is independent of the chosen $\ell^{\infty}$-IPD.
Since $h_1,h_2 \in H$ were arbitrary, this completes the proof of this item.

\ref{item.bddlin} First, $k$-linearity of $b \mapsto (I^{\boldsymbol{a}}\varphi)[b]$ is clear by linearity of pointwise Pettis integrals.
Second, we note that Eq. \eqref{eq.MOIbd} gives
\begin{align*}
    \big\|\big(I^{\boldsymbol{a}}\varphi\big)[b]\big\| & = \Bigg\|\int_{\Sigma} F\,d\rho\Bigg\| = \sup\Bigg\{\Bigg|\Bigg\la \Bigg(\int_{\Sigma} F\,d\rho\Bigg)h_1,h_2 \Bigg\ra\Bigg| : \|h_1\|,\|h_2\| \leq 1\Bigg\} \\
    & \leq  \sup\Bigg\{\int_{\Sigma} |\la F(\sigma) h_1,h_2 \ra| \,\rho(d\sigma) : \|h_1\|,\|h_2\| \leq 1\Bigg\} \\
    & \leq \|b_1\|\cdots\|b_k\|\int_{\Sigma}\|\varphi_1(\cdot,\sigma)\|_{\ell^{\infty}(\sigma(a_1))}\cdots\|\varphi_{k+1}(\cdot,\sigma)\|_{\ell^{\infty}(\sigma(a_{k+1}))}\,\rho(d\sigma)
\end{align*}
whenever $(\Sigma,\rho,\varphi_1,\ldots,\varphi_{k+1})$ is a $\ell^{\infty}$-IPD of $\varphi$.
Taking the infimum over $\ell^{\infty}$-IPDs of $\varphi$ therefore gives
\[
\big\|I^{\boldsymbol{a}}\varphi\big\|_{B_k(\cM^k;\cM)} \leq \|\varphi\|_{\ell^{\infty}(\sigma(a_1),\cB_{\sigma(a_1)}\hspace{-0.1mm}) \iotimes \cdots \iotimes \ell^{\infty}(\sigma(a_{k+1}),\cB_{\sigma(a_{k+1})}\hspace{-0.1mm})}.
\]
It remains to prove that $\varphi \mapsto I^{\boldsymbol{a}}\varphi$ is linear.
It is easy to see that $I^{\boldsymbol{a}}(c \,\varphi) = c\,I^{\boldsymbol{a}}\varphi$.
Additivity follows from two facts.
First, for $j \in \{1,2\}$, suppose that $(\Sigma_j,\sH_j,\rho_j)$ is a measure space and $F_j \colon \Sigma_j \to B(H)$ is pointwise Pettis integrable.
If $(\Sigma,\sH,\rho)$ is the disjoint union of $(\Sigma_1,\sH_1,\rho_1)$ and $(\Sigma_2,\sH_2,\rho_2)$, and $F_0 \colon \Sigma \to B(H)$ is defined by $F_0|_{\Sigma_j} = F_j$ for $j \in \{1,2\}$, then $F$ is pointwise Pettis integrable and
\[
\int_{\Sigma} F_0 \,d\rho = \int_{\Sigma_1} F_1\,d\rho_1 + \int_{\Sigma_2} F_2\,d\rho_2.
\]
Second, suppose that $(\Sigma_1,\rho_1,\varphi_1,\ldots,\varphi_{k+1})$ and $(\Sigma_2,\rho_2,\psi_1,\ldots,\psi_{k+1})$ are $\ell^{\infty}$-IPDs of $\varphi$ and $\psi$, respectively.
If $(\Sigma,\sH,\rho)$ is again the disjoint union of $(\Sigma_1,\sH_1,\rho_1)$ and $(\Sigma_2,\sH_2,\rho_2)$, and, for $j \in \{1,\ldots,k+1\}$, we define $\chi_j \colon \sigma(a_j) \times \Sigma \to \C$ by $\chi_j|_{\sigma(a_j) \times \Sigma_1} = \varphi_j$ and $\chi_j|_{\sigma(a_j) \times \Sigma_2} = \psi_j$, then $(\Sigma,\rho,\chi_1,\ldots,\chi_{k+1})$ is a $\ell^{\infty}$-IPD of $\varphi + \psi$.
(Examine the proof of Proposition \ref{prop.babyIPTP}\ref{item.Banachstaralg}.)
We leave it to the reader to combine these two facts to conclude that $I^{\boldsymbol{a}}(\varphi+\psi) = I^{\boldsymbol{a}}\varphi+I^{\boldsymbol{a}}\psi$.
This completes the proof.
\end{proof}

\begin{ex}[Matrix MOIs]\label{ex.fdMOI}
If $S$ and $T$ are sets, then we shall write $T^S$ for the set of all functions $S \to T$.
We claim that if $S_1,\ldots,S_m$ are finite discrete spaces, then $\ell^{\infty}(S_1,\cB_{S_1}) \iotimes \cdots \iotimes \ell^{\infty}(S_m,\cB_{S_m})$ is all of $\C^{S_1 \times \cdots \times S_m}$.
Indeed, if $\varphi \colon S_1 \times \cdots \times S_m \to \C$ is any function, then
\[
\varphi(s_1,\ldots,s_m) = \sum_{(t_1,\ldots,t_m) \in S_1 \times \cdots \times S_m} \varphi(t_1,\ldots,t_m) \, 1_{\{t_1\}}(s_1) \cdots 1_{\{t_m\}}(s_m),
\]
for all $(s_1,\ldots,s_m) \in S_1 \times \cdots \times S_m$, which is easily seen to give a $\ell^{\infty}$-IPD of $\varphi$.
This particular $\ell^{\infty}$-IPD also helps us to compute matrix MOIs.
Specifically, if $n \in \N$, $\boldsymbol{A} = (A_1,\ldots,A_{k+1}) \in \MnC_{\nu}^{k+1}$, and $\varphi \colon \sigma(A_1) \times \cdots \times \sigma(A_{k+1}) \to \C$ is any function, then (as we encourage the reader to verify)
\[
\big(I^{\boldsymbol{A}}\varphi\big)[B] = \sum_{\blambda \in \sigma(A_1) \times \cdots \times \sigma(A_{k+1})} \varphi(\blambda) \, P_{\lambda_1}^{A_1}B_1\cdots P_{\lambda_k}^{A_k}B_kP_{\lambda_{k+1}}^{A_{k+1}},
\]
for all $B = (B_1,\ldots,B_k) \in \MnC^k$.
This agrees with Notation \ref{nota.fdMOI} in the appendix.
\end{ex}

\begin{ex}[Polynomials]\label{ex.polyMOI}
Fix $P(\blambda) = \sum_{|\gamma| \leq d} c_{\gamma}\,\blambda^{\gamma} \in \C[\lambda_1,\ldots,\lambda_{k+1}]$.
By Example \ref{ex.mvarpoly} and the definition of MOIs, if $\boldsymbol{a} = (a_1,\ldots,a_{k+1}) \in \cM_{\nu}^{k+1}$ and $b = (b_1,\ldots,b_k) \in \cM^k$, then
\[
\big(I^{\boldsymbol{a}}P\big)[b] = \sum_{|\gamma| \leq d} c_{\gamma}\,a_1^{\gamma_1}b_1\cdots a_k^{\gamma_k}b_k\,a_{k+1}^{\gamma_{k+1}}.
\]
In particular, by Example \ref{ex.polydivdiff}, if $n \in \N_0$ and $p_n(\lambda) \coloneqq \lambda^n$, then
\[
    \big(I^{\boldsymbol{a}}p_n^{[k]}\big)[b] = \sum_{|\gamma|=n-k}a_1^{\gamma_1}b_1\cdots a_k^{\gamma_k}b_k\,a_{k+1}^{\gamma_{k+1}}.
\]
This formula will be important for our application to differentiating operator functions.
\end{ex}

\begin{ex}\label{ex.WienerMOI}
Fix $f = \int_{\R}e^{i\boldsymbol{\cdot}\xi}\,\mu(d\xi) \in W_k(\R)$.
By the proof of Lemma \ref{lem.WkcCk} and the definition of MOIs,
\[
\big(I^{\boldsymbol{a}}f^{[k]}\big)[b] = \int_{\R \times \Delta_k}e^{it_1\xi a_1}b_1\cdots e^{it_k\xi a_k}b_k\,e^{it_{k+1}\xi a_{k+1}} (i\xi)^k\,(\mu\otimes\rho_j)(d\xi,d \boldsymbol{t}),
\]
for all $\boldsymbol{a} = (a_1,\ldots,a_{k+1}) \in \cM_{\sa}^{k+1}$ and $b = (b_1,\ldots,b_k) \in \cM^k$.
\end{ex}

\begin{rem}[Another Useful Approach]
While MOIs defined by separation of variables suffice for our present purposes, it is worth mentioning that there is another approach to defining MOIs --- due originally to B. S. Pavlov \cite{pavlov} --- that has also been used frequently in the literature, e.g., to differentiate operator functions in the Schatten $p$-norms \cite{coine,lemerdymcdonald,lemerdyskripka}.
(We refer the reader again to \cite{skripka} for more applications.)
When both approaches apply, MOIs constructed using Pavlov's approach agree with MOIs constructed using the separation of variables approach;
please see \cite{coinelemerdysukochev,nikitopoulosMOI} for proofs and additional discussion.
\end{rem}

\subsection{Proof of Theorem \ref{thm.main2}}\label{sec.diffopfunccalc}

We are now finally ready to differentiate operator functions.
In doing so, we shall take the formulas for the derivatives of the $m^{\text{th}}$ power map in a Banach algebra as a starting point.
Though these formulas are well-known, we supply a proof for the convenience of the reader.

\begin{prop}[Derivatives of the Power Map]\label{prop.polyderiv}
If $\cA$ is a Banach algebra over $\F \in \{\R,\C\}$, $m \in \N$, and $F_m(a) \coloneqq a^m$ for $a \in \cA$, then $F_m \in BC_{\loc}^{\infty}(\cA;\cA)$ and
\[
D^kF_m(a)[b_1,\ldots,b_k] = \partial_{b_1}\cdots\partial_{b_k}a^m = \sum_{\pi \in S_k}\sum_{|\gamma| =m-k} a^{\gamma_1}b_{\pi(1)}\cdots a^{\gamma_k}b_{\pi(k)}a^{\gamma_{k+1}}, \numberthis \label{eq.polyderiv}
\]
for all $k \in \N$ and $a,b_1,\ldots,b_k \in \cA$.
\end{prop}
\begin{proof}
Fix $m \in \N$ and $a \in \cA$, and write $\|\cdot\| = \|\cdot\|_{\cA}$.
We prove Eq. \eqref{eq.polyderiv} by induction on $k$.
To this end, we note that if $c \in \cA$, then
\begin{align*}
    F_m(a+c)-F_m(a) & = (a+c)^m-a^m = \sum_{j=0}^{m-1}\big((a+c)^{j+1}a^{m-(j+1)}-(a+c)^ja^{m-j}\big) \\
    & = \sum_{j=0}^{m-1}(a+c)^j(a+c-a)a^{m-1-j} = \sum_{\gamma_1+\gamma_2=m-1}(a+c)^{\gamma_1}c\,a^{\gamma_2} \numberthis\label{eq.pert}
\end{align*}
by telescoping the sum.
Thus
\begin{align*}
    \frac{1}{\|c\|}\Bigg\|F_m(a+c)-F_m(a) - \sum_{\gamma_1+\gamma_2=m-1}a^{\gamma_1}c\,a^{\gamma_2}\Bigg\| & = \frac{1}{\|c\|}\Bigg\|\sum_{\gamma_1+\gamma_2=m-1}\big((a+c)^{\gamma_1}-a^{\gamma_1}\big)c\,a^{\gamma_2}\Bigg\| \\
    & \leq \sum_{\gamma_1+\gamma_2=m-1}\big\|(a+c)^{\gamma_1}-a^{\gamma_1}\big\|\,\|a\|^{\gamma_2} \to 0
\end{align*}
as $c \to 0$.
This proves Eq. \eqref{eq.polyderiv} when $k=1$.

Now, assume Eq. \eqref{eq.polyderiv} holds for $k \in \N$ (and all $m \in \N$ and $a,b_1,\ldots,b_k \in \cA$).
In order to prove Eq. \eqref{eq.polyderiv} with $k$ replaced by $k+1$, we set some notation and make a few preliminary observations.
For $j \in \N$, $\gamma \in \N_0^{j+1}$, $b = (b_1,\ldots,b_j) \in \cA^j$, $\pi \in S_j$, $\boldsymbol{a} = (a_1,\ldots,a_{j+1}) \in \cA^{j+1}$, a set $S$, and $s \in S$, we define
\[
b^{\pi} \coloneqq (b_{\pi(1)},\ldots,b_{\pi(j)}), \; \Pi^{\gamma}(\boldsymbol{a},b) \coloneqq a_1^{\gamma_1}b_1\cdots a_j^{\gamma_j}b_j a_{j+1}^{\gamma_{j+1}}, \text{ and } \, s_{(j)} \coloneqq (s,\ldots,s) \in S^j.
\]
Using this notation with $j=k$, the inductive hypothesis may be rewritten as
\[
D^kF_m(a)[b] = \sum_{\pi \in S_k} \sum_{|\gamma|=m-k} \Pi^{\gamma}\big(a_{(k+1)},b^{\pi}\big), \numberthis\label{eq.ind}
\]
for all $b \in \cA^k$.
Now, if $c \in \cA$, then
\begin{align*}
    \Pi^{\gamma}\big((a+c)_{(k+1)},b\big) - \Pi^{\gamma}\big(a_{(k+1)},b\big) & = \sum_{j=1}^{k+1} \Big(\Pi^{\gamma}\big(\big((a+c)_{(j)},a_{(k+1-j)}\big),b\big) - \Pi^{\gamma}\big(\big((a+c)_{(j-1)},a_{(k+2-j)}\big),b\big)\Big) \\
    & = \sum_{j=1}^{k+1}\Bigg(\prod_{\ell=1}^{j-1} (a+c)^{\gamma_{\ell}}b_{\ell}\Bigg) \big((a+c)^{\gamma_j}-a^{\gamma_j}\big)\prod_{\ell=j}^k b_{\ell}\,a^{\gamma_{\ell+1}} \\
    & = \sum_{j=1}^{k+1}\,\sum_{\delta_{1,j}+\delta_{2,j}=\gamma_j-1}\Bigg(\prod_{\ell=1}^{j-1} (a+c)^{\gamma_{\ell}}b_{\ell}\Bigg) (a+c)^{\delta_{1,j}}c\,a^{\delta_{2,j}}\prod_{\ell=j}^k b_{\ell}\,a^{\gamma_{\ell+1}} \numberthis\label{eq.1pert}
\end{align*}
by telescoping another sum and applying Eq. \eqref{eq.pert}.
(Empty products are defined to be multiplication by $1 \in \F$.)
Next, observe that if $1 \leq j \leq k+1$, $\e = (\e_1,\ldots,\e_j) \in \N_0^j$, and $c_1,\ldots,c_j \in \cA$, then
\begin{align*}
    \prod_{\ell=1}^j (a+c)^{\e_{\ell}}c_{\ell} - \prod_{\ell=1}^j a^{\e_{\ell}}c_{\ell} & = \sum_{q=1}^j\Bigg(\Bigg(\prod_{\ell_1=1}^q(a+c)^{\e_{\ell_1}}c_{\ell_1}\Bigg)\prod_{\ell_2=q+1}^ja^{\e_{\ell_2}}c_{\ell_2} - \Bigg(\prod_{\ell_1=1}^{q-1}(a+c)^{\e_{\ell_1}}c_{\ell_1}\Bigg)\prod_{\ell_2=q}^ja^{\e_{\ell_2}}c_{\ell_2}\Bigg) \\
    & = \sum_{q=1}^j\Bigg(\prod_{\ell=1}^{q-1}(a+c)^{\e_{\ell}}c_{\ell}\Bigg)\big((a+c)^{\e_q}-a^{\e_q}\big)c_q\prod_{\ell=q+1}^j a^{\e_{\ell}}c_{\ell}\numberthis\label{eq.2pert}
\end{align*}
by telescoping yet another sum.
It is then easy to see, by combining Eqs. \eqref{eq.1pert} and \eqref{eq.2pert} (with the multi-index $\e = (\gamma_1,\ldots,\gamma_{j-1},\delta_{1,j})$ and the elements $(c_1,\ldots,c_j) = (b_1,\ldots,b_{j-1},c)$), that
\begin{align*}
    & \frac{1}{\|c\|}\sup_{\|b_1\|,\ldots,\|b_k\| \leq 1} \Bigg\|\Pi^{\gamma}\big((a+c)_{(k+1)},b\big) - \Pi^{\gamma}\big(a_{(k+1)},b\big) \\
    & \hspace{50mm} - \sum_{j=1}^{k+1}\,\sum_{\delta_{1,j}+\delta_{2,j} = \gamma_j-1}\Bigg(\prod_{\ell=1}^{j-1} a^{\gamma_{\ell}}b_{\ell}\Bigg) a^{\delta_{1,j}}c\,a^{\delta_{2,j}}\prod_{\ell=j}^k b_{\ell}\,a^{\gamma_{\ell+1}}\Bigg\| \to 0 \numberthis\label{eq.diff1}
\end{align*}
as $c \to 0$.
After noting that
\begin{align*}
    & \sum_{\tilde{\pi} \in S_{k+1}}\,\sum_{\tilde{\gamma} \in \N_0^{k+2} : |\tilde{\gamma}| = m-k-1}\Pi^{\tilde{\gamma}}\big(a_{(k+2)},(c,b)^{\tilde{\pi}}\big) \\
    & \hspace{20mm} = \sum_{\pi \in S_k}\sum_{|\gamma| = m-k}\Bigg(\sum_{j=1}^{k+1}\,\sum_{\delta_{1,j}+\delta_{2,j} = \gamma_j-1}\Bigg(\prod_{\ell=1}^{j-1} a^{\gamma_{\ell}}b_{\pi(\ell)}\Bigg) a^{\delta_{1,j}}c\,a^{\delta_{2,j}}\prod_{\ell=j}^k b_{\pi(\ell)}\,a^{\gamma_{\ell+1}}\Bigg),
\end{align*}
it follows from Eqs. \eqref{eq.ind} and \eqref{eq.diff1} that
\begin{align*}
    \frac{1}{\|c\|} \Bigg\|D^kF_m(a+c) - D^kF_m(a) - \sum_{\tilde{\pi} \in S_{k+1}}\,\sum_{\tilde{\gamma} \in \N_0^{k+2} : |\tilde{\gamma}| = m-k-1}\Pi^{\tilde{\gamma}}\big(a_{(k+2)},(c,\boldsymbol{\cdot})^{\tilde{\pi}}\big)\Bigg\|_{B_k(\cA^k;\cA)} \to 0
\end{align*}
as $c \to 0$.
By definition of the Fr\'{e}chet derivative, we have shown that if $\tilde{b} \coloneqq (b_0,\ldots,b_k) \in \cA^{k+1}$, then
\[
D^{k+1}F_m(a)\big[\tilde{b}\big] = \sum_{\tilde{\pi} \in S_{k+1}}\,\sum_{\tilde{\gamma} \in \N_0^{k+2} : |\tilde{\gamma}| = m-k-1}\Pi^{\tilde{\gamma}}\big(a_{(k+2)},\tilde{b}^{\tilde{\pi}}\big).
\]
This is the claimed formula, so this completes the proof of Eq. \eqref{eq.polyderiv}.
Since it is clear from Eq. \eqref{eq.polyderiv} that the derivatives of $F_m$ are bounded on bounded sets, we are done.
\end{proof}

We now prove the result that motivates our definition of $NC^k(\R)$:
Theorem \ref{thm.main2}, the infinite-dimensional analog --- proof and all --- of Theorem \ref{thm.diffmatfunc}.

\begin{lem}\label{lem.CstarMOI}
Fix a von Neumann algebra $\cM$, a natural number $k \in \N$, and a function $f \in NC^k(\R)$.
If $\boldsymbol{a} = (a_1,\ldots,a_{k+1}) \in \cM_{\sa}^{k+1}$ and $b = (b_1,\ldots,b_k) \in \cM^k$, then
\[
\big(I^{\boldsymbol{a}}f^{[k]}\big)[b] \in C^*(1,a_1,\ldots,a_{k+1},b_1,\ldots,b_k) \subseteq \cM,
\]
where, for $\cS \subseteq \cM$, $C^*(\cS) \subseteq \cM$ denotes the smallest $C^*$-subalgebra of $\cM$ containing $\cS$.
\end{lem}
\begin{proof}
Fix a sequence $(q_n)_{n \in \N}$ of polynomials such that $q_n \to f$ in $NC^k(\R)$ as $n \to \infty$.
First, it is clear from Example \ref{ex.polyMOI} that
\[
\big(I^{\boldsymbol{a}}q_n^{[k]}\big)[b] \in C^*(1,a_1,\ldots,a_{k+1},b_1,\ldots,b_k),
\]
for all $n \in \N$.
Now, if we define $r \coloneqq \max\{\|a_j\| : 1 \leq j \leq k+1\}$, then
\begin{align*}
    \big\|\big(I^{\boldsymbol{a}}f^{[k]}\big)[b] - \big(I^{\boldsymbol{a}}q_n^{[k]}\big)[b]\big\| & = \big\|\big(I^{\boldsymbol{a}}\big(f^{[k]}-q_n^{[k]}\big)\big)[b]\big\| \\
    & \leq \big\|f^{[k]} - q_n^{[k]}\big\|_{\ell^{\infty}(\sigma(a_1),\cB_{\sigma(a_1)}\hspace{-0.1mm}) \iotimes \cdots \iotimes \ell^{\infty}(\sigma(a_{k+1}),\cB_{\sigma(a_{k+1})}\hspace{-0.1mm})}\|b_1\|\cdots\|b_k\| \\
    & \leq \big\|f^{[k]} - q_n^{[k]}\big\|_{r,k+1}\|b_1\|\cdots\|b_k\| \leq \|f-q_n\|_{\cC^{[k]},r}\|b_1\|\cdots\|b_k\| \to 0 
\end{align*}
as $n \to \infty$ by Theorem \ref{thm.babyMOI}\ref{item.bddlin} and the fact that $q_n \to f$ in $NC^k(\R)$ as $n \to \infty$.
The result follows.
\end{proof}

\begin{proof}[Proof of Theorem \ref{thm.main2}]\hspace{-0.3mm}
We first set some notation.
If $V$ and $W$ are normed vector spaces over $\F \in \{\R,\C\}$ and $T \in B_k(V^k;W)$, then we define $\operatorname{Sym}(T)[v_1,\ldots,v_k] \coloneqq \sum_{\pi \in S_k}T(v_{\pi(1)},\ldots,v_{\pi(k)})$, for all $v_1,\ldots,v_k \in V$.
Using this notation, we may rewrite Eq. \eqref{eq.main2} as
\[
D^kf_{\mathsmaller{\cA}}(a) = \operatorname{Sym}\big(\big(I^{a,\ldots,a}f^{[k]}\big)|_{\cA_{\sa}^k}\big), \numberthis\label{eq.opfunccalcderrewrite}
\]
for all $a \in \cA_{\sa}$.

Now, for $n \in \N$, define $p_n(\lambda) \coloneqq \lambda^n$ as usual.
Then $(p_n)_{\mathsmaller{\cA}} = F_n|_{\cA_{\sa}}$ in the notation of Proposition \ref{prop.polyderiv}.
By Proposition \ref{prop.polyderiv} and Example \ref{ex.polyMOI}, Eq. \eqref{eq.opfunccalcderrewrite} holds when $f=p_n$.
Also, it is obvious that Eq. \eqref{eq.opfunccalcderrewrite} holds when $f=p_0 \equiv 1$.
Therefore, by linearity, Eq. \eqref{eq.opfunccalcderrewrite} holds for all $f \in \C[\lambda]$.

Finally, fix $f \in NC^k(\R)$ and a sequence $(q_n)_{n \in \N}$ of polynomials converging to $f$ in $NC^k(\R)$.
By Lemma \ref{lem.CstarMOI}, if $\boldsymbol{a} \in \cA_{\sa}^{k+1}$ and $b \in \cA^k$, then $(I^{\boldsymbol{a}}f^{[k]})[b] \in \cA$.
We shall take this for granted in our notation.
Now, for any $r >0$, define
\[
\cA_{\sa,r} \coloneqq \{a \in \cA_{\sa} : \|a\| \leq r\}.
\]
If $a \in \cA_{\sa,r}$, then the spectral radius $r(a)$ of $a$ is equal to $\|a\|$, which is at most $r$.
Therefore,
\[
\sup_{a \in \cA_{\sa,r}}\|f(a)-q_n(a)\| = \|f-q_n\|_{\ell^{\infty}([-r,r])} \to 0\pagebreak
\]
as $n \to \infty$.
Also, writing $\|\cdot\|_j \coloneqq \|\cdot\|_{B_j(\cA_{\sa}^j;\cA)}$ for $j \in \N$, Theorem \ref{thm.babyMOI}\ref{item.bddlin} and the previous paragraph give that if $a \in \cA_{\sa,r}$ and $j \in \{1,\ldots,k\}$, then
\begin{align*}
    \Big\|\operatorname{Sym}\big(I^{a,\ldots,a}f^{[j]}\big) - D^j(q_n)_{\mathsmaller{\cA}}(a)\Big\|_j & = \Big\|\operatorname{Sym}\big(I^{a,\ldots,a}\big(f^{[j]}-q_n^{[j]}\big)\big)\Big\|_j \leq j!\Big\|I^{a,\ldots,a}\big(f^{[j]}-q_n^{[j]}\big)\Big\|_j \\
    & \leq j! \big\|f^{[j]}-q_n^{[j]}\big\|_{\ell^{\infty}(\sigma(a),\cB_{\sigma(a)})^{\iotimes(j+1)}} \leq j!\big\|f^{[j]}-q_n^{[j]}\big\|_{r,j+1}.
\end{align*}
In particular,
\[
\max_{1 \leq j \leq k}\sup_{a \in \cA_{\sa,r}}\Big\|\operatorname{Sym}\big(I^{a,\ldots,a}f^{[j]}\big) - D^j(q_n)_{\mathsmaller{\cA}}(a)\Big\|_j \leq k! \|f-q_n\|_{\cC^{[k]},r} \to 0
\]
as $n \to \infty$.
Since $r > 0$ was arbitrary, we conclude from Proposition \ref{prop.Frechcomp} that $f_{\mathsmaller{\cA}} \in BC_{\loc}^k(\cA_{\sa};\cA)$ and that Eq. \eqref{eq.opfunccalcderrewrite} holds for all $a \in \cA_{\sa}$.
\end{proof}

\subsection{Demonstration that \texorpdfstring{$NC^k(\R) \subsetneq C^k(\R)$}{}}\label{sec.NCknonex}

In Section \ref{sec.NCkexs2}, we saw that $NC^k(\R)$ is ``close" to $C^k(\R)$ in the sense that a function only has to be ``slightly better than $C^k$" to belong to $NC^k(\R)$.
The goal of this section is to show that nevertheless $NC^k(\R) \subsetneq C^k(\R)$, for all $k \in \N$.
Specifically, we combine Schatten estimates for Taylor remainders of operator functions (Proposition \ref{prop.remainderestim}) with a construction of D. Potapov et al. from \cite{potapov} (Theorem \ref{thm.Remainder}) to prove the following.

\begin{thm}\label{thm.NCkcounterex}
If $k \in \N$, then $NC^k(\R) \subsetneq C^k(\R)$.
Specifically, we have the following counterexample.
Fix $\eta \in C_c^{\infty}(\R)$ such that $\eta \equiv 1$ on $\big[-\frac{1}{2}-\frac{1}{e},\frac{1}{e}+\frac{1}{2}\big]$ and $\supp \eta \subseteq \big[-\frac{3}{5}-\frac{1}{e},\frac{1}{e}+\frac{3}{5}\big]$, and define
\[
h(x) \coloneqq 1_{(0,1)}(|x|)\frac{\eta(x)\,|x|}{\sqrt{\log|\log|x|-1|}}, \; x \in \R.
\]
If $f_k(x) \coloneqq x^{k-1}h(x)$, for all $x \in \R$, then $f_k \in C^k(\R) \setminus NC^k(\R)$.
\end{thm}

To begin, we set some notation for Taylor remainders.

\begin{defi}[Taylor Remainder]
Let $V$ and $W$ be normed vector spaces and $k \in \N$.
If $F \colon V \to W$ is $(k-1)$-times Fr\'{e}chet differentiable and $p,h \in V$, then we define
\[
R_{k,F,p}(h) \coloneqq F(p+h)-F(p) - \sum_{j=1}^{k-1} \frac{1}{j!}\partial_h^jF(p) = F(p+h) - F(p) - \sum_{j=1}^{k-1}\frac{1}{j!}D^jF(p)[\underbrace{h,\ldots,h}_{j \text{ times}}] \in W.
\]
We call $R_{k,F,p} \colon V \to W$ the \textbf{$\boldsymbol{k^{\text{\textbf{th}}}}$ Taylor Remainder of $\boldsymbol{F}$ at $\boldsymbol{p}$}.
\end{defi}

Recall that if $f \in C^k(\R)$, then $f \in NC^{k-1}(\R)$ (if we take $NC^0(\R) \coloneqq C(\R)$) by Theorem \ref{thm.WkNCk}.
Therefore, if $\cA$ is a unital $C^*$-algebra, then $f_{\mathsmaller{\cA}} \in C^{k-1}(\cA_{\sa};\cA)$ by Theorem \ref{thm.main2}.
In particular, $R_{k,f_{\text{\scalebox{0.7}{$\mathsmaller{\cA}$}}},a}(b) \in \cA$ makes sense whenever $f \in C^k(\R)$ and $a,b \in \cA_{\sa}$.
Now, we state one of the key ingredients of the proof of Theorem \ref{thm.NCkcounterex}.
If $H$ is a Hilbert space and $1 \leq p < \infty$, then we use the notation $(\cS_p(H),\|\cdot\|_{\cS_p})$ for the ideal of Schatten $p$-class operators (Definition 2.2.1 in \cite{nikitopoulosMOI}) on $H$ and $(\cS_{\infty}(H),\|\cdot\|_{\cS_{\infty}})$ for $(B(H),\|\cdot\|_{H \to H})$.

\begin{thm}[Potapov--Skripka--Sukochev--Tomskova \cite{potapov}]\label{thm.Remainder}
Let $k \in \N$ and $f_k \colon \R \to \C$ be as in Theorem \ref{thm.NCkcounterex}.
There exists a separable complex Hilbert space $H$ and bounded linear operators $a,b \in B(H)_{\sa}$ such that $b \in \cS_k(H)$ and $R_{k,(f_k)_{\text{\scalebox{0.7}{$\mathsmaller{B(H)}$}}},a}(b) \not\in \cS_1(H)$.
\end{thm}

This is Theorem 5.1 in \cite{potapov}.
Next, we work toward the Taylor remainder estimates that will help to disqualify $f_k$ from belonging to $NC^k(\R)$.
For more information about applications of MOI theory to Taylor remainders of operator functions, please see Section 5.4 of \cite{skripka}.

\begin{lem}[Schatten Estimates for MOIs]\label{lem.SpMOI}
Suppose $p,p_1,\ldots,p_k \in [1,\infty]$ satisfy $\frac{1}{p}=\frac{1}{p_1}+\cdots+\frac{1}{p_k}$.
If $H$ is a complex Hilbert space, $\boldsymbol{a} = (a_1,\ldots,a_{k+1}) \in B(H)_{\nu}^{k+1}$, and
\[
\varphi \in \ell^{\infty}(\sigma(a_1),\cB_{\sigma(a_1)}) \iotimes \cdots \iotimes \ell^{\infty}(\sigma(a_{k+1}),\cB_{\sigma(a_{k+1})}),
\]
then
\[
\big\|\big(I^{\boldsymbol{a}}\varphi \big)[b]\big\|_{\cS_p} \leq \|\varphi\|_{\ell^{\infty}(\sigma(a_1),\cB_{\sigma(a_1)}\hspace{-0.1mm}) \iotimes \cdots \iotimes \ell^{\infty}(\sigma(a_{k+1}),\cB_{\sigma(a_{k+1})}\hspace{-0.1mm})}\|b_1\|_{\cS_{p_1}}\cdots\|b_k\|_{\cS_{p_k}},
\]
for all $b = (b_1,\ldots,b_k) \in B(H)^k$.
(As usual, $0 \cdot \infty \coloneqq 0$.)
\end{lem}
\pagebreak
\begin{proof}
Fix $b = (b_1,\ldots,b_k) \in B(H)^k$ and a $\ell^{\infty}$-IPD $(\Sigma,\rho,\varphi_1,\ldots,\varphi_{k+1})$ of $\varphi$.
Also, if $\sigma \in \Sigma$, then we define $F(\sigma) \coloneqq \varphi_1(a_1,\sigma)\,b_1\cdots \varphi_k(a_k,\sigma)\,b_k\,\varphi_{k+1}(a_{k+1},\sigma)$.
If $p_1,\ldots,p_k,p \in [1,\infty]$ satisfy $\frac{1}{p_1}+\cdots+\frac{1}{p_k}=\frac{1}{p}$, then
\begin{align*}
    \|F(\sigma)\|_{\cS_p} & \leq \|\varphi_1(a_1,\sigma)\|_{\cS_{\infty}}\|b_1\|_{\cS_{p_1}}\cdots\|\varphi_k(a_k,\sigma)\|_{\cS_{\infty}}\|b_k\|_{\cS_{p_k}}\|\varphi_{k+1}(a_{k+1},\sigma)\|_{\cS_{\infty}} \\
    & \leq \|\varphi_1(\cdot,\sigma)\|_{\ell^{\infty}(\sigma(a_1))}\,\|b_1\|_{\cS_{p_1}}\cdots\|\varphi_k(\cdot,\sigma)\|_{\ell^{\infty}(\sigma(a_k))}\,\|b_k\|_{\cS_{p_k}}\|\varphi_{k+1}(\cdot,\sigma)\|_{\ell^{\infty}(\sigma(a_{k+1}))}
\end{align*}
by H\"{o}lder's Inequality for the Schatten norms.
Therefore, by the Schatten $p$-norm Minkowski Integral Inequality (Theorem 3.4.3 in \cite{nikitopoulosMOI}), we have
\[
\big\|\big(I^{\boldsymbol{a}}\varphi\big)[b]\big\|_{\cS_p}  = \Bigg\|\int_{\Sigma} F\,d\rho\Bigg\|_{\cS_p} \leq \|b_1\|_{\cS_{p_1}}\cdots\|b_k\|_{\cS_{p_k}}\int_{\Sigma}\|\varphi_1(\cdot,\sigma)\|_{\ell^{\infty}(\sigma(a_1))}\cdots\|\varphi_{k+1}(\cdot,\sigma)\|_{\ell^{\infty}(\sigma(a_{k+1}))}\,\rho(d\sigma).
\]
Taking the infimum over $\ell^{\infty}$-IPDs of $\varphi$ then gives the desired estimate.
\end{proof}

\begin{prop}[Taylor Remainder Estimates]\label{prop.remainderestim}
Fix a complex Hilbert space $H$, a natural number $k \in \N$, and a function $f \in NC^k(\R)$.
If $a,b \in B(H)_{\sa}$ and we define $r \coloneqq \max\{\|a+tb\| : 0 \leq t \leq 1\}$, then
\[
\big\|R_{k,f_{\text{\scalebox{0.7}{$\mathsmaller{B(H)}$}}},a}(b)\big\|_{\cS_p} \leq \big\|f^{[k]}\big\|_{r,k+1} \|b\|_{\cS_{kp}}^k, \numberthis\label{eq.remestim}
\]
for all $p \in [1,\infty)$.
In particular, if in addition $b \in \cS_{kp}(H)$, then $R_{k,f_{\text{\scalebox{0.7}{$\mathsmaller{B(H)}$}}},a}(b) \in \cS_p(H)$.
\end{prop}
\begin{proof}
We begin by recalling one form of Taylor's Theorem (e.g., Theorem 107 of Section 1.4 in \cite{hajek}):
if $V$ is a normed vector space, $W$ is a Banach space, and $F \in C^k(V;W)$, then
\[
R_{k,F,p}(h) = \frac{1}{(k-1)!}\int_0^1(1-t)^{k-1}\partial_h^kF(p+th)\, dt,
\]
for all $p,h \in V$, where the integral above is a vector-valued Riemann integral.
In particular, if $f \in NC^k(\R)$, then, by Theorem \ref{thm.main2} (with $\cA = \cM = B(H)$), we have
\[
R_{k,f_{\text{\scalebox{0.7}{$\mathsmaller{B(H)}$}}},a}(b) = k\int_0^1(1-t)^{k-1}\big(I^{a+tb,\ldots,a+tb}f^{[k]}\big)[b,\ldots,b]\, dt. \numberthis\label{eq.integrem}
\]
In this case, the integral above is also a pointwise Pettis integral, as we urge the reader to check.
Now, if $t \in [0,1]$, then $\sigma(a+tb) \subseteq [-r,r]$.
Therefore, if $p \in [1,\infty)$, then Lemma \ref{lem.SpMOI} gives
\[
\big\|\big(I^{a+tb,\ldots,a+tb}f^{[k]}\big)[b,\ldots,b]\big\|_{\cS_p} \leq \big\|f^{[k]}\big\|_{\ell^{\infty}(\sigma(a+tb),\cB_{\sigma(a+tb)}\hspace{-0.1mm})^{\iotimes (k+1)}}\|b\|_{\cS_{kp}}^k \leq \big\|f^{[k]}\big\|_{r,k+1}\|b\|_{\cS_{kp}}^k.
\]
Thus
\begin{align*}
    \big\|R_{k,f_{\text{\scalebox{0.7}{$\mathsmaller{B(H)}$}}},a}(b)\big\|_{\cS_p} & \leq k \big\|f^{[k]}\big\|_{r,k+1}\|b\|_{\cS_{kp}}^k\int_0^1(1-t)^{k-1}\,dt = \big\|f^{[k]}\big\|_{r,k+1}\|b\|_{\cS_{kp}}^k
\end{align*}
by Eq. \eqref{eq.integrem} and the Schatten $p$-norm Minkowski Integral Inequality (Theorem 3.4.3 in \cite{nikitopoulosMOI}).
\end{proof}

We are now ready.

\begin{proof}[Proof of Theorem \ref{thm.NCkcounterex}]
Fix $k \in \N$.
It is shown in Appendix A of \cite{potapov} that $f_k \in C^k(\R)$.
Now, let $H$ be a complex Hilbert space and $a,b \in B(H)_{\sa}$. If $f \in NC^k(\R)$ and $b \in \cS_k(H)$, then $R_{k,f_{\text{\scalebox{0.7}{$\mathsmaller{B(H)}$}}},a}(b) \in \cS_1(H)$ by Proposition \ref{prop.remainderestim}.
Therefore, $f_k \not\in NC^k(\R)$ by Theorem \ref{thm.Remainder}.
\end{proof}

In view of Theorem \ref{thm.main2}, there is another possible approach to proving $f_k \not\in NC^k(\R)$.

\begin{conj}\label{conj.fk}
If $k \in \N$ and $f_k$ is as in Theorem \ref{thm.NCkcounterex}, then $(f_k)_{\mathsmaller{B(\ell^{\text{\scalebox{0.7}{$2$}}}(\N))}} \colon B(\ell^2(\N))_{\sa} \to B(\ell^2(\N))$ is not $k$-times Fr\'{e}chet differentiable.
\end{conj}

If this conjecture is correct, then we would immediately conclude $f_k \not\in NC^k(\R)$ from Theorem \ref{thm.main2}.
From private correspondence with E. McDonald and F. A. Sukochev, it seems possible that ideas from \cite{potapov} could be adapted to prove Conjecture \ref{conj.fk}, but --- to the author's knowledge --- this has never been carried out.
Interestingly, for $k \geq 2$, it even seems to be the case that the literature lacks explicit examples of functions $f \in C^k(\R)$ such that $f_{\mathsmaller{B(\ell^{\text{\scalebox{0.7}{$2$}}}(\N))}} \colon B(\ell^2(\N))_{\sa} \to B(\ell^2(\N))$ has been confirmed not to be $k$-times Fr\'{e}chet differentiable, though it is widely accepted that such functions should exist.
(As we mentioned in the introduction, when $k=1$, any $f \in C^1(\R) \setminus \dot{B}_1^{1,1}(\R)_{\loc}$ would do, by results of Peller \cite{peller0}.)
\pagebreak

We end this section by briefly sketching a slight upgrade to Theorem \ref{thm.NCkcounterex} for the interested reader:
that actually $f_k \in C^k(\R) \setminus \cC^{[k]}(\R)$.
To prove this, we shall use perturbation formulas from \cite{nikitopoulosOpDer} (Equations (12) and (13) therein), which --- in keeping with the rest of our development --- we have avoided thus far.

\begin{prop}[Taylor Remainder Formula]\label{prop.Taylorexp}
Let $\cA$ be a unital $C^*$-algebra, $\cM$ be a von Neumann algebra containing $\cA$ as a unital $C^*$-subalgebra (e.g., $\cM = \cA^{**}$), and $k \in \N$.
If $f \in \cC^{[k]}(\R)$, then
\[
R_{k,f_{\text{\scalebox{0.7}{$\mathsmaller{\cA}$}}},a}(b) = \underbrace{\int_{\sigma(a)}\cdots\int_{\sigma(a)}}_{k \; \mathrm{ times}}\int_{\sigma(a+b)}f^{[k]}(\blambda) \,P^{a+b}(d\lambda_1)\,b\,P^a(d\lambda_2)\cdots b \,P^a(d\lambda_{k+1}), \numberthis\label{eq.remainderform}
\]
for all $a,b \in \cA_{\sa}$, where the right hand side of Eq. \eqref{eq.remainderform} is a MOI in $\cM$.
\end{prop}
\begin{proof}[Sketch of proof]
First, by a smooth cutoff argument, it suffices to assume $f \in \cC^{[k]}(\R)$ is compactly supported, in which case $f^{[j]} \in \ell^{\infty}(\R,\cB_{\R})^{\iotimes(j+1)}$ whenever $0 \leq j \leq k$.
(This will ensure that we can apply the perturbation formulas from \cite{nikitopoulosOpDer}.)
Under this assumption, we prove Eq. \eqref{eq.remainderform} by induction on $k$.

To begin, we have
\[
R_{1,f_{\text{\scalebox{0.7}{$\mathsmaller{\cA}$}}},a}(b) = f(a+b)-f(a) = \int_{\sigma(a)}\int_{\sigma(a+b)}f^{[1]}(\lambda_1,\lambda_2)\,P^{a+b}(d\lambda_1)\,b\,P^a(d\lambda_2)
\]
by Equation (12) in \cite{nikitopoulosOpDer}.
Now, if $k \in \N$, and we assume that Eq. \eqref{eq.remainderform} holds, then
\begin{align*}
    R_{k+1,f_{\text{\scalebox{0.7}{$\mathsmaller{\cA}$}}},a}(b) & = R_{k,f_{\text{\scalebox{0.7}{$\mathsmaller{\cA}$}}},a}(b) - \frac{1}{k!}\partial_b^kf_{\mathsmaller{\cA}}(a) = \big(I^{a+b,a,\ldots,a}f^{[k]}\big)[b,\ldots,b] - \big(I^{a,\ldots,a}f^{[k]}\big)[b,\ldots,b] \\
    & = \int_{\sigma(a)}\cdots\int_{\sigma(a)}\int_{\sigma(a+b)}f^{[k+1]}(\lambda_1,\ldots,\lambda_{k+2}) \,P^{a+b}(d\lambda_1)\,b\,P^a(d\lambda_2)\cdots b \,P^a(d\lambda_{k+2})
\end{align*}
by the inductive hypothesis, Theorem \ref{thm.main2}, and Equation (13) in \cite{nikitopoulosOpDer}.
This completes the proof.
\end{proof}

\begin{cor}[Upgraded Taylor Remainder Estimates]\label{cor.remestim}
Fix a complex Hilbert space $H$, a natural number $k \in \N$, and a function $f \in \cC^{[k]}(\R)$.
If $a,b \in B(H)_{\sa}$, then
\[
\big\|R_{k,f_{\text{\scalebox{0.7}{$\mathsmaller{B(H)}$}}},a}(b)\big\|_{\cS_p} \leq \big\|f^{[k]}\big\|_{\ell^{\infty}(\sigma(a+b),\cB_{\sigma(a+b)}\hspace{-0.1mm}) \iotimes \ell^{\infty}(\sigma(a),\cB_{\sigma(a)}\hspace{-0.1mm})^{\iotimes k}}  \|b\|_{\cS_{kp}}^k,
\]
for all $p \in [1,\infty)$.
\end{cor}
\begin{proof}
Combine Proposition \ref{prop.Taylorexp} (with $\cA = \cM = B(H)$) and Lemma \ref{lem.SpMOI}.
\end{proof}

By repeating the argument from the proof of Theorem \ref{thm.NCkcounterex}, we see that Theorem \ref{thm.Remainder} and Corollary \ref{cor.remestim} together imply that $f_k \not\in \cC^{[k]}(\R)$.

\appendix
\section{Loose ends}\label{app.looseends}

\subsection{Proof of Theorem \ref{thm.diffmatfunc}}\label{sec.dermatfunc}

In this section, we present a version of Hiai's approach from \cite{hiai} to differentiating matrix functions.
Both the derivative formula and its proof will serve as motivation for the ``infinite-dimensional case" proven in Section \ref{sec.diffopfunccalc}.
We shall freely use information and notation from Sections \ref{sec.divdiff} and \ref{sec.frechder}.
In addition, we shall use the formula for the derivatives of the power map in a Banach algebra (Proposition \ref{prop.polyderiv}).

To begin, we introduce notation for finite-dimensional multiple operator integrals (MOIs).
For the rest of this section, fix $n,k \in \N$ and $\bA = (A_1,\ldots,A_{k+1}) \in \MnC_{\nu}^{k+1}$.

\begin{nota}\label{nota.fdMOI}
If $\varphi \colon \sigma(A_1) \times \cdots \times \sigma(A_{k+1}) \to \C$ is any function, then we write
\[
\big(I^{\bA}\varphi\big)[B] \coloneqq \sum_{\blambda \in \sigma(A_1) \times \cdots \times \sigma(A_{k+1})} \varphi(\blambda) \, P_{\lambda_1}^{A_1}B_1\cdots P_{\lambda_k}^{A_k}B_kP_{\lambda_{k+1}}^{A_{k+1}} \in \MnC,
\]
for all $B = (B_1,\ldots,B_k) \in \MnC^k$.
\end{nota}

\begin{ex}\label{ex.polyfdMOI}
For each $j \in \{1,\ldots,k+1\}$, fix a function $\varphi_j \colon \sigma(A_j) \to \C$, and define
\[
\varphi(\blambda) \coloneqq \varphi_1(\lambda_1)\cdots\varphi_{k+1}(\lambda_{k+1})
\]
for $\blambda \coloneqq (\lambda_1,\ldots,\lambda_{k+1}) \in \sigma(A_1) \times \cdots \times \sigma(A_{k+1})$.
If $B = (B_1,\ldots,B_k) \in \MnC^k$, then
\begin{align*}
    \big(I^{\boldsymbol{A}}\varphi\big)[B] & = \sum_{\blambda \in \sigma(A_1) \times \cdots \times \sigma(A_{k+1})} \varphi_1(\lambda_1)\cdots\varphi_{k+1}(\lambda_{k+1})\,P_{\lambda_1}^{A_1}B_1\cdots P_{\lambda_k}^{A_k}B_kP_{\lambda_{k+1}}^{A_{k+1}} \\
    & = \Bigg(\sum_{\lambda_1 \in \sigma(A_1)}\varphi_1(\lambda_1)\,P_{\lambda_1}^{A_1}\Bigg)B_1\cdots \Bigg(\sum_{\lambda_k \in \sigma(A_k)}\varphi_k(\lambda_k)\,P_{\lambda_k}^{A_k}\Bigg)B_k\sum_{\lambda_{k+1} \in \sigma(A_{k+1})}\varphi_{k+1}(\lambda_{k+1})\,P_{\lambda_{k+1}}^{A_{k+1}} \\
    & = \varphi_1(A_1)B_1\cdots \varphi_k(A_k)B_k\varphi_{k+1}(A_{k+1})
\end{align*}
by Eq. \eqref{eq.fdfunccalc}.
In particular, if $m \in \N_0$ and $p_m(\lambda) \coloneqq \lambda^m$, then
\[
\big(I^{\boldsymbol{A}}p_m^{[k]}\big)[B] = \sum_{|\gamma| = m-k} A_1^{\gamma_1}B_1\cdots A_k^{\gamma_k}B_kA_{k+1}^{\gamma_{k+1}} \numberthis\label{eq.fdpoldivdiffMOI}
\]
by Example \ref{ex.polydivdiff}.
\end{ex}

To prove Theorem \ref{thm.diffmatfunc}, we need just two ingredients:
an operator norm estimate on $I^{\boldsymbol{A}}\varphi$ and density of polynomials in the space $C^k(\R)$ with the $C^k$ topology (first paragraph of Section \ref{sec.NCkdefandprop}).

\begin{lem}\label{lem.fdMOIestim}
Remain in the setting of Notation \ref{nota.fdMOI}.
Then
\[
\big\|I^{\bA}\varphi\big\|_{B_k(\MnC^k;\MnC)} \leq n^k \max\{|\varphi(\blambda)| : \blambda \in \sigma(A_1) \times \cdots \times \sigma(A_{k+1})\}, \numberthis\label{eq.crudefdMOIestim}
\]
where $\MnC$ is given the operator norm $\|\cdot\| = \|\cdot\|_{\C^n \to \C^n}$.
\end{lem}
\begin{proof}
If $B = (B_1,\ldots,B_k) \in \MnC^k$, then Eq. \eqref{eq.fdfunccalc} gives
\begin{align*}
    \big\|I^{\bA}\varphi(B)\big\| & = \Bigg\|\sum_{\lambda_2 \in \sigma(A_2),\ldots,\lambda_{k+1} \in \sigma(A_{k+1})} \varphi(A_1,\lambda_2,\ldots,\lambda_{k+1}) \, B_1P_{\lambda_2}^{A_2}\cdots B_kP_{\lambda_{k+1}}^{A_{k+1}}\Bigg\| \\
    & \leq \sum_{\lambda_2 \in \sigma(A_2),\ldots,\lambda_{k+1} \in \sigma(A_{k+1})}\|\varphi(A_1,\lambda_2,\ldots,\lambda_{k+1})\| \, \|B_1\|\,\big\|P_{\lambda_2}^{A_2}\big\| \cdots \|B_k\| \, \big\|P_{\lambda_{k+1}}^{A_{k+1}}\big\|\\
    & = \sum_{\lambda_2 \in \sigma(A_2),\ldots,\lambda_{k+1} \in \sigma(A_{k+1})}\max_{\lambda_1 \in \sigma(A_1)}|\varphi(\lambda_1,\ldots,\lambda_{k+1})| \, \|B_1\|\cdots \|B_k\| \\
    & \leq n^k\max\{|\varphi(\blambda)| : \blambda \in \sigma(A_1) \times \cdots \times \sigma(A_{k+1})\}\,\|B_1\|\cdots\|B_k\|
\end{align*}
because $A_j$ has at most $n$ distinct eigenvalues, for all $j \in \{1,\ldots,k+1\}$.
\end{proof}
\begin{rem}
It turns out (Proposition 4.1.3 in \cite{skripka}) that
\[
\big\|I^{\bA}\varphi\big\|_{B_k((\MnC,\|\cdot\|_{\HS})^k;(\MnC,\|\cdot\|_{\HS}))} = \max\{|\varphi(\blambda)| : \blambda \in \sigma(A_1) \times \cdots \times \sigma(A_{k+1})\},
\]
where $\|\cdot\|_{\HS}$ is the Hilbert--Schmidt norm.
Due to the inequality $\|\cdot\| \leq \|\cdot\|_{\HS} \leq \sqrt{n} \|\cdot\|$, we may therefore replace the $n^k$ in Eq. \eqref{eq.crudefdMOIestim} with $n^{\frac{k}{2}}$.
Note that even this sharper estimate depends on the dimension $n$ in an unbounded way, which suggests difficulties with the infinite-dimensional case.
\end{rem}

\begin{lem}\label{lem.polydense}
If $k \in \N$, then $\C[\lambda]$ is dense in $C^k(\R)$ with the $C^k$ topology.
\end{lem}
\begin{proof}
We first prove that if $r > 0$ and $f \in C^k(\R)$, then there exists a sequence $(q_n)_{n \in \N}$ of polynomials such that, for all $j \in \{0,\ldots,k\}$, $q_n^{(j)} \to f^{(j)}$ uniformly on $[-r,r]$ as $n \to \infty$.
To this end, use the classical Weierstrass Approximation Theorem to find a sequence $(q_{0,n})_{n \in \N}$ of polynomials such that $q_{0,n} \to f^{(k)}$ uniformly on $[-r,r]$ as $n \to \infty$.
Now, for $\ell \in \{1,\ldots,k\}$ and $n \in \N$, recursively define
\[
q_{\ell,n}(\lambda) \coloneqq f^{(k-\ell)}(0) + \int_0^{\lambda} q_{\ell-1,n}(t)\,dt, \; \lambda \in \R.
\]
Note that $q_{\ell,n} \in \C[\lambda]$.
By an induction argument using the Dominated Convergence Theorem and the Fundamental of Calculus, the sequence $(q_n)_{n \in \N} \coloneqq (q_{k,n})_{n \in \N}$ accomplishes the stated goal.
\pagebreak

Next, fix $f \in C^k(\R)$.
By the previous paragraph, if $N \in \N$, then there exists $q_N \in \C[\lambda]$ such that
\[
\max_{0 \leq j \leq k}\big\|(f-q_N)^{(j)}\big\|_{\ell^{\infty}([-N,N])} < \frac{1}{N}.
\]
Then the sequence $(q_N)_{N \in \N}$ of polynomials converges to $f$ in the $C^k$ topology.
This completes the proof.
\end{proof}

We are now ready.

\begin{proof}[Proof of Theorem \ref{thm.diffmatfunc}]
First, if $V$ and $W$ are normed vector spaces over $\F \in \{\R,\C\}$ and $T \in B_k(V^k;W)$, then we define $\operatorname{Sym}(T)[v_1,\ldots,v_k] \coloneqq \sum_{\pi \in S_k}T(v_{\pi(1)},\ldots,v_{\pi(k)})$, for all $v_1,\ldots,v_k \in V$.
Using this notation, we may rewrite Eq. \eqref{eq.matfunccalcderform} as
\[
D^kf_{\mathsmaller{\MnC}}(A) = \operatorname{Sym}\big(\big(I^{A,\ldots,A}f^{[k]}\big)|_{\MnC_{\sa}^k}\big), \numberthis\label{eq.matfunccalcderrewrite}
\]
for all $A \in \MnC_{\sa}$.

Now, for $m \in \N$, define $p_m(\lambda) \coloneqq \lambda^m$. Then $(p_m)_{\mathsmaller{\MnC}} = F_m|_{\MnC_{\sa}}$ in the notation of Proposition \ref{prop.polyderiv}.
By Proposition \ref{prop.polyderiv} and Example \ref{ex.polyfdMOI}, Eq. \eqref{eq.matfunccalcderrewrite} holds when $f=p_m$.
Also, it is obvious that Eq. \eqref{eq.matfunccalcderrewrite} holds when $f=p_0 \equiv 1$.
Therefore, by linearity, Eq. \eqref{eq.matfunccalcderrewrite} holds for all $f \in \C[\lambda]$.

Finally, fix $f \in C^k(\R)$ and a sequence $(q_N)_{N \in \N}$ of polynomials converging to $f$ in $C^k(\R)$.
Such a sequence exists by Lemma \ref{lem.polydense}.
Now, for any $r >0$, define
\[
\MnC_{\sa,r} \coloneqq \{A \in \MnC_{\sa} : \|A\| \leq r\}.
\]
If $A \in \MnC_{\sa,r}$, then the spectral radius $r(A)$ of $A$ is equal to $\|A\|$, which is at most $r$.
Therefore,
\[
\sup_{A \in \MnC_{\sa,r}}\|f(A)-q_N(A)\| = \|f-q_N\|_{\ell^{\infty}([-r,r])} \to 0
\]
as $N \to \infty$.
Also, writing $\|\cdot\|_j \coloneqq \|\cdot\|_{B_j(\MnC_{\sa}^j;\MnC)}$ for $j \in \N$, Lemma \ref{lem.fdMOIestim} and the previous paragraph give that if $A \in \MnC_{\sa,r}$ and $j \in \{1,\ldots,k\}$, then
\begin{align*}
    \Big\|\operatorname{Sym}\big(I^{A,\ldots,A}f^{[j]}\big) - D^j(q_N)_{\mathsmaller{\MnC}}(A)\Big\|_j & = \Big\|\operatorname{Sym}\big(I^{A,\ldots,A}\big(f^{[j]}-q_N^{[j]}\big)\big)\Big\|_j \leq j!\Big\|I^{A,\ldots,A}\big(f^{[j]}-q_N^{[j]}\big)\Big\|_j \\
    & \leq j! \,n^j\big\|f^{[j]}-q_N^{[j]}\big\|_{\ell^{\infty}(\sigma(A)^{j+1})} \leq j!\,n^j\big\|f^{[j]}-q_N^{[j]}\big\|_{\ell^{\infty}([-r,r]^{j+1})}.
\end{align*}
In particular, by Corollary \ref{cor.easyboundfk},
\[
\max_{1 \leq j \leq k}\sup_{A \in \MnC_{\sa,r}}\Big\|\operatorname{Sym}\big(I^{A,\ldots,A}f^{[j]}\big) - D^j(q_N)_{\mathsmaller{\MnC}}(A)\Big\|_j \leq n^k\max_{1 \leq j \leq k} \big\|f^{(j)}-q_N^{(j)}\big\|_{\ell^{\infty}([-r,r])} \to 0
\]
as $N \to \infty$.
Since $r > 0$ was arbitrary, we conclude from Proposition \ref{prop.Frechcomp} that $f_{\mathsmaller{\MnC}} \in C^k(\MnC_{\sa};\MnC)$ and that Eq. \eqref{eq.matfunccalcderrewrite} holds for all $A \in \MnC_{\sa}$.
This completes the proof.
\end{proof}

The reason this proof works is that the finite-dimensional MOI $I^{\boldsymbol{A}}f^{[k]}$ satisfies a (dimension-dependent) operator norm estimate involving the uniform norm of $f^{[k]}$.
In the infinite-dimensional case, the uniform norm is too weak for this operator norm estimate.
However, there is a stronger norm, the $\ell^{\infty}$-\textit{integral projective tensor norm} (Section \ref{sec.IPTP}), that works.
This motivates our definition of $NC^k(\R)$ in Section \ref{sec.NCkdefandprop}, since it gives us an infinite-dimensional analog (Theorem \ref{thm.main2}) of Theorem \ref{thm.diffmatfunc} and its proof.

\subsection{Proof of Eq. \texorpdfstring{\eqref{eq.I2assymp}}{}}\label{sec.IBP}

Fix $\psi,\chi \in C_c^{\infty}(\R)$ such that $\psi \equiv 1$ on $[-1,1]$, $\supp \psi \subseteq [-2,2]$, $\chi \equiv 1$ on $\big[\frac{3}{4},\frac{3}{2}\big]$, and $\supp \chi \subseteq \big[\frac{1}{2},2\big]$.
Define
\[
\phi(y) \coloneqq y+y^{-1} \; \text{ and } \; g_{\zeta}(y) \coloneqq (1-\chi(y))\,\psi\big(\zeta^{-1}y\big) \, \text{ for } y,\zeta > 0.
\]
We aim to show that
\[
I_2(\zeta) \coloneqq \int_0^{\infty} y^{\frac{1}{2}}\,g_{\zeta}(y)\,e^{-i\zeta \phi(y)}\,dy = O\big(\zeta^{-1}\big) \, \text{ as } \zeta \to \infty.
\]
To do so, we shall need to integrate by parts three times.
We record a few derivatives for this purpose.
First, $\phi'(y) = 1-y^{-2} = \frac{y^2-1}{y^2}$ and $\phi''(y) = 2y^{-3}$.
Second,
\begin{align*}
    \frac{d}{dy}\big(y^{\frac{1}{2}}g_{\zeta}(y)\big) & = \frac{1}{2}y^{-\frac{1}{2}}g_{\zeta}(y) + y^{\frac{1}{2}}g_{\zeta}'(y), \\
    \frac{d^2}{dy^2}\big(y^{\frac{1}{2}}g_{\zeta}(y)\big) & = -\frac{1}{4}y^{-\frac{3}{2}}g_{\zeta}(y) + y^{-\frac{1}{2}}g_{\zeta}'(y) + y^{\frac{1}{2}}g_{\zeta}''(y), \; \text{ and} \\
    \frac{d^3}{dy^3}\big(y^{\frac{1}{2}}g_{\zeta}(y)\big) & = \frac{3}{8}y^{-\frac{5}{2}}g_{\zeta}(y) -\frac{3}{4}y^{-\frac{3}{2}}g_{\zeta}'(y) + \frac{3}{2}y^{-\frac{1}{2}}g_{\zeta}''(y) + y^{\frac{1}{2}}\,g_{\zeta}'''(y).
\end{align*}
Recall now that $\phi'(y) \neq 0$ for $y \in \supp g_{\zeta}$ (since $g_{\zeta} \equiv 0$ near $1$), and note that $\sup_{\zeta \geq 1}\big\|g_{\zeta}^{(k)}\big\|_{\ell^{\infty}(\R)} < \infty$, for all $k \in \N_0$.
Therefore, as $\zeta \to \infty$, we have
\begin{align*}
    I_2(\zeta) & = \int_0^{\infty} y^{\frac{1}{2}}\,g_{\zeta}(y)\,e^{-i\zeta \phi(y)}\,dy = \frac{1}{i\zeta}\int_0^{\infty}\frac{d}{dy}\Bigg(\frac{y^{\frac{1}{2}}\,g_{\zeta}(y)}{\phi'(y)}\Bigg)e^{-i\zeta \phi(y)}\,dy \\
    & = \frac{1}{i\zeta}\int_0^{\infty}\Bigg(\frac{-\phi''(y)}{\phi'(y)^2}y^{\frac{1}{2}}\,g_{\zeta}(y) + \frac{1}{\phi'(y)}\frac{d}{dy}\big(y^{\frac{1}{2}}\,g_{\zeta}(y)\big)\Bigg)\,e^{-i\zeta \phi(y)}\,dy \\
    & = -\frac{2}{i\zeta}\int_0^{\infty}\frac{g_{\zeta}(y)\,y^{\frac{3}{2}}}{(y^2-1)^2}\,e^{-i\zeta \phi(y)}\,dy + \frac{1}{i\zeta}\int_0^{\infty}\frac{1}{\phi'(y)}\frac{d}{dy}\big(y^{\frac{1}{2}}\,g_{\zeta}(y)\big)\,e^{-i\zeta \phi(y)}\,dy \\
    & = O(\zeta^{-1}) - \frac{1}{\zeta^2}\int_0^{\infty}\frac{d}{dy}\Bigg(\frac{1}{\phi'(y)^2}\frac{d}{dy}\big(y^{\frac{1}{2}}\,g_{\zeta}(y)\big)\Bigg)\,e^{-i\zeta \phi(y)}\,dy \\
    & = O(\zeta^{-1}) - \frac{1}{\zeta^2}\int_0^{\infty}\Bigg(\frac{-2\phi''(y)}{\phi'(y)^3}\frac{d}{dy}\big(y^{\frac{1}{2}}\,g_{\zeta}(y)\big)+\frac{1}{\phi'(y)^2}\frac{d^2}{dy^2}\big(y^{\frac{1}{2}}\,g_{\zeta}(y)\big)\Bigg)\,e^{-i\zeta \phi(y)}\,dy \\
    & = O(\zeta^{-1}) + \frac{2}{\zeta^2}\int_0^{\infty}\frac{y^{\frac{5}{2}}g_{\zeta}(y) +2y^{\frac{7}{2}}g_{\zeta}'(y)}{(y^2-1)^3}\,e^{-i\zeta \phi(y)}\,dy - \frac{1}{\zeta^2}\int_0^{\infty}\frac{1}{\phi'(y)^2}\frac{d^2}{dy^2}\big(y^{\frac{1}{2}}\,g_{\zeta}(y)\big)\,e^{-i\zeta \phi(y)}\,dy \\
    & = O(\zeta^{-1})+O(\zeta^{-2})-\frac{1}{i\zeta^3}\int_0^{\infty}\frac{d}{dy}\Bigg(\frac{1}{\phi'(y)^3}\frac{d^2}{dy^2}\big(y^{\frac{1}{2}}\,g_{\zeta}(y)\big)\Bigg)\,e^{-i\zeta \phi(y)}\,dy \\
    & = O(\zeta^{-1})-\frac{1}{i\zeta^3}\int_0^{\infty}\Bigg(\frac{-3\phi''(y)}{\phi'(y)^4}\frac{d^2}{dy^2}\big(y^{\frac{1}{2}}\,g_{\zeta}(y)\big)+\frac{1}{\phi'(y)^3}\frac{d^3}{dy^3}\big(y^{\frac{1}{2}}\,g_{\zeta}(y)\big)\Bigg)\,e^{-i\zeta \phi(y)}\,dy \\
    & = O(\zeta^{-1})-\frac{1}{i\zeta^3}\int_0^{\infty}\Bigg(\frac{\frac{3}{2}y^{\frac{7}{2}}g_{\zeta}(y)-6y^{\frac{9}{2}}g_{\zeta}'(y)-6y^{\frac{11}{2}}g_{\zeta}''(y)}{(y^2-1)^4}+\frac{1}{\phi'(y)^3}\frac{d^3}{dy^3}\big(y^{\frac{1}{2}}\,g_{\zeta}(y)\big)\Bigg)\,e^{-i\zeta \phi(y)}\,dy \\
    & = O(\zeta^{-1})+O(\zeta^{-3})-\frac{1}{i\zeta^3}\int_0^{\infty}\frac{1}{\phi'(y)^3}\frac{d^3}{dy^3}\big(y^{\frac{1}{2}}\,g_{\zeta}(y)\big)\,e^{-i\zeta \phi(y)}\,dy.
\end{align*}
(We leave it to the reader to confirm that there are no boundary terms at zero.)
But
\[
\frac{1}{\phi'(y)^3}\frac{d^3}{dy^3}\big(y^{\frac{1}{2}}\,g_{\zeta}(y)\big) = \frac{y^6}{(y^2-1)^3}\Bigg(\frac{3}{8}y^{-\frac{5}{2}}g_{\zeta}(y) -\frac{3}{4}y^{-\frac{3}{2}}g_{\zeta}'(y) + \frac{3}{2}y^{-\frac{1}{2}}g_{\zeta}''(y) + y^{\frac{1}{2}}\,g_{\zeta}'''(y)\Bigg)
\]
and $g_{\zeta}(y) = 0$ if $y \geq 2\zeta$.
It follows --- because of the dominant $y^{\frac{1}{2}}$ term --- that
\[
\int_0^{\infty}\frac{1}{\phi'(y)^3}\frac{d^3}{dy^3}\big(y^{\frac{1}{2}}\,g_{\zeta}(y)\big)\,e^{-i\zeta \phi(y)}\,dy = O(\zeta^{\frac{3}{2}}) \, \text{ as } \zeta \to \infty.
\]
Thus
\[
I_2(\zeta) = O(\zeta^{-1}) + \zeta^{-3}O(\zeta^{\frac{3}{2}}) = O(\zeta^{-1}) \text{ as } \zeta \to \infty,
\]
as desired. \qed

\begin{ack}
\phantomsection
\addcontentsline{toc}{section}{Acknowledgements}
I acknowledge support from NSF grant DGE 2038238 and partial support from NSF grants DMS 1253402 and DMS 1800733.
I also have many people to thank for help with this paper.
First, I am grateful to Bruce Driver, Adrian Ioana, and Todd Kemp for many helpful conversations and/or guidance about writing.
Second, I would like to thank David Jekel for bringing to my attention his parallel work on $C_{\operatorname{nc}}^k(\R)$ in Section 18 of his dissertation \cite{jekel}.
Third, I would like to thank Jacob Sterbenz for suggesting $\kappa$ (from Section \ref{sec.Wkloc}) as a compactly supported H\"{o}lder continuous function with nonintegrable Fourier transform.
Finally, special thanks go to Edward McDonald.
The first version of this paper worked only with $\cA = B(H)$ but included a full development of the ``separation of variables" approach to defining MOIs in $B(H)$ with $H$ not necessarily separable.
Edward informed me that the latter development was of independent interest and encouraged me to publish it on its own.
This led me to write the paper \cite{nikitopoulosMOI} and allowed me to make room for the general $C^*$-algebra version of the results in this paper.
\end{ack}

\end{document}